\documentclass[10pt]{amsart}
\usepackage[T1]{fontenc}
\usepackage{geometry}
\usepackage[latin1] {inputenc}
\usepackage{amsmath}
\usepackage{amsfonts, amssymb, textcomp}
\usepackage[colorlinks=flase, linkcolor=red,urlcolor=green, citecolor=blue]{hyperref}
\usepackage{subeqnarray}

\usepackage{latexsym}
\usepackage{fancyhdr}
\usepackage{longtable}
\usepackage{amsmath, amssymb}
\usepackage{graphicx}
\usepackage{enumerate}

\numberwithin{equation}{section}

\theoremstyle{plain}
\newtheorem{proposition}{Proposition}[section]
\newtheorem{theorem}[proposition]{Theorem}
\newtheorem{lemma}[proposition]{Lemma}
\newtheorem{corollary}[proposition]{Corollary}
\newtheorem{definition}[proposition]{Definition}
\newtheorem{example}[proposition]{Example}

\newtheorem{remark}[proposition]{Remark}

\newcommand{\RR}{\mathbb{R}}
\newcommand{\CC}{\mathbb{C}}
\newcommand{\NN}{\mathbb{N}}

\newcommand{\id}{\operatorname{id}}

\let\on=\operatorname

\newenvironment{demo}[1]{\par\smallskip\noindent{\bf #1.}}{\par\smallskip}

\makeatletter
\@namedef{subjclassname@2020}{%
  \textup{2020} Mathematics Subject Classification}
\makeatother

\newsavebox{\fmbox}
\newenvironment{fmpage}[1]
 {\begin{lrbox}{\fmbox}\begin{minipage}{#1}}
 {\end{minipage}\end{lrbox}\fbox{\usebox{\fmbox}}}

\title[On the conjugate weight function]
{On the conjugate weight function and ultradifferentiable classes of entire functions}

\author[G.~Schindl]{Gerhard Schindl}

\address{G.~Schindl: Fakult\"at f\"ur Mathematik, Universit\"at Wien, Oskar-Morgenstern-Platz~1, A-1090 Wien, Austria.}
\email{gerhard.schindl@univie.ac.at}

\begin{document}

\begin{abstract}
We introduce the new notion of a conjugate weight function and provide a detailed study of this operation and its properties. Then we apply this knowledge to study classes of ultradifferentiable functions defined in terms of fast growing weight functions in the sense of Braun-Meise-Taylor and hence violating standard regularity requirements. Therefore, we transfer recent results shown by the author and D.N. Nenning from the weight sequence to the weight function framework. In order to proceed and to complete the picture we also define the conjugate associated weight matrix and investigate the relation to conjugate weight sequences via the corresponding associate weight functions. Finally, as it has already been done in the weight sequence case, we generalize results by M. Markin from the small Gevrey-setting and show how the corresponding non-standard ultradifferentiable function classes can be used to detect boundedness of normal linear operators on Hilbert spaces (associated with an evolution equation problem). On the one hand, when involving the weight matrix here the crucial information concerning regularity of the weak solutions can be expressed in terms of only one weight, namely of the given weight function. But, on the other hand, for the connection to the weighted entire setting the required conditions on the weight function are too restrictive in the general case.
\end{abstract}

\thanks{This research was funded in whole or in part by the Austrian Science Fund (FWF) 10.55776/PAT9445424}
\keywords{Weight functions and weight sequences, associated weight functions, conjugate weight function, growth and regularity properties for functions and sequences, boundedness of linear operators}
\subjclass[2020]{26A12, 26A48, 26A51, 30D15, 34G10, 46A13, 46E10, 47B02}
\date{\today}

\maketitle

\section{Introduction}
This article is the direct continuation of the author's (joint) works \cite{microclasses}, \cite{genLegendreconj}, \cite{genLegendreconjBMT}, and the aim is to transfer the main notions and results from the weight sequence case treated in \cite{microclasses} to the (ultradifferentiable) weight function setting.\vspace{6pt}

In \cite{microclasses} we have been interested in introducing and studying ultradifferentiable function classes defined in terms of non-standard (small) weight sequences $\mathbf{M}$. The goal was to emphasize the meaning of such ``exotic classes'' which are strictly smaller than the space of all real-analytic functions and, in particular, to generalize results by M. Markin from \cite{Markin01} dealing with \emph{small Gevrey sequences,} i.e. when $\mathbf{G}^s:=(p!^s)_{p\in\NN}$ with $0<s<1$, to more general (families) of weight sequences. In order to proceed we have introduced and investigated the notion of the \emph{conjugate weight sequence $\mathbf{M}^{\ast}$,} see \cite[Sect. 2.5 \& 2.6]{microclasses} and the summary in Section \ref{conjsequsection}. It has turned out that the ultradifferentiable class defined in terms of a \emph{regular small $\mathbf{M}$} can be described (as l.c.v.s.) as the weighted space of entire functions defined in terms of the \emph{associated weight function $\omega_{\mathbf{M}^{\ast}}$,} see \cite[Thm. 3.4]{microclasses}. Moreover, when considering on a Hilbert space $H$ with norm $\|\cdot\|$ the evolution equation
\begin{equation}\label{evolequ}
y'(t)=Ay(t),
\end{equation}
with \emph{$A$ being a normal (and closed) operator on $H$,} then a priori regularity of the (weak) solutions $y$ of \eqref{evolequ} expressed in terms of a family $\mathfrak{F}$ of regular small weight sequences can be used to detect boundedness of the operator $A$; see the main results \cite[Thm. 4.8 \& 4.9]{microclasses}.\vspace{6pt}

Classically, in the ultradifferentiable framework there also exists a second approach, namely dealing with a weight function $\omega:[0,+\infty)\rightarrow[0,+\infty)$ (in the sense of Braun-Meise-Taylor) instead. This setting is more recent, its modern form is based on \cite{BraunMeiseTaylor90} and it is known that both methods are in general mutually distinct; see \cite{BonetMeiseMelikhov07} and \cite{compositionpaper}. Taking into account this observation it is then natural to study the following problem:
\begin{itemize}
\item[$(Q)$] Transfer the information, definitions and notions from \cite{microclasses} to the weight function setting and then investigate the similarities and differences w.r.t. the weight sequence case.
\end{itemize}
The content of this article is to answer this question in detail. Indeed, the following two motivations are of particular interest:

\emph{First,} in order to proceed it is required to introduce the notion of the conjugate weight function $\omega^{\ast}$ and this operation seems to be new in the literature; see \eqref{conjugate}. The definition has been motivated by the Gevrey-type weights studied in Example \ref{conjugateexample}. We work in a general setting, more precisely we consider weight functions in the general sense from \cite{index} and \cite{genLegendreconj}, see Definition \ref{weightfctdef}. This gives the possibility to use the obtained information dealing with the conjugate operation $\omega\mapsto\omega^{\ast}$ in different contexts, hence the new results admit applications to other weighted settings as well and we expect that $\omega^{\ast}$ has meanings in various directions. However, in order to get that $\omega^{\ast}$ is again a weight function one crucially requires $t=o(\omega(t))$ as $t\rightarrow+\infty$, see Lemma \ref{conjugatewelldeflemma}, thus $\omega$ has to be growing fast and this condition is non-standard especially in the BMT-weight function setting.

\emph{Second,} recall that with any weight function $\omega$ satisfying basic growth properties (``BMT-weights'') one can associate a weight matrix $\mathcal{M}_{\omega}=\{\mathbf{W}^{(\ell)}: \ell>0\}$, i.e. a one-parameter family of weight sequences with parameter $\ell>0$. $\mathcal{M}_{\omega}$ can be used equivalently to describe the corresponding ultradifferentiable spaces; this is shown in \cite{compositionpaper}, \cite{dissertation} and see Section \ref{assomatrixsection} for precise definitions and a summary. Thus it seems to be possible and reasonable to consider $\mathfrak{F}=\mathcal{M}_{\omega}$ in the results from \cite[Sect. 4]{microclasses} generalizing M. Markin's approach. In this situation the crucial information concerning regularity of (weak) solutions from \eqref{evolequ} can be expressed in terms of one given weight function $\omega$ in a unified way. From this point of view the weight function setting seems to be more natural for the problem studied by M. Markin since by involving the matrix parameter the weight $\omega$ admits naturally more flexibility in order to test regularity of the solutions from \eqref{evolequ}.

Related to these comments and in order to complete the picture it is also the problem to see how the two notions of conjugation in the weight sequence and weight function setting are related to each other and how they are interacting. More precisely, the goal is to understand the connections between the functions $\omega^{\ast}_{\mathbf{M}}$ and $\omega_{\mathbf{M}^{\ast}}$ and between the matrices $\mathcal{M}^{\ast}_{\omega}:=\{(\mathbf{W}^{(\ell)})^{\ast}: \ell>0\}$ and $\mathcal{M}_{\omega^{\ast}}$.\vspace{6pt}

The paper is structured as follows: First, in Section \ref{weightgrowthsection} we gather all relevant definitions and growth conditions concerning weight functions. In Section \ref{conjuagtesection} we define the conjugate weight function $\omega^{\ast}$ and provide a detailed study of this operation: In Section \ref{conjugatedefsection} we investigate basic properties of $\omega^{\ast}$ and introduce also the double-conjugate $\omega^{\ast\ast}$; in Section \ref{conjugategrowthrelsection} we show how the conjugate operation is interacting with the known growth relations between weight functions from Section \ref{growthsection}. Then, in Section \ref{conjLegendresection} we study the interaction of the conjugate with the \emph{generalized lower and upper Legendre envelopes $\check{\star}$ and $\widehat{\star}$} from \cite{genLegendreconj}, see the summary in Section \ref{generalizedLegendresection} and the main result Proposition \ref{conjLegendreenvelopeprop}. In Section \ref{indexsection} we show how the conjugate operation is modifying the growth indices $\gamma(\cdot)$ and $\overline{\gamma}(\cdot)$ for weight functions defined and studied in \cite{index}, see Theorem \ref{indexlemma} (for precise definitions of $\gamma(\cdot)$ and $\overline{\gamma}(\cdot)$ see Section \ref{growthindexsect}). Section \ref{weightsequsection} is dedicated to the weight sequence setting; i.e. when $\omega\equiv\omega_{\mathbf{M}}$ and $\mathbf{M}\in\RR_{>0}^{\NN}$ satisfies basic growth properties. In Corollary \ref{conjwelldefcor} it is shown that both $\omega^{\ast}_{\mathbf{M}}$ and $\omega_{\mathbf{M}^{\ast}}$ are simultaneously well defined. The main result is Theorem \ref{conjequivthm} (and its Corollary \ref{conjequivthmcor1}) where the expected relation between $\omega^{\ast}_{\mathbf{M}}$ and $\omega_{\mathbf{M}^{\ast}}$ is established. Analogously as before, in Lemma \ref{growthrellemma} we investigate the interaction with growth relations (between weight sequences) and Theorem \ref{conjsequthm} provides information about the interplay with the generalized lower and upper Legendre envelopes.

Then, in Section \ref{BMTsection} we apply all this information to a BMT-weight function $\omega$ and its associated matrix $\mathcal{M}_{\omega}$. We introduce the corresponding \emph{conjugate associated matrix} $\mathcal{M}^{\ast}_{\omega}$ (see \eqref{conjumatrix}) and study its behavior in Proposition \ref{matrixadmissibleprop}: On the one hand this set is always well defined but, on the other hand, the identity $\mathcal{M}^{\ast}_{\omega}=\mathcal{M}_{\sigma}$ for a BMT-weight $\sigma$ already implies the restriction that both matrices $\mathcal{M}_{\omega}$ and $\mathcal{M}^{\ast}_{\omega}$ are \emph{constant,} i.e. all sequence are equivalent, and hence one loses flexibility w.r.t. the matrix parameter.

Similarly, this happens in the main results Theorems \ref{mainBMTthm} and \ref{mainBMTthm1} where it is shown that assumption $\overline{\gamma}(\omega)<1$ is convenient to obtain more information but also restrictive. In view of this it turns out that whenever the conjugate objects enjoy desired standard growth conditions, then this forces already for the originally given weight (matrix) the restriction mentioned before and the fact that the (ultradifferentiable) classes defined via $\omega$ (resp. equivalently $\mathcal{M}_{\omega}$) can already be described and defined via a single weight sequence; see the comments in Remark \ref{rigidrem} and Section \ref{Thm34trivialremsection}. All results mentioned so far can have more applications for other analogously defined weighted settings since we exclusively investigate weights and their growth properties.

Moreover, in Section \ref{Markinsection} we provide a concrete application and transfer the results from \cite[Sect. 4]{microclasses} to the BMT-weight function setting. On the one hand it turns out that some statements hold without any further assumption which is basically due to the fact that (each) $\mathbf{W}^{(\ell)}$ is automatically very regular; i.e. it belongs to the set \hyperlink{LCset}{$\mathcal{LC}$} (see Section \ref{sequdefsection}). But, on the other hand, the technical issue described in the previous paragraph (concerning Section \ref{BMTsection}) also yields difficulties and restrictions of the generality in this direction: By the explanations in $(d)$ in Section \ref{Markinsection} the technical result \cite[Prop. 4.5]{microclasses} has to be modified resp. generalized accordingly; see Proposition \ref{Prop45new}. This result is required for constructing a crucial specific (uniform bound) sequence $\mathbf{a}$ for the matrix $\mathcal{M}_{\omega}$ which is then used to obtain Corollary \ref{uniformboundcor} and Theorem \ref{Thm48new} even in the general (non-constant) case. (Proposition \ref{Prop45new} (also) slightly improves the results in \cite[Sect. 4]{microclasses}.) In Proposition \ref{Thm48newexample} we construct a ``non-constant example'' of a weight function to which Theorem \ref{Thm48new} can be applied; here $\omega=\omega_{\mathbf{M}}$ with $\mathbf{M}$ satisfying appropriate growth conditions and Theorem \ref{Thm48new} holds for the non-constant matrix $\mathcal{M}_{\omega_{\mathbf{M}}}$. However, in the second main result Theorem \ref{Thm49new} when combining the information from Theorem \ref{Thm48new} and the equivalence from \cite[Thm. 3.4]{microclasses}, see also Corollary \ref{conjequivthmcor} and Section \ref{Thm34trivialremsection}, then by the required assumptions on the weight in the latter result establishing the isomorphism to certain weighted spaces of entire functions, again we have to restrict to the constant case. Therefore, the constructed example also illustrates the different nature of Theorems \ref{Thm48new} and  \ref{Thm49new}.

Finally, in Appendix \ref{slowlyvaryingappendix} we study the property that $\omega_{\mathbf{M}}$ is \emph{slowly varying} (see \eqref{slowlyvarying}) in terms of the defining sequence $\mathbf{M}$; see Proposition \ref{lowvarassofctcharact}. This question is inspired by the reflexivity of a crucial growth relation from Section \ref{growthsection} and concrete applications of this technical statement are expected, especially when constructing (counter-)examples.\vspace{6pt}

\textbf{Acknowledgements.} The author of this article thanks the two anonymous referees for the careful reading and the valuable suggestions which have improved the presentation of the results.

\section{Weight functions and growth properties}\label{weightgrowthsection}

\subsection{General notation}
We write $\NN:=\{0,1,2,\dots\}$, $\NN_{>0}:=\{1,2,\dots\}$, occasionally also $\RR_{>0}:=(0,+\infty)$, and for $x\in\RR^d$ the expression $|x|$ denotes the usual Euclidean norm on $\RR^d$.

\subsection{Basic definitions and growth relations}\label{growthsection}
First, let us recall the general notion of a weight function in the sense of \cite{index} and \cite{genLegendreconj}.

\begin{definition}\label{weightfctdef}
$\omega:[0,+\infty)\rightarrow[0,+\infty)$ is called a \emph{weight function} if
\begin{itemize}
\item[$(*)$] $\omega$ is non-decreasing and

\item[$(*)$] $\lim_{t\rightarrow+\infty}\omega(t)=+\infty$.
\end{itemize}
\end{definition}
We are going to use the following notation for any (weight) function $\omega:[0,+\infty)\rightarrow[0,+\infty)$ and arbitrary $\alpha>0$: set $\omega^{\iota}(t):=\omega(\frac{1}{t})$, $t>0$, and $\omega^{1/\alpha}(t):=\omega(t^{1/\alpha})$, $t\ge 0$. Note that $\omega^{1/\alpha}$ is obtained by a power-substitution; $\omega$ is a weight function according to the previous definition if and only if $\omega^{1/\alpha}$ is so for some/any $\alpha>0$. Moreover, $(\omega^{\iota})^{1/\alpha}(t)=(\omega^{1/\alpha})^{\iota}(t)$ for all $t>0$ and finally write $\id^{1/\alpha}$ for $t\mapsto t^{1/\alpha}$. The function $\id^{1/\alpha}$ is also called the \emph{Gevrey weight} with index $\alpha>0$.

The next two conditions are natural for weight functions and appear frequently in the literature for \emph{Braun-Meise-Taylor weight functions,} see e.g. \cite{BonetMeiseMelikhov07}, and the abbreviations are motivated by the notation already used in \cite{dissertation}. We also refer to Section \ref{BMTsection} for more details:

\begin{itemize}
\item[\hypertarget{om1}{$(\omega_1)$}] $\omega(2t)=O(\omega(t))$ as $t\rightarrow+\infty$, i.e. $\exists\;L\ge 1\;\forall\;t\ge 0:\;\;\;\omega(2t)\le L(\omega(t)+1)$.
	
\item[\hypertarget{om6}{$(\omega_6)$}] $\exists\;H\ge 1\;\forall\;t\ge 0:\;2\omega(t)\le\omega(H t)+H$.
\end{itemize}

Next, let us recall some crucial growth relations between weight functions, see \cite[Sect. 2.3]{ultradifferentiablecomparison}: We write
\begin{itemize}
\item[$(*)$] $\sigma\hypertarget{ompreceq}{\preceq}\tau$ if
\begin{equation}\label{bigOrelation}
\tau(t)=O(\sigma(t))\;\text{as}\;t\rightarrow+\infty,	
\end{equation}
\item[$(*)$] $\sigma\hypertarget{omvartriangle}{\vartriangleleft}\tau$ if
\begin{equation}\label{smallOrelation}
	\tau(t)=o(\sigma(t))\;\text{as}\;t\rightarrow+\infty.	
\end{equation}
\item[$(*)$] $\sigma$ and $\tau$ are called \emph{equivalent,} written $\sigma\hypertarget{sim}{\sim}\tau$, if $\sigma\hyperlink{ompreceq}{\preceq}\tau$ and $\tau\hyperlink{ompreceq}{\preceq}\sigma$.

\item[$(*)$] \hyperlink{ompreceq}{$\preceq$} is obviously reflexive, relation \hyperlink{omvartriangle}{$\vartriangleleft$} implies \hyperlink{ompreceq}{$\preceq$} but \hyperlink{omvartriangle}{$\vartriangleleft$} is not reflexive because each weight functions tends to infinity.
\end{itemize}

Moreover, we write
\begin{itemize}
\item[$(*)$] $\sigma\hypertarget{ompreceqc}{\preceq_{\mathfrak{c}}}\tau$ if
$$\exists\;h,C>0\;\forall\;t\ge 0:\;\;\;\tau(t)\le\sigma(ht)+C,$$
\item[$(*)$] $\sigma\hypertarget{omvartrianglec}{\vartriangleleft_{\mathfrak{c}}}\tau$ if
$$\forall\;h>0\;\exists\;C_h>0\;\forall\;t\ge 0:\;\;\;\tau(t)\le\sigma(ht)+C_h.$$
\item[$(*)$] We write $\sigma\hypertarget{simc}{\sim_{\mathfrak{c}}}\tau$ if $\sigma\hyperlink{ompreceqc}{\preceq_{\mathfrak{c}}}\tau$ and $\tau\hyperlink{ompreceqc}{\preceq_{\mathfrak{c}}}\sigma$.
\item[$(*)$] \hyperlink{ompreceqc}{$\preceq_{\mathfrak{c}}$} is reflexive and \hyperlink{omvartrianglec}{$\vartriangleleft_{\mathfrak{c}}$} implies \hyperlink{ompreceqc}{$\preceq_{\mathfrak{c}}$} but \hyperlink{omvartrianglec}{$\vartriangleleft_{\mathfrak{c}}$} is in general not reflexive: If $\sigma\hyperlink{omvartrianglec}{\vartriangleleft_{\mathfrak{c}}}\sigma$ holds then $\sigma$ has to be \emph{slowly varying,} i.e.
\begin{equation}\label{slowlyvarying}
\forall\;u>0:\;\;\;\lim_{t\rightarrow+\infty}\frac{\sigma(tu)}{\sigma(t)}=1;
\end{equation}
see \cite[$(1.2.1)$]{regularvariation}, \cite[Sect. 2.2]{index}, the comments at the end of \cite[Sect. 2.3]{ultradifferentiablecomparison} and Appendix \ref{slowlyvaryingappendix} for more details.
\end{itemize}

\subsection{Growth indices $\gamma(\omega)$ and $\overline{\gamma}(\omega)$}\label{growthindexsect}
We briefly recall the definitions of the growth indices $\gamma(\omega)$ and $\overline{\gamma}(\omega)$; see \cite[Sect. 2.3 $(6)$, Sect. 2.4 $(7)$]{index} and the references in these sections.

Weight functions in the sense of Definition \ref{weightfctdef} are sufficient and convenient in order to introduce and work with the crucial indices $\gamma(\omega)$ and $\overline{\gamma}(\omega)$. Let $\omega$ be a weight function and $\gamma>0$. We say that $\omega$ has property $(P_{\omega,\gamma})$ if
\begin{equation}\label{newindex1}
\exists\;K>1:\;\;\;\limsup_{t\rightarrow+\infty}\frac{\omega(K^{\gamma}t)}{\omega(t)}<K.
\end{equation}
If $(P_{\omega,\gamma})$ holds for some $K>1$, then $(P_{\omega,\gamma'})$ is satisfied for all $\gamma'\le\gamma$ with the same $K$ since $\omega$ is non-decreasing. Moreover we can restrict to $\gamma>0$, because for $\gamma\le 0$ condition $(P_{\omega,\gamma})$ is satisfied for any weight function $\omega$, again because $\omega$ is non-decreasing and $K>1$. Then put
\begin{equation}\label{newindex2}
\gamma(\omega):=\sup\{\gamma>0:\;\;(P_{\omega,\gamma})\;\;\text{is satisfied}\},
\end{equation}
and if none condition $(P_{\omega,\gamma})$ (with $\gamma>0$) holds then set $\gamma(\omega):=0$.

Analogously, for $\gamma>0$ we say that $\omega$ has property $(\overline{P}_{\omega,\gamma})$ if
\begin{equation}\label{newindex3}
\exists\;A>1:\;\;\;\liminf_{t\rightarrow+\infty}\frac{\omega(A^{\gamma}t)}{\omega(t)}>A.
\end{equation}
If $(\overline{P}_{\omega,\gamma})$ holds for some $A>1$, then $(\overline{P}_{\omega,\gamma'})$ is satisfied for all $\gamma'\ge\gamma$ with the same $A$ since $\omega$ is non-decreasing. Moreover, we can restrict to $\gamma>0$ because for $\gamma\le 0$ condition $(\overline{P}_{\omega,\gamma})$ is never satisfied for any weight function ($\omega$ is assumed to be non-decreasing and $A>1$). Then set
\begin{equation}\label{newindex4}
\overline{\gamma}(\omega):=\inf\{\gamma>0: \;\;(\overline{P}_{\omega,\gamma})\;\;\text{is satisfied}\}.
\end{equation}
We obtain $(0\le)\gamma(\omega)\le\overline{\gamma}(\omega)$, see \cite[Sect. 2.3 \& 2.4]{index}, and by definition the following identities are valid:
\begin{equation}\label{newindex5}
\forall\;\alpha>0:\;\;\;\gamma(\omega^{1/\alpha})=\alpha\gamma(\omega),\hspace{15pt}\overline{\gamma}(\omega^{1/\alpha})=\alpha\overline{\gamma}(\omega).
\end{equation}

The indices $\gamma(\cdot)$ and $\overline{\gamma}(\cdot)$ are preserved under \hyperlink{sim}{$\sim$}; this follows from the characterizations in the main results \cite[Thm. 2.11 \& 2.16]{index}, see also \cite[Rem. 2.12]{index} and the comments after \cite[Thm. 2.16]{index}. Moreover, conditions \hyperlink{om1}{$(\omega_1)$} and \hyperlink{om6}{$(\omega_6)$} are equivalent reformulations of these growth indices as follows.

\begin{remark}\label{om1om6indexrem}
\emph{Let $\omega$ be a weight function. Then $\omega$ has \hyperlink{om1}{$(\omega_1)$} if and only if $\gamma(\omega)>0$, see \cite[Thm. 2.11, Rem. 2.12, Cor. 2.14]{index}. And by \cite[Thm. 2.16, $(7)$, Cor. 2.17]{index} condition \hyperlink{om6}{$(\omega_6)$} holds if and only if $\overline{\gamma}(\omega)<+\infty$.}
\end{remark}

\subsection{Generalized lower and upper Legendre conjugates}\label{generalizedLegendresection}
Let $\sigma,\tau$ be weight functions according to Definition \ref{weightfctdef}. We introduce the \emph{generalized lower Legendre conjugate} by
\begin{equation}\label{wedgeformula}
\sigma\check{\star}\tau(t):=\inf_{s>0}\{\sigma(s)+\tau(t/s)\},\;\;\;t\in[0,+\infty),
\end{equation}
and the \emph{generalized upper Legendre conjugate} by
\begin{equation}\label{widehatformula}
\sigma\widehat{\star}\tau(t):=\sup_{s\ge 0}\{\sigma(s)-\tau(s/t)\},\;\;\;t\in(0,+\infty).
\end{equation}
We also have $\sigma\check{\star}\tau(0)=\sigma(0)+\tau(0)$ (see \cite[Rem. 3.1]{genLegendreconj}) and set $\sigma\widehat{\star}\tau(0):=\sigma(0)-\tau(0)$ which is justified by \cite[Lemma 4.1 $(b)$]{genLegendreconj}. These operations and their effects on the growth indices $\gamma(\cdot)$, $\overline{\gamma}(\cdot)$ are studied in detail in \cite[Sect. 3 \& 4]{genLegendreconj} for (general) weight functions and in \cite[Sect. 5]{genLegendreconj} the weight sequence setting is studied; i.e. when $\sigma=\omega_{\mathbf{M}}$, $\tau=\omega_{\mathbf{N}}$ and so focusing on \emph{associated weight functions.} Concerning this setting we refer to Section \ref{assofunctionsection} for more details.\vspace{6pt}

$\check{\star}$ yields again a weight function but, in general, for $\widehat{\star}$ this is not clear; see requirements $(A)$ and $(B)$ in \cite[Sect. 4.2]{genLegendreconj}: First, in order to ensure $\sigma\widehat{\star}\tau(t)\ge 0$ for all $t$ it suffices to assume that $\tau(0)=0$, see \cite[Lemma 4.2]{genLegendreconj}. Second, in order to guarantee that $\sigma\widehat{\star}\tau$ is well defined, i.e. that $\sigma\widehat{\star}\tau(t)<+\infty$ for all $t$, it is required to have
\begin{equation}\label{equ39}
\forall\;t\in(0,+\infty)\;\exists\;D_t>0\;\forall\;s\ge 0:\;\;\;\sigma(s)-\tau(s/t)\le D_t.
\end{equation}
This condition corresponds to $(B)$ in \cite[Sect. 4.2]{genLegendreconj} resp. to \cite[$(4.4)$]{genLegendreconj} with $t_0=+\infty$. It is immediate and mentioned at the end of \cite[Sect. 2.3]{ultradifferentiablecomparison} that we have \eqref{equ39} if and only if $\tau\hyperlink{omvartrianglec}{\vartriangleleft_{\mathfrak{c}}}\sigma$; i.e. this relation gives the fact that $\sigma\widehat{\star}\tau$ is well defined on $[0,+\infty)$.


\section{The conjugate weight function}\label{conjuagtesection}
\subsection{Definition and basic properties}\label{conjugatedefsection}
Let $\omega$ be a weight function according to Definition \ref{weightfctdef}. Then introduce the \emph{conjugate function}
\begin{equation}\label{conjugate}
\omega^{\ast}(s):=\sup_{t\ge 0}\{st-\omega(t)\},\;\;\;s\in[0,+\infty),
\end{equation}
and this implies the following properties:

\begin{itemize}
\item[$(*)$] By definition $\omega^{\ast}$ is non-decreasing and convex since it is given as the supremum of convex functions (indeed the supremum of affine lines).

\item[$(*)$] One has
\begin{equation}\label{conjugatefastgrowthequ}
\forall\;s\ge 0:\;\;\;\omega^{\ast}(s)\ge\max\{s^2-\omega(s),s-\omega(1)\}.
\end{equation}
\eqref{conjugatefastgrowthequ} implies, in particular, that $\log(s)=o(\omega^{\ast}(s))$ as $s\rightarrow+\infty$ and hence $\lim_{s\rightarrow+\infty}\omega^{\ast}(s)=+\infty$.

\item[$(*)$] $\omega^{\ast}(0)=-\omega(0)$ because $\omega$ is non-decreasing and therefore, since also $\omega^{\ast}$ is non-decreasing, in order to ensure $\omega^{\ast}(s)\ge 0$ for all $s\ge 0$, one has to assume $\omega(0)=0$. In this case $\omega^{\ast}(0)=0$ holds, too.
\end{itemize}

The next result characterizes the situation when $\omega^{\ast}$ is well defined on $[0,+\infty)$.

\begin{lemma}\label{conjugatewelldeflemma}
Let $\omega$ be a weight function, then the following are equivalent:
\begin{itemize}
\item[$(i)$] $\omega^{\ast}(s)<+\infty$ for all $s\ge 0$.

\item[$(ii)$] $t=o(\omega(t))$ as $t\rightarrow+\infty$; i.e. $\omega\hyperlink{omvartriangle}{\vartriangleleft}\id$.
\end{itemize}
\end{lemma}

\demo{Proof}
$(i)\Rightarrow(ii)$ By assumption
$$\forall\;s\ge 0\;\exists\;D_s>0\;\forall\;t\ge 0:\;\;\;st-\omega(t)\le D_s,$$
thus $\limsup_{t\rightarrow+\infty}\frac{st}{\omega(t)}\le 1$ for all $s\ge 0$ since $\omega$ is a weight function. As $s\rightarrow+\infty$ this verifies $(ii)$.\vspace{6pt}

$(ii)\Rightarrow(i)$ The relation $t=o(\omega(t))$ precisely means
$$\forall\;\epsilon>0\;\exists\;D_{\epsilon}>0\;\forall\;t\ge 0:\;\;\;t\le\epsilon\omega(t)+D_{\epsilon}.$$
Fix $\epsilon>0$, then this estimate implies $\frac{t}{\epsilon}-\omega(t)\le\frac{D_{\epsilon}}{\epsilon}$ for all $t\ge 0$ and so, by definition, one has $\omega^{\ast}(\epsilon^{-1})\le\frac{D_{\epsilon}}{\epsilon}$. Since $\epsilon>0$ is arbitrary we are done.
\qed\enddemo

Let $\omega$ be a given weight function, then in view of the above observations and Lemma \ref{conjugatewelldeflemma} in order to ensure that $\omega^{\ast}$ is a weight function according to Definition \ref{weightfctdef} and satisfying $\omega^{\ast}(0)=0$, one shall assume the following properties:

\begin{itemize}
\item[$(\mathfrak{C}_1)$] $\omega(0)=0$,

\item[$(\mathfrak{C}_2)$] $t=o(\omega(t))$ as $t\rightarrow+\infty$, i.e. $\omega\hyperlink{omvartriangle}{\vartriangleleft}\id$.
\end{itemize}

Note that $(\mathfrak{C}_2)$ is obviously preserved under equivalence of weight functions and it is also stable w.r.t. \hyperlink{simc}{$\sim_{\mathfrak{c}}$}. The next result states that $\omega^{\ast}$ automatically has $(\mathfrak{C}_2)$.

\begin{lemma}\label{conjugatewelldeflemma1}
Let $\omega$ be a weight function, then $t=o(\omega^{\ast}(t))$ as $t\rightarrow+\infty$ is valid.
\end{lemma}

\demo{Proof}
The desired relation is equivalent to $\sup_{t\ge 0}\{t-\omega(t)/s\}=\frac{\omega^{\ast}(s)}{s}\rightarrow+\infty$ as $s\rightarrow+\infty$. And this holds since for any given $C\ge 1$ large there exists some $t_C\ge 0$ such that $t_C-\omega(t_C)/s\ge C\Leftrightarrow s(t_C-C)\ge\omega(t_C)$ for all sufficiently large $s\ge 0$: e.g. take $t_C:=2C$.
\qed\enddemo

\emph{Note:} Formally, for the proof of Lemma \ref{conjugatewelldeflemma1} it is not required to assume $(\mathfrak{C}_2)$ for $\omega$. If $(\mathfrak{C}_2)$ fails then $\omega^{\ast}(s)=+\infty$ for all sufficiently large $s$, see Lemma \ref{conjugatewelldeflemma} ($\omega^{\ast}$ is not a weight function) and in this case $t=o(\omega^{\ast}(t))$ is formally trivial.

\begin{example}\label{conjugateexample}
\emph{Let $\sigma:=\id^{1/\alpha}$ for $0<\alpha<1$ and so both $(\mathfrak{C}_1)$ and $(\mathfrak{C}_2)$ are satisfied. Then $\sigma^{\ast}(s)=\sup_{t\ge 0}\{st-t^{1/\alpha}\}$ and we compute this supremum by direct computations:}

\emph{For $s=0$ we have $\sigma^{\ast}(0)=0$ since $\sigma(0)=0$ and for any $s>0$ let $f_s(t):=st-t^{1/\alpha}$, $t\in[0,+\infty)$. Then $f'_s(t)=s-\frac{1}{\alpha}t^{\frac{1}{\alpha}-1}=s-\frac{1}{\alpha}t^{\frac{1-\alpha}{\alpha}}$, so $f'_s(t)=0$ if and only if $(s\alpha)^{\frac{\alpha}{1-\alpha}}=t$ and this value gives the (unique) point where the maximum of $f_s$ is attained on $[0,+\infty)$. And since $\sigma^{\ast}(s)=f_s((s\alpha)^{\frac{\alpha}{1-\alpha}})=s(s\alpha)^{\frac{\alpha}{1-\alpha}}-(s\alpha)^{\frac{1}{1-\alpha}}=s^{\frac{1}{1-\alpha}}(\alpha^{\frac{\alpha}{1-\alpha}}-\alpha^{\frac{1}{1-\alpha}})$ we see that $\sigma^{\ast}\hyperlink{sim}{\sim}\id^{\frac{1}{1-\alpha}}$.}
	
\emph{Note that by the general assumption on the parameter $\alpha$ we have that $0<1-\alpha<1$ as well and hence $\id^{1/(1-\alpha)}$ is a weight function, too. Note also that $\alpha^{\frac{\alpha}{1-\alpha}}-\alpha^{\frac{1}{1-\alpha}}>0\Leftrightarrow\alpha^{\alpha}>\alpha\Leftrightarrow 1>\alpha^{1-\alpha}$. This example implies, in particular, that up to equivalence the Gevrey-weights with indices $0<\alpha<1$ are preserved under taking the conjugate operation and that $\id^2: t\mapsto t^2$ is self-conjugate. This last fact should be compared with the information from \cite[Lemma 2.6]{microclasses} and \cite[Ex. 2.9]{genLegendreconj}.}
\end{example}

By combining Lemmas \ref{conjugatewelldeflemma} and \ref{conjugatewelldeflemma1} and motivated by the previous example, when $\omega$ is a given weight function satisfying $(\mathfrak{C}_1)$ and $(\mathfrak{C}_2)$, then it is natural to consider the \emph{double-conjugate $\omega^{\ast\ast}$} given by
$$\omega^{\ast\ast}(t):=(\omega^{\ast})^{\ast}(t)=\sup_{s\ge 0}\{ts-\omega^{\ast}(s)\},\;\;\;t\in[0,+\infty).$$
Thus $\omega^{\ast\ast}$ is again a weight function satisfying $\omega^{\ast\ast}(0)=0$. Moreover, one immediately infers that
\begin{equation}\label{doubleconjugatetrivequ}
\forall\;t\ge 0:\;\;\;\omega^{\ast\ast}(t)=\sup_{s\ge 0}\{ts-\sup_{u\ge 0}\{su-\omega(u)\}\}\underbrace{\le}_{u=t}\omega(t).
\end{equation}

\begin{remark}
\emph{Let $\omega$ be a weight function such that $(\mathfrak{C}_2)$ is violated. Then there exists some $s_0\in[0,+\infty)$ such that $\omega^{\ast}(s)=+\infty$ for all $s>s_0$. Therefore one might set in this case $\omega^{\ast\ast}(t):=\sup_{s\in[0,s_0)}\{ts-\omega^{\ast}(s)\}$ for all $t\ge 0$ and \eqref{doubleconjugatetrivequ} still holds. Finally, Lemma \ref{conjugatewelldeflemma1} becomes in this case trivial.}
\end{remark}

Next we investigate the converse estimate in \eqref{doubleconjugatetrivequ} and ask: Under which assumptions on $\omega$ one can expect that at least $\omega^{\ast\ast}\hyperlink{sim}{\sim}\omega$ holds? We prove the following result which is inspired by (the proof of) \cite[Prop. 1.6]{PetzscheVogt}:

\begin{lemma}\label{doublestarequlemma}
Let $\omega$ be a weight function satisfying $(\mathfrak{C}_1)$ and $(\mathfrak{C}_2)$. Moreover, assume that $\omega$ is \emph{convex,} then
\begin{equation}\label{doublestarequlemmaequ}
\forall\;t\ge 0:\;\;\;\omega(t)=\omega^{\ast\ast}(t).
\end{equation}
\end{lemma}

\emph{Note:} The weights treated in Example \ref{conjugateexample} are convex and thus this lemma is consistent with the computations made for these particular weights.

\demo{Proof}
First, by $(\mathfrak{C}_1)$ and $(\mathfrak{C}_2)$ the function $\omega^{\ast}$ is well defined and note that by the above comments and assumptions $\omega(0)=0=\omega^{\ast\ast}(0)$. So let now $t>0$ be given and fixed.

Let $u\ge 0$ with $u\neq t$, then by convexity one has that $u\mapsto\frac{\omega(u)-\omega(t)}{u-t}$ is non-decreasing and since $\omega$ is non-decreasing the quotient is non-negative (for any $u\neq t$). There exists $s\ge 0$ such that
$$\inf_{u>t}\frac{\omega(u)-\omega(t)}{u-t}=\omega'_{+}(t)\ge s\ge\omega'_{-}(t)=\sup_{0\le u<t}\frac{\omega(u)-\omega(t)}{u-t}.$$
Multiplying these estimates with $u-t$ we get that $s(u-t)\le\omega(u)-\omega(t)$ for all $u\ge 0$ and note that for $u=t$ this estimate is becoming trivial. Consequently $\inf_{s\ge 0}\sup_{u\ge 0}\{s(u-t)-\omega(u)\}\le-\omega(t)$, equivalently $-\inf_{s\ge 0}\sup_{u\ge 0}\{s(u-t)-\omega(u)\}\ge\omega(t)$ for all $t>0$. The left-hand side yields $-\inf_{s\ge 0}\sup_{u\ge 0}\{s(u-t)-\omega(u)\}=\sup_{s\ge 0}\{-\sup_{u\ge 0}\{s(u-t)-\omega(u)\}\}=\sup_{s\ge 0}\{st-\sup_{u\ge 0}\{su-\omega(u)\}\}=\omega^{\ast\ast}(t)$ and so the converse of \eqref{doubleconjugatetrivequ} for all $t>0$. Thus \eqref{doublestarequlemmaequ} is shown.
\qed\enddemo

\subsection{On growth relations}\label{conjugategrowthrelsection}
The aim is to study how the growth relations between weight functions are transferred under applying the conjugate operation. First, given a weight function $\omega$ with $(\mathfrak{C}_1)$ and $(\mathfrak{C}_2)$, then writing $\omega_c(t):=\omega(ct)$ for $c>0$ one immediately infers
$$\forall\;c,d>0\;\forall\;s\ge 0:\;\;\;(d\omega_c)^{\ast}(s)=\sup_{t\ge 0}\{st-d\omega(ct)\}=d\sup_{u\ge 0}\left\{\frac{s}{dc}u-\omega(u)\right\}=d\omega^{\ast}\left(\frac{s}{cd}\right).$$

\begin{lemma}\label{equivalencelemma}
Let $\sigma$ and $\tau$ be weight functions. Assume that both weights satisfy $(\mathfrak{C}_1)$ and such that $(\mathfrak{C}_2)$ is valid for $\tau$.
\begin{itemize}
\item[$(i)$] If $\sigma\hyperlink{ompreceq}{\preceq}\tau$, then $\sigma^{\ast}$ is also well defined, i.e. $\sigma$ satisfies $(\mathfrak{C}_2)$ too, and
\begin{equation}\label{equivalencelemmaequ}
\exists\;C\ge 1\;\forall\;s\ge 0:\;\;\;C\sigma^{\ast}(s)\le\tau^{\ast}(sC)+C.
\end{equation}
Conversely, \eqref{equivalencelemmaequ} implies $\sigma^{\ast\ast}\hyperlink{ompreceq}{\preceq}\tau^{\ast\ast}$ and the fact that $\sigma^{\ast}$ is well defined.

\item[$(ii)$] If $\sigma\hyperlink{omvartriangle}{\vartriangleleft}\tau$, then $\sigma^{\ast}$ is also well defined and
\begin{equation}\label{equivalencelemmaequ1}
\forall\;0<c<1\;\exists\;D_c\ge 1\;\forall\;s\ge 0:\;\;\;c\sigma^{\ast}(s)\le\tau^{\ast}(sc)+D_c.
\end{equation}
Conversely, \eqref{equivalencelemmaequ1} implies $\sigma^{\ast\ast}\hyperlink{omvartriangle}{\vartriangleleft}\tau^{\ast\ast}$ and the fact that $\sigma^{\ast}$ is well defined.
\end{itemize}

Moreover:
\begin{itemize}
\item[$(a)$] $\sigma\hyperlink{ompreceq}{\preceq}\tau$ implies the mixed variant of \hyperlink{om6}{$(\omega_6)$} between $\sigma^{\ast}$ and $\tau^{\ast}$. If either $\sigma^{\ast}$ and $\tau^{\ast}$ satisfies in addition \hyperlink{om1}{$(\omega_1)$}, then $\sigma^{\ast}$ and $\tau^{\ast}$ are related by the mixed variant of \hyperlink{om1}{$(\omega_1)$}.

\item[$(b)$] $\sigma\hyperlink{omvartriangle}{\vartriangleleft}\tau$ implies both the mixed variant of \hyperlink{om6}{$(\omega_6)$} and \hyperlink{om1}{$(\omega_1)$} between $\sigma^{\ast}$ and $\tau^{\ast}$.
\end{itemize}
\end{lemma}

\demo{Proof}
$(i)$ By assumption $\tau(t)\le C\sigma(t)+C$ for some $C\ge 1$ and for all $t\ge 0$ and $\tau^{\ast}$ is well defined. Thus, for any $s\ge 0$:
\begin{align*}
\tau^{\ast}(s)&=\sup_{t\ge 0}\{st-\tau(t)\}\ge\sup_{t\ge 0}\{st-C\sigma(t)\}-C=C\sup_{t\ge 0}\{(s/C)t-\sigma(t)\}-C
\\&
=C\sigma^{\ast}(s/C)-C.
\end{align*}
This estimate implies, in particular, that $\sigma^{\ast}$ is well defined; see also Lemma \ref{conjugatewelldeflemma} applied to $\sigma$. Moreover, \eqref{equivalencelemmaequ} also implies that $\sigma^{\ast}$ is well defined since $\tau^{\ast}$ is so and for any $s\ge 0$ one has:
\begin{align*}
\tau^{\ast\ast}(s)&=\sup_{t\ge 0}\{st-\tau^{\ast}(t)\}\le\sup_{t\ge 0}\{st-C\sigma^{\ast}(t/C)\}+C=C\sup_{t\ge 0}\{\frac{st}{C}-\sigma^{\ast}(t/C)\}+C
\\&
=C\sup_{u\ge 0}\{su-\sigma^{\ast}(u)\}+C=C\sigma^{\ast\ast}(s)+C.
\end{align*}

$(ii)$ By assumption for any $1>c>0$ (small) there exists some $D_c\ge 1$ such that $\tau(t)\le c\sigma(t)+D_c$ for all $t\ge 0$. Consequently, the conclusions follow then by the same computations as in $(i)$.\vspace{6pt}

$(a)$ W.l.o.g. one can assume that $C\ge 2$ in relation $\sigma\hyperlink{ompreceq}{\preceq}\tau$ and hence \eqref{equivalencelemmaequ} implies $2\sigma^{\ast}(s)\le\tau^{\ast}(sC)+C$ for all $s\ge 0$ which is the mixed version of \hyperlink{om6}{$(\omega_6)$}.\vspace{6pt}

If $\tau^{\ast}$ satisfies in addition \hyperlink{om1}{$(\omega_1)$}, then choose $n\in\NN_{>0}$ such that $2^{n-1}\ge C$ with $C$ denoting the constant appearing in \eqref{equivalencelemmaequ}. Iteration of \hyperlink{om1}{$(\omega_1)$} implies
$C\sigma^{\ast}(2s)\le\tau^{\ast}(2Cs)+C\le\tau^{\ast}(2^ns)+C\le L\tau^{\ast}(s)+L+C$ for some $L\ge 1$ and all $s\ge 0$. Similarly, if $\sigma^{\ast}$ satisfies in addition \hyperlink{om1}{$(\omega_1)$}, then $C\sigma^{\ast}(2s)\le CL\sigma^{\ast}(C^{-1}s)+CL\le L\tau^{\ast}(s)+CL+CL$ for some $L\ge 1$ and all $s\ge 0$. Both derived estimates yield the conclusion.\vspace{6pt}

$(b)$ Since $\sigma\hyperlink{omvartriangle}{\vartriangleleft}\tau$ is stronger than $\sigma\hyperlink{ompreceq}{\preceq}\tau$, the mixed variant of \hyperlink{om6}{$(\omega_6)$} follows by $(a)$ and the mixed variant of \hyperlink{om1}{$(\omega_1)$} holds by choosing $c:=\frac{1}{2}$ in \eqref{equivalencelemmaequ1}.
\qed\enddemo

\emph{Note:} In the proof of $(ii)(a)$ \emph{any choice} $n\in\NN_{>0}$ such that $2^{n-1}\ge C$ holds is sufficient to conclude and, indeed, one may select the \emph{minimal} $n$. This comment applies to any argument involving iteration techniques; in particular we use this in several proofs in this article.

\begin{corollary}\label{equivalencelemmacor}
Let $\omega$ be a weight function satisfying $(\mathfrak{C}_1)$ and $(\mathfrak{C}_2)$. Moreover, assume that $\omega$ is equivalent to a convex weight function $\sigma$ satisfying $(\mathfrak{C}_1)$. Then $\omega$ and $\omega^{\ast\ast}$ are equivalent; more precisely:
\begin{equation}\label{equivalencelemmacorequ}
\omega\hyperlink{sim}{\sim}\sigma=\sigma^{\ast\ast}\hyperlink{sim}{\sim}\omega^{\ast\ast}.
\end{equation}
\end{corollary}

\demo{Proof}
First, recall that $(\mathfrak{C}_2)$ is preserved under equivalence and hence $\sigma$ shares this property, too. Then $(i)$ in Lemma \ref{equivalencelemma} gives that $\omega^{\ast\ast}$ and $\sigma^{\ast\ast}$ are equivalent and by Lemma \ref{doublestarequlemma} we infer that $\sigma(t)=\sigma^{\ast\ast}(t)$ for any $t\ge 0$. Combining everything yields the desired relations stated in \eqref{equivalencelemmacorequ}.
\qed\enddemo

Next we treat \hyperlink{ompreceq}{$\preceq_{\mathfrak{c}}$}, \hyperlink{omvartriangle}{$\vartriangleleft_{\mathfrak{c}}$} and verify that these relations are reversed when taking the conjugate. Thus these relations can be considered to be more natural when dealing with conjugate weight functions.

\begin{lemma}\label{equivalencelemma1}
Let $\sigma$ and $\tau$ be weight functions. Assume that both weights satisfy $(\mathfrak{C}_1)$ and such that $(\mathfrak{C}_2)$ is valid for $\tau$.
\begin{itemize}
\item[$(i)$] If $\sigma\hyperlink{ompreceqc}{\preceq_{\mathfrak{c}}}\tau$, then
\begin{equation}\label{equivalencelemma1equ}
\exists\;h,C\ge 1\;\forall\;s\ge 0:\;\;\;\sigma^{\ast}(s)\le\tau^{\ast}(sh)+C;
\end{equation}
i.e. $\tau^{\ast}\hyperlink{ompreceqc}{\preceq_{\mathfrak{c}}}\sigma^{\ast}$. Moreover, $\sigma^{\ast}$ is also well defined, i.e. $\sigma$ satisfies $(\mathfrak{C}_2)$, too.

Conversely, \eqref{equivalencelemma1equ} implies $\sigma^{\ast\ast}\hyperlink{ompreceqc}{\preceq_{\mathfrak{c}}}\tau^{\ast\ast}$ and that $\sigma^{\ast}$ is also well defined.

\item[$(ii)$] If $\sigma\hyperlink{omvartrianglec}{\vartriangleleft_{\mathfrak{c}}}\tau$, then
\begin{equation}\label{equivalencelemma1equ1}
\forall\;0<h<1\;\exists\;C_h\ge 1\;\forall\;s\ge 0:\;\;\;\sigma^{\ast}(s)\le\tau^{\ast}(sh)+C_h;
\end{equation}
i.e. $\tau^{\ast}\hyperlink{omvartrianglec}{\vartriangleleft_{\mathfrak{c}}}\sigma^{\ast}$. Moreover, $\sigma^{\ast}$ is also well defined. Conversely, \eqref{equivalencelemma1equ1} implies that $\sigma^{\ast}$ is well defined and $\sigma^{\ast\ast}\hyperlink{omvartrianglec}{\vartriangleleft_{\mathfrak{c}}}\tau^{\ast\ast}$.
\end{itemize}
In addition we get:
\begin{itemize}
\item[$(a)$] If $\sigma\hyperlink{ompreceqc}{\preceq_{\mathfrak{c}}}\tau$ and either $\sigma^{\ast}$ and $\tau^{\ast}$ satisfies in addition \hyperlink{om1}{$(\omega_1)$}, then $\sigma^{\ast}$ and $\tau^{\ast}$ are related by the mixed variant of \hyperlink{om1}{$(\omega_1)$}.

\item[$(b)$] $\sigma\hyperlink{omvartrianglec}{\vartriangleleft_{\mathfrak{c}}}\tau$ implies the mixed variant of \hyperlink{om1}{$(\omega_1)$} between $\sigma^{\ast}$ and $\tau^{\ast}$.
\end{itemize}
\end{lemma}

\demo{Proof}
$(i)$ By assumption $\tau(t)\le\sigma(ht)+C$ for some $C,h\ge 1$ and for all $t\ge 0$ and $\tau^{\ast}$ is well defined. Then, for any $s\ge 0$:
\begin{align*}
\tau^{\ast}(s)&=\sup_{t\ge 0}\{st-\tau(t)\}\ge\sup_{t\ge 0}\{st-\sigma(ht)\}-C=\sup_{u\ge 0}\{(s/h)u-\sigma(u)\}-C
\\&
=\sigma^{\ast}(s/h)-C.
\end{align*}
This estimate implies that $\sigma^{\ast}$ is well defined. Conversely, we get for any $s\ge 0$:
\begin{align*}
\tau^{\ast\ast}(s)&=\sup_{t\ge 0}\{st-\tau^{\ast}(t)\}\le\sup_{t\ge 0}\{st-\sigma^{\ast}(t/h)\}+C=\sup_{u\ge 0}\{(sh)u-\sigma^{\ast}(u)\}+C
\\&
=\sigma^{\ast\ast}(sh)+C.
\end{align*}

$(ii)$ By assumption even for any $1>h>0$ (small) there exists some $C_h\ge 1$ such that $\tau(t)\le\sigma(ht)+C_h$ for all $t\ge 0$. Consequently, the conclusions follow by the same computations as in $(i)$.\vspace{6pt}

$(a)$ Note that in general $h>2$ and hence the conclusion does not follow automatically. However, by iterating \hyperlink{om1}{$(\omega_1)$} for the particular weight one is able to conclude analogously as in the proof of Lemma \ref{equivalencelemma}.\vspace{6pt}

$(b)$ We choose $h:=\frac{1}{2}$.
\qed\enddemo

\subsection{On the generalized Legendre envelopes}\label{conjLegendresection}
We study the interaction of the conjugate operation and the generalized Legendre envelope operators $\check{\star}$ and $\widehat{\star}$.

\begin{proposition}\label{conjLegendreenvelopeprop}
Let $\sigma$ and $\tau$ be given weight functions such that $\sigma(0)=0=\tau(0)$, i.e. both weights satisfy $(\mathfrak{C}_1)$. Moreover, assume that $\sigma\check{\star}\tau$ satisfies $(\mathfrak{C}_2)$, then
\begin{equation}\label{conjLegendreenvelopepropequ}
\forall\;t\ge 0:\;\;\;\tau^{\ast}\widehat{\star}\sigma(t)=(\sigma\check{\star}\tau)^{\ast}(t)=\sigma^{\ast}\widehat{\star}\tau(t),
\end{equation}
and
\begin{equation}\label{conjLegendreenvelopepropequ2}
\forall\;t\ge 0:\;\;\;\sigma^{\ast}\widehat{\star}\tau(t)=((\id\widehat{\star}\tau)\widehat{\star}\sigma)(t).
\end{equation}
\end{proposition}

\demo{Proof}
$\sigma\check{\star}\tau$ is a weight function such that $\sigma\check{\star}\tau(0)=\sigma(0)+\tau(0)=0$, see \cite[Rem. 3.1, Lemma 3.2 $(a)\&(b)$]{genLegendreconj}, and by assumption $(\sigma\check{\star}\tau)^{\ast}$ is well defined. Note that $(\mathfrak{C}_2)$ for $\sigma\check{\star}\tau$ precisely means $t=o(\sigma\check{\star}\tau(t))$ as $t\rightarrow+\infty$ and, in view of \cite[Lemma 3.2 $(a)$]{genLegendreconj}, this property implies $(\mathfrak{C}_2)$ for $\sigma$ and $\tau$. Let $t>0$ be given, then:
\begin{align*}
(\sigma\check{\star}\tau)^{\ast}(t)&=\sup_{s\ge 0}\{st-\sigma\check{\star}\tau(s)\}=\sup_{s\ge 0}\{st-\inf_{u>0}\{\sigma(u)+\tau(s/u)\}\}
\\&
=\sup_{s\ge 0}\{st+\sup_{u>0}\{-\sigma(u)-\tau(s/u)\}\}=\sup_{s\ge 0}\sup_{u>0}\{st-\sigma(u)-\tau(s/u)\}
\\&
=\sup_{u>0}\sup_{s\ge 0}\{st-\sigma(u)-\tau(s/u)\}=\sup_{u>0}\{-\sigma(u)+\sup_{s\ge 0}\{st-\tau(s/u)\}\}
\\&
=\sup_{u>0}\{-\sigma(u)+\sup_{v\ge 0}\{uvt-\tau(v)\}\}=\sup_{u>0}\{-\sigma(u)+\tau^{\ast}(ut)\}
\\&
=\sup_{w>0}\{\tau^{\ast}(w)-\sigma(w/t)\}=\sup_{w\ge 0}\{\tau^{\ast}(w)-\sigma(w/t)\}=\tau^{\ast}\widehat{\star}\sigma(t).
\end{align*}
Note that for $t=0$ we have $(\sigma\check{\star}\tau)^{\ast}(0)=-\sigma\check{\star}\tau(0)=0$, whereas by \cite[Lemma 4.1 $(b)$]{genLegendreconj} and $(\mathfrak{C}_1)$ it holds that $\tau^{\ast}\widehat{\star}\sigma(0)=-\tau(0)-\sigma(0)=0$. Summarizing, the first part in \eqref{conjLegendreenvelopepropequ} is verified.

And since $\check{\star}$ is commutative, see \cite[Lemma 3.2 $(c)$]{genLegendreconj}, analogously one infers the second one, too. Finally, by this computation$(\mathfrak{C}_2)$ for $\sigma\check{\star}\tau$ implies that both $\tau^{\ast}\widehat{\star}\sigma$ and $\sigma^{\ast}\widehat{\star}\tau$ are well defined.\vspace{6pt}

Concerning \eqref{conjLegendreenvelopepropequ2} let $t>0$, then
\begin{align*}
(\sigma^{\ast}\widehat{\star}\tau)(t)&=\sup_{s\ge 0}\{\sigma^{\ast}(s)-\tau(s/t)\}=\sup_{s\ge 0}\{\sup_{u\ge 0}\{su-\sigma(u)\}-\tau(s/t)\}
\\&
=\sup_{u\ge 0}\{\sup_{s\ge 0}\{su-\sigma(u)\}-\tau(s/t)\}=\sup_{u\ge 0}\{\sup_{s\ge 0}\{su-\tau(s/t)\}-\sigma(u)\}
\\&
=\sup_{u>0}\{\sup_{v\ge 0}\{v-\tau(v/(ut))\}-\sigma(u)\}=\sup_{u>0}\{\id\widehat{\star}\tau(ut)-\sigma(u)\}
\\&
=\sup_{w>0}\{\id\widehat{\star}\tau(w)-\sigma(w/t)\}=\sup_{w\ge 0}\{\id\widehat{\star}\tau(w)-\sigma(w/t)\}=((\id\widehat{\star}\tau)\widehat{\star}\sigma)(t).
\end{align*}
For $t=0$, again by $(\mathfrak{C}_1)$ and \cite[Lemma 4.1 $(b)$]{genLegendreconj} we infer $(\sigma^{\ast}\widehat{\star}\tau)(0)=-\sigma(0)-\tau(0)=0$, $((\id\widehat{\star}\tau)\widehat{\star}\sigma)(0)=(\id(0)-\tau(0))-\sigma(0)=0$, and hence equality for all $t\ge 0$.
\qed\enddemo

\subsection{The conjugate operation and growth indices}\label{indexsection}
We study now the effect of taking the conjugate operation defined in \eqref{conjugate} on the growth indices $\gamma(\omega)$ and $\overline{\gamma}(\omega)$ defined in Section \ref{growthindexsect}. First, note that $\sigma\hyperlink{ompreceq}{\preceq}\sigma$ (reflexivity) is clear and one can take in the proof of $(i)$ in Lemma \ref{equivalencelemma} any constant $C\ge 1$ and $\sigma=\tau$. Thus, by applying $(a)$ in this result to $\sigma=\tau$ and in view of Remark \ref{om1om6indexrem} we immediately get the following result:

\begin{lemma}\label{indexlemma0}
Let $\sigma$ be a weight function satisfying $(\mathfrak{C}_1)$ and $(\mathfrak{C}_2)$. Then $\sigma^{\ast}$ satisfies \hyperlink{om6}{$(\omega_6)$}, i.e. $\overline{\gamma}(\sigma^{\ast})<+\infty$.
\end{lemma}

Moreover, the second part of $(a)$ in Lemma \ref{equivalencelemma} becomes trivial, however $(b)$ cannot be applied to $\sigma=\tau$ since relation \hyperlink{omvartriangle}{$\vartriangleleft$} is not reflexive.

On the other hand we get that assumption $(\mathfrak{C}_2)$ also implies a certain restriction of the value $\gamma(\sigma)$:

\begin{lemma}\label{indexlemma1}
Let $\sigma$ be a weight function. Then $(\mathfrak{C}_2)$ implies $\gamma(\sigma)\le 1$; i.e. $\gamma(\sigma)\in[0,1]$.
\end{lemma}

\demo{Proof}
By \cite[Thm. 2.11 \& Rem. 2.15]{index} it follows that $\gamma(\sigma)>1$ implies $\omega(t)=O(t^{\alpha})$ as $t\rightarrow+\infty$ for some $\alpha\in(0,1)$ which obviously contradicts $(\mathfrak{C}_2)$.
\qed\enddemo

The next main result relates the growth indices of $\sigma$ and $\sigma^{\ast}$ in a precise manner and should be compared with \cite[$(A.1)\&(A.2)$]{microclasses}.

\begin{theorem}\label{indexlemma}
Let $\sigma$ be a weight function satisfying $(\mathfrak{C}_1)$ and $(\mathfrak{C}_2)$.
\begin{itemize}
\item[$(i)$] It holds that $\max\{0,1-\overline{\gamma}(\sigma)\}\le\gamma(\sigma^{\ast})$.

\item[$(ii)$] It holds that $\overline{\gamma}(\sigma^{\ast})\le 1-\gamma(\sigma)$.

\item[$(iii)$] Assume that in addition $\sigma^{\ast\ast}\hyperlink{sim}{\sim}\sigma$ holds. Then the following is valid:
\begin{itemize}
\item[$(a)$] It holds that  $\gamma(\sigma)=1-\overline{\gamma}(\sigma^{\ast})$.

\item[$(b)$] It holds that $1-\overline{\gamma}(\sigma)=\gamma(\sigma^{\ast})$.

\item[$(c)$] It holds that $\overline{\gamma}(\sigma)=1-\gamma(\sigma^{\ast})$.

\item[$(d)$] It holds that $\gamma(\sigma)=1-\overline{\gamma}(\sigma^{\ast})$.
\end{itemize}
In this situation all indices belong to the closed unit interval $[0,1]$.
\end{itemize}
\end{theorem}

\emph{Note:} If $\overline{\gamma}(\omega)<1$, then the estimate in $(i)$ in Theorem \ref{indexlemma} implies, in particular, that $\gamma(\sigma^{\ast})>0$ and so \hyperlink{om1}{$(\omega_1)$} for $\sigma^{\ast}$. If $\overline{\gamma}(\sigma)\ge 1$, then this property for $\sigma^{\ast}$ is not clear.

In view of Lemma \ref{indexlemma1} in $(ii)$ one has $\gamma(\sigma)\in[0,1]$ and hence the estimate is natural and consistent with the growth limitation obtained in Lemma \ref{indexlemma1}: $\gamma(\sigma)>1$ would imply $\overline{\gamma}(\sigma^{\ast})<0$ but this is impossible for any weight function.

\demo{Proof}
First, note that by assumption on $\sigma$ the conjugate $\sigma^{\ast}$ is also a weight function satisfying $(\mathfrak{C}_1)$ and $(\mathfrak{C}_2)$ and $\sigma^{\ast\ast}\hyperlink{sim}{\sim}\sigma$ implies $\gamma(\sigma^{\ast\ast})=\gamma(\sigma)$, $\overline{\gamma}(\sigma^{\ast\ast})=\overline{\gamma}(\sigma)$.\vspace{6pt}

$(i)$ Since $\sigma^{\ast}$ is a weight function we get $\gamma(\sigma^{\ast})\ge 0$ and hence the case $\overline{\gamma}(\sigma)\ge 1$ is trivial. So let now $(0\le)\overline{\gamma}(\sigma)<1$ and let $\gamma>0$ be such that $\gamma<1-\overline{\gamma}(\sigma)(\le 1)$ and such a choice is possible by assumption on the index $\overline{\gamma}(\sigma)$. Since $1-\gamma>\overline{\gamma}(\sigma)$ by definition of this index
$$\exists\;A>1\;\exists\;\epsilon_0\in(0,1)\;\exists\;H\ge 0\;\forall\;\epsilon\in(0,\epsilon_0)\;\forall\;t\ge 0:\;\;\;\sigma(A^{1-\gamma}t)\ge A^{1+\epsilon}\sigma(t)-H,$$
and we estimate as follows for all $s\ge 0$:
\begin{align*}
\sigma^{\ast}(A^{\gamma-1}s)&=\sup_{t\ge 0}\{(A^{\gamma-1}s)t-\sigma(t)\}\underbrace{=}_{u=A^{\gamma-1}t}\sup_{u\ge 0}\{su-\sigma(A^{1-\gamma}u)\}
\\&
\le\sup_{u\ge 0}\{su-A^{1+\epsilon}\sigma(u)\}+H=A^{1+\epsilon}\sup_{u\ge 0}\{(A^{-1-\epsilon}s)u-\sigma(u)\}+H
\\&
=A^{1+\epsilon}\sigma^{\ast}(A^{-1-\epsilon}s)+H.
\end{align*}
Summarizing, this estimate gives $\sigma^{\ast}(A^{\gamma+\epsilon}s)\le A^{1+\epsilon}\sigma^{\ast}(s)+H$ for all $s\ge 0$. Now fix $\epsilon'\in(0,\epsilon_0)$ and let $\epsilon>0$ be such that $\epsilon<\gamma\epsilon'<\epsilon'<\epsilon_0$. Then set $K:=A^{1+\epsilon'}(>A^{1+\epsilon})$, $\gamma':=\frac{\gamma+\epsilon}{1+\epsilon'}$ and so $\gamma'<\gamma\Leftrightarrow\epsilon<\gamma\epsilon'$ holds. Therefore, the verified estimate implies $(P_{\sigma^{\ast},\gamma'})$ with this choice for $K$ and so by definition $\gamma(\sigma^{\ast})\ge\gamma'$. (Indeed, in view of \cite[Thm. 2.11]{index} we even have a strict inequality here.) Then let $\epsilon'\rightarrow 0$, consequently $\epsilon\rightarrow0$ and hence $\gamma(\sigma^{\ast})\ge\gamma$ holds. Since $0<\gamma<1-\overline{\gamma}(\sigma)$ was chosen arbitrarily we have shown that $1-\overline{\gamma}(\sigma)\le\gamma(\sigma^{\ast})$.\vspace{6pt}

$(ii)$ First let $\gamma(\sigma)\in(0,1)$ and let $0<\gamma<1$ be such that $\gamma(\sigma)>1-\gamma(>0)$, i.e. $\gamma>1-\gamma(\sigma)$, and note that such values $\gamma$ exist by assumption on $\gamma(\sigma)$. Then by definition we get
$$\exists\;K>1\;\exists\;\epsilon_0\in(0,1)\;\exists\;H\ge 0\;\forall\;\epsilon\in(0,\epsilon_0)\;\forall\;t\ge 0:\;\;\;\sigma(K^{1-\gamma}t)\le K^{1-\epsilon}\sigma(t)+H.$$
Hence for all $s\ge 0$:
\begin{align*}
\sigma^{\ast}(K^{\gamma-1}s)&=\sup_{t\ge 0}\{(K^{\gamma-1}s)t-\sigma(t)\}\underbrace{=}_{u=K^{\gamma-1}t}\sup_{u\ge 0}\{su-\sigma(K^{1-\gamma}u)\}
\\&
\ge\sup_{u\ge 0}\{su-K^{1-\epsilon}\sigma(u)\}-H=K^{1-\epsilon}\sup_{u\ge 0}\{(K^{\epsilon-1}s)u-\sigma(u)\}-H
\\&
=K^{1-\epsilon}\sigma^{\ast}(K^{\epsilon-1}s)-H.
\end{align*}
This proves $\sigma^{\ast}(K^{\gamma-\epsilon}s)\ge K^{1-\epsilon}\sigma^{\ast}(s)-H$ for all $s\ge 0$.

Let $\epsilon'\in(0,\epsilon_0)$ be fixed, then choose $\epsilon>0$ such that $\epsilon<\gamma\epsilon'<\epsilon'<\epsilon_0$ and put $A:=K^{1-\epsilon'}(<K^{1-\epsilon})$, $\gamma':=\frac{\gamma-\epsilon}{1-\epsilon'}$. Then $A>1$ and $\gamma'>\gamma\Leftrightarrow\epsilon'\gamma>\epsilon$ and the above estimate verifies that $(\overline{P}_{\sigma^{\ast},\gamma'})$ is valid for this choice of $A$. Then let $\epsilon'\rightarrow 0$, hence $\epsilon\rightarrow 0$ too and so $\gamma'\rightarrow\gamma$ which verifies $\overline{\gamma}(\sigma^{\ast})\le\gamma$. Finally, we let $\gamma\rightarrow 1-\gamma(\sigma)$ and therefore $\overline{\gamma}(\sigma^{\ast})\le 1-\gamma(\sigma)$ is shown.\vspace{6pt}

If $\gamma(\sigma)=1$, then the same proof gives $\overline{\gamma}(\sigma^{\ast})=0$: Indeed, here in the arguments we can consider \emph{any} $\gamma>0$ and we restrict again to $\gamma\in(0,1)$ and then let $\gamma\rightarrow 1-\gamma(\sigma)=0$ which gives $\overline{\gamma}(\sigma^{\ast})=1-\gamma(\sigma)=0$. Finally, if $\gamma(\sigma)=0$, then follow the proof above with $\gamma>1$ and choose $\epsilon,\epsilon'\in(0,\epsilon_0)$ such that $\epsilon<\epsilon'<\epsilon'\gamma<\epsilon_0$. Again, as $\epsilon',\epsilon\rightarrow 0$ we get $\overline{\gamma}(\sigma^{\ast})\le\gamma$ and then let $\gamma\rightarrow 1(=1-\gamma(\sigma))$.\vspace{6pt}

$(iii)(a)$ First, $(ii)$ implies $\overline{\gamma}(\sigma^{\ast})\le 1-\gamma(\sigma)\le 1$. Then apply $(i)$ to $\sigma^{\ast}$ which gives $1-\overline{\gamma}(\sigma^{\ast})\le\gamma(\sigma^{\ast\ast})=\gamma(\sigma)$. For these arguments note that $\sigma^{\ast}$ is a weight function satisfying $(\mathfrak{C}_1)$ and $(\mathfrak{C}_2)$ (recall Lemma \ref{conjugatewelldeflemma1}). Summarizing, $1-\overline{\gamma}(\sigma^{\ast})=\gamma(\sigma)$ holds as desired.

$(iii)(b)$ We apply $(ii)$ to $\sigma^{\ast}$ and which gives $(0\le)\overline{\gamma}(\sigma)=\overline{\gamma}(\sigma^{\ast\ast})\le 1-\gamma(\sigma^{\ast})(\le 1)$. Second, using this $(i)$ implies $0\le 1-\overline{\gamma}(\sigma)\le\gamma(\sigma^{\ast})$ and combining everything yields $1-\overline{\gamma}(\sigma)=\gamma(\sigma^{\ast})$.

$(iii)(c)$ Replace $\sigma$ by $\sigma^{\ast}$ and $\sigma^{\ast}$ by $\sigma^{\ast\ast}$ in the proof of $(iii)(a)$.

$(iii)(d)$ Replace $\sigma$ by $\sigma^{\ast}$ and $\sigma^{\ast}$ by $\sigma^{\ast\ast}$ in the proof of $(iii)(b)$.
\qed\enddemo

\begin{remark}
\emph{Consider the weight(s) from Example \ref{conjugateexample}: First, note that $\gamma(\id^{1/\alpha})=\alpha=\overline{\gamma}(\id^{1/\alpha})$ for any $\alpha>0$ and thus Theorem \ref{indexlemma} holds for resp. can be applied to all $\alpha\in(0,1)$. However, as verified in this Example, since for such indices $\alpha$ one infers $(\id^{1/\alpha})^{\ast}\hyperlink{sim}{\sim}\id^{1/(1-\alpha)}$ it follows that the (limiting) case $\alpha=1$ cannot be considered in Theorem \ref{indexlemma}.}
\end{remark}

\section{The weight sequence setting}\label{weightsequsection}
In this section we focus on the case $\omega=\omega_{\mathbf{M}}$, i.e. when the weight function is expressed in terms of a given weight sequence $\mathbf{M}$. $\omega_{\mathbf{M}}$ denotes the corresponding \emph{associated weight function} and due to this additional information we can consider $\omega_{\mathbf{M}}^{\ast}$ and $\omega_{\mathbf{M}^{\ast}}$, $\mathbf{M}^{\ast}$ denoting the \emph{conjugate sequence} introduced and used in \cite{microclasses}. Thus the goal is to establish a precise connection between both notions: the conjugate in the weight function and in the weight sequence framework.\vspace{6pt}

Frequently, we use the following estimates (with the convention $0^0:=1$):
\begin{equation}\label{Stirlinglike}
\forall\;p\in\NN:\;\;\;p!\le p^p\le e^pp!.
\end{equation}
The first estimate is clear and the second one follows, as recognized by one of the anonymous referees in the report, by the power series expansion of the exponential function since $e^p=\sum_{k=0}^{+\infty}\frac{p^k}{k!}\ge\frac{p^p}{p!}$ for all $p\in\NN$. Indeed, formally \emph{Stirling's formula} is not required to prove \eqref{Stirlinglike} and note that these useful estimates appear regularly in this context in the literature and also in some of the author's works.

\subsection{Weight sequences}\label{sequdefsection}
Given a sequence $\mathbf{M}=(M_p)_p\in\RR_{>0}^{\NN}$ we also use the notation $\mu=(\mu_p)_p$ with $\mu_p:=\frac{M_p}{M_{p-1}}$, $p\ge 1$, $\mu_0:=1$. Moreover, we set $m_p:=\frac{M_p}{p!}$ and obtain the sequence $\mathbf{m}=(m_p)_{p\in\NN}$. Analogously these conventions are used for all other sequences under consideration. $\mathbf{M}$ is called \emph{normalized,} if $1=M_0\le M_1$ is valid. For the \emph{pointwise product sequence} we simply write $\mathbf{M}\cdot\mathbf{N}$ and $\frac{\mathbf{M}}{\mathbf{N}}$ for the \emph{pointwise quotient sequence;} i.e. these sequences are given by $M_p\cdot N_p$ and $\frac{M_p}{N_p}$ respectively.\vspace{6pt}

$\mathbf{M}$ is called \emph{log-convex} if
$$\forall\;p\in\NN_{>0}:\;M_p^2\le M_{p-1} M_{p+1},$$
equivalently if $\mu$ is non-decreasing. $\mathbf{M}$ is called \emph{strong log-convex} if even $\mathbf{m}$ is log-convex which is equivalent to $p\mapsto\frac{\mu_p}{p}$ being non-decreasing. If $\mathbf{M}$ is log-convex and normalized, then $\mu_p\ge 1$ for all $p\in\NN$, both mappings $p\mapsto M_p$ and $p\mapsto(M_p)^{1/p}$ are non-decreasing, $(M_p)^{1/p}\le\mu_p$ for all $p\in\NN_{>0}$ and finally $M_pM_q\le M_{p+q}$ for all $p,q\in\NN$; see e.g. \cite[Lemmas 2.0.4 \& 2.0.6]{diploma}.

$\mathbf{M}$ satisfies \emph{moderate growth,} abbreviated in the following by \hypertarget{mg}{$(\on{mg})$}, if
$$\exists\;C\ge 1\;\forall\;p,q\in\NN:\;M_{p+q}\le C^{p+q+1} M_p M_q.$$
In view of \cite{Komatsu73} this condition is also known under \emph{(M.2) or stability under ultradifferential operators.} In case $M_0=1$, then it suffices to consider $C^{p+q}$ in the estimate.\vspace{6pt}

We introduce the value
$$\mathbf{M}_{\iota}:=\liminf_{p\rightarrow+\infty}(M_p/M_0)^{1/p}=\liminf_{p\rightarrow+\infty}(M_p)^{1/p},$$
and recall the fact that for any log-convex sequence one has $\mathbf{M}_{\iota}=\lim_{p\rightarrow+\infty}(M_p)^{1/p}=\lim_{p\rightarrow+\infty}\mu_p\in(0,+\infty]$, see e.g. \cite[Sect. 2.2, Lemma 2.1]{regularnew}, and we comment on a minor technical issue in this result:

\begin{remark}\label{Lemma21rem}
\emph{By inspecting the proof of \cite[Lemma 2.1]{regularnew} it turns out that one shall assume there $\lim_{p\rightarrow+\infty}\mu_p=C\in(0,+\infty)$ and note that the case $C=0$ is not excluded by the general assumption $\mathbf{M}\in\RR_{>0}^{\NN}$. However, when $\lim_{p\rightarrow+\infty}\mu_p=0$ then the result is still valid by following the arguments and only considering the upper estimates with $C=0$: So for any $\epsilon>0$ there exists $p_{\epsilon}\in\NN_{>0}$ such that for all $i\in\NN_{>0}$ one obtains $(0<)(M_{p_{\epsilon}+i})^{\frac{1}{p_{\epsilon}+i}}\le\epsilon^{\frac{i}{p_{\epsilon}+i}}(M_{p_{\epsilon}})^{\frac{1}{p_{\epsilon}+i}}$. Then, as $i\rightarrow+\infty$ we infer $\epsilon^{\frac{i}{p_{\epsilon}+i}}(M_{p_{\epsilon}})^{\frac{1}{p_{\epsilon}+i}}\rightarrow\epsilon$ and since $\epsilon>0$ was chosen arbitrary it follows that $\lim_{p\rightarrow+\infty}(M_p)^{1/p}=0$, too.}

\emph{But, finally, we recall that for any log-convex $\mathbf{M}\in\RR_{>0}^{\NN}$ this case can never occur since $p\mapsto\mu_p$ is non-decreasing and each $\mu_p$ is non-negative.}
\end{remark}

Let $\mathbf{M},\mathbf{N}\in\RR_{>0}^{\NN}$ be given, then write
\begin{itemize}
\item[$(*)$] $\mathbf{M}\le\mathbf{N}$ if $M_p\le N_p$ for all $p\in\NN$,

\item[$(*)$] $\mathbf{M}\hypertarget{preceq}{\preceq}\mathbf{N}$ if
$$\sup_{p\in\NN_{>0}}\left(\frac{M_p}{N_p}\right)^{1/p}<+\infty,$$

\item[$(*)$] $\mathbf{M}\hypertarget{mtriangle}{\vartriangleleft}\mathbf{N}$ if $\lim_{p\rightarrow+\infty}\left(\frac{M_p}{N_p}\right)^{1/p}=0$; i.e. if
\begin{equation*}\label{triangleestim}
\forall\;h>0\;\exists\;C_h\ge 1\;\forall\;p\in\NN:\;\;\;M_p\le C_hh^pN_p.
\end{equation*}
\item[$(*)$] $\mathbf{M}$ and $\mathbf{N}$ are called \emph{equivalent,} formally denoted by $\mathbf{M}\hypertarget{approx}{\approx}\mathbf{N}$, if $\mathbf{M}\hyperlink{preceq}{\preceq}\mathbf{N}$ and $\mathbf{N}\hyperlink{preceq}{\preceq}\mathbf{M}$.

\item[$(*)$] Obviously, $\mathbf{M}\hyperlink{mtriangle}{\vartriangleleft}\mathbf{N}$ implies $\mathbf{M}\hyperlink{preceq}{\preceq}\mathbf{N}$, but \hyperlink{mtriangle}{$\vartriangleleft$} is never reflexive.
\end{itemize}

Condition moderate growth is obviously preserved under equivalence of sequences. An important example, also for the considerations in this work, are the \emph{Gevrey sequences} $\mathbf{G}^s:=(p!^s)_{p\in\NN}$, $s>0$. Note that $\mathbf{G}^s$ is equivalent to $\overline{\mathbf{G}}^s$ with $\overline{\mathbf{G}}^s:=(p^{ps})_{p\in\NN}$ by \eqref{Stirlinglike} and each $\mathbf{G}^s$ is normalized, log-convex and satisfies \hyperlink{mg}{$(\on{mg})$}.\vspace{6pt}

For our purpose it is convenient to consider weight sequences in the following (general) sense:

\begin{definition}\label{defweightsequ}
$\mathbf{M}=(M_p)_{p\in\NN}$ is called a \emph{weight sequence} if
$$(\mathbf{M}_{\iota}=)\lim_{p\rightarrow+\infty}(M_p)^{1/p}=+\infty.$$
If $\mathbf{M}$ is a log-convex weight sequence, then in view of \cite[Sect. 2.2, Lemma 2.1]{regularnew} and see also the above Remark \ref{Lemma21rem}, we get $\lim_{p\rightarrow+\infty}\mu_p=+\infty$, too.
\end{definition}

Moreover, occasionally we consider the following (sub)set of weight sequences:
$$\hypertarget{LCset}{\mathcal{LC}}:=\{\mathbf{M}\in\RR_{>0}^{\NN}:\;\mathbf{M}\;\text{is normalized, log-convex},\;\lim_{p\rightarrow+\infty}(M_p)^{1/p}=+\infty\}.$$
We see that $\mathbf{M}\in\hyperlink{LCset}{\mathcal{LC}}$ if and only if $1=\mu_0\le\mu_1\le\dots$ and $\lim_{p\rightarrow+\infty}\mu_p=+\infty$ (see e.g. \cite[p. 104]{compositionpaper}) and there is a one-to-one correspondence between $\mathbf{M}$ and $\mu=(\mu_p)_p$ by taking $M_p:=\prod_{i=0}^p\mu_i$. Obviously, $\mathbf{G}^s\in\hyperlink{LCset}{\mathcal{LC}}$ for any $s>0$.\vspace{6pt}

For any given sequence $\mathbf{a}=(a_p)_p\in\RR_{>0}^{\NN}$ the \emph{upper Matuszewska index} $\alpha(\mathbf{a})$ is defined by
$$\alpha(\mathbf{a}):=\inf\{\alpha\in\RR: \frac{a_p}{p^{\alpha}}\;\text{is almost decreasing}\}=\inf\{\alpha\in\RR: \exists\;H\ge 1\;\forall\;1\le p\le q:\;\;\;\frac{a_q}{q^{\alpha}}\le H\frac{a_p}{p^{\alpha}}\},$$
and the \emph{lower Matuszewska index} $\beta(\mathbf{a})$ by
$$\beta(\mathbf{a}):=\sup\{\beta\in\RR: \frac{a_p}{p^{\beta}}\;\text{is almost increasing}\}=\sup\{\beta\in\RR: \exists\;H\ge 1\;\forall\;1\le p\le q:\;\;\;\frac{a_p}{p^{\beta}}\le H\frac{a_q}{q^{\beta}}\}.$$
$\beta(\mathbf{a})>0$ implies, in particular, $\lim_{p\rightarrow+\infty}a_p=+\infty$. For these definitions we refer to \cite[Sect. 3.2]{index} and the citations there. Crucially we also mention \cite[Thm. 3.10]{index} where the connection $\gamma(\mathbf{M})=\beta(\mu)$ is established even for any arbitrary $\mathbf{M}\in\RR_{>0}^{\NN}$. Here, $\gamma(\mathbf{M})$ denotes the growth index introduced by V. Thilliez in \cite{Thilliezdivision}. Note that in \cite{index} concerning the sequence of quotients a different notation has been used: There, $\mathbf{m}$ denotes $\mu$, more precisely one has the correspondence $m_p=\mu_{p+1}$ for all $p$. But this index shift does not effect our considerations; we refer to \cite[Rem. 3.8]{index} and the comments given after \cite[Cor. 3.12]{index}.

\subsection{Associated weight functions}\label{assofunctionsection}
For the following definition we refer to \cite[Chapitre I]{mandelbrojtbook}, \cite[Sect. 3]{Komatsu73} and the more recent work \cite{regularnew} in which more properties and explanations are summarized and non-standard cases are treated. Let $\mathbf{M}=(M_p)_{p\in\NN}\in\RR_{>0}^{\NN}$ be given, then the \emph{associated function} $\omega_{\mathbf{M}}: [0,+\infty)\rightarrow[0,+\infty)\cup\{+\infty\}$ is defined as follows:
\begin{equation}\label{assofunc}
\omega_{\mathbf{M}}(t):=\sup_{p\in\NN}\log\left(\frac{M_0t^p}{M_p}\right),\qquad t\ge 0,
\end{equation}
with the conventions $0^0:=1$ and $\log(0)=-\infty$. This ensures $\omega_{\mathbf{M}}(0)=0$ and $\omega_{\mathbf{M}}(t)\ge 0$ for any $t\ge 0$ since $\frac{t^0M_0}{M_0}=1$ for all $t\ge 0$. \eqref{assofunc} corresponds to \cite[$(3.1)$]{Komatsu73} and we immediately have that $\omega_{\mathbf{M}}$ is non-decreasing and satisfies $\lim_{t\rightarrow+\infty}\omega_{\mathbf{M}}(t)=+\infty$. However, the fact that $\omega_{\mathbf{M}}: [0,+\infty)\rightarrow[0,+\infty)$, i.e. that $\omega_{\mathbf{M}}$ is well defined, is equivalent to assume $\lim_{p\rightarrow+\infty}(M_p)^{1/p}=+\infty$; see \cite[Lemma 2.2]{regularnew} and compare this with Definition \ref{defweightsequ}. For any $s>0$ we have $\id^{1/s}\hyperlink{sim}{\sim}\omega_{\mathbf{G}^s}$; for a proof we refer e.g. to \cite[Ex. 2.9]{genLegendreconj}.\vspace{6pt}

Let $\mathbf{M}$ be a log-convex weight sequence, then one can compute $\mathbf{M}$ by involving $\omega_{\mathbf{M}}$ as follows, see \cite[Chapitre I, 1.4, 1.8]{mandelbrojtbook}, \cite[Prop. 3.2]{Komatsu73} and more recently \cite{regularnew}:
\begin{equation}\label{Prop32Komatsu}
M_p=M_0\sup_{t\ge 0}\frac{t^p}{\exp(\omega_{\mathbf{M}}(t))},\;\;\;p\in\NN.
\end{equation}
If $\mathbf{M}$ is a weight sequence which is not log-convex, then the right-hand side of \eqref{Prop32Komatsu} yields $\mathbf{M}^{\on{lc}}_p$ with $\mathbf{M}^{\on{lc}}$ denoting the \emph{log-convex minorant} of $\mathbf{M}$ and then $\mathbf{M}^{\on{lc}}$ is a log-convex weight sequence.

Moreover, for any log-convex weight sequence let us introduce the \emph{counting function}
\begin{equation*}\label{counting}
\Sigma_{\mathbf{M}}(t):=|\{p\ge 1:\;\;\;\mu_p\le t\}|,
\end{equation*}
and recall the following integral representation formula
\begin{equation}\label{assointrepr}
\omega_{\mathbf{M}}(t)=\int_0^t\frac{\Sigma_{\mathbf{M}}(u)}{u}du=\int_{\mu_1}^t\frac{\Sigma_{\mathbf{M}}(u)}{u}du;
\end{equation}
see \cite[1.8. III]{mandelbrojtbook}, \cite[$(3.11)$]{Komatsu73} and \cite[Lemma 2.5]{regularnew}. Consequently, \begin{equation*}\label{assovanishing}
\forall\;t\in[0,\mu_1]:\;\;\;\omega_{\mathbf{M}}(t)=0,
\end{equation*}
and so for any $\mathbf{M}\in\hyperlink{LCset}{\mathcal{LC}}$ the function $\omega_{\mathbf{M}}$ vanishes on $[0,1]$. Next let us recall \cite[Lemma 2.4]{regularnew}.

\begin{lemma}\label{lemma1}
Let $\mathbf{M}=(M_p)_{p\in\NN}$ be a log-convex weight sequence and hence $\mathbf{M}_{\iota}=\lim_{p\rightarrow+\infty}(M_p)^{1/p}=\lim_{p\rightarrow+\infty}\mu_p=+\infty$. Then
\begin{equation}\label{lemma1equ}
\omega_{\mathbf{M}}(t)=0,\;\text{for}\;t\in[0,\mu_1],\;\;\;\;\omega_{\mathbf{M}}(t)=\log\left(\frac{M_0t^p}{M_p}\right)\;\text{for}\;t\in[\mu_p,\mu_{p+1}],\;p\ge 1.
\end{equation}
\end{lemma}

For concrete applications, when one is (only) interested in the associated function w.l.o.g. one can assume that the defining sequence is log-convex. This fact holds since the geometric regularizing procedure introducing $\omega_{\mathbf{M}}$ gives, in particular, that $\omega_{\mathbf{M}}=\omega_{\mathbf{M}^{\on{lc}}}$; see again \cite{regularnew}.

\begin{remark}\label{indexweightsequrem}
\emph{A log-convex weight sequence with $M_0=1$ is a weight sequence in the sense of \cite{index}; see \cite[Sect. 3.1]{index}. However, if $\mathbf{M}$ is a (log-convex) weight sequence with $M_0\neq 1$, then introduce $\widetilde{\mathbf{M}}$ with $\widetilde{M}_p:=\frac{M_p}{M_0}$. So $\widetilde{\mathbf{M}}$ is a (log-convex) weight sequence with $\widetilde{M}_0=1$ and $\widetilde{\mathbf{M}}\hyperlink{approx}{\approx}\mathbf{M}$ holds but even $\widetilde{\mu}_p=\mu_p$ for all $p$ is valid.}

\emph{Moreover, $\omega_{\widetilde{\mathbf{M}}}=\omega_{\mathbf{M}}$ since}
$$\forall\;t\ge 0:\;\;\;\omega_{\widetilde{\mathbf{M}}}(t)=\sup_{p\in\NN}\log\left(\frac{\widetilde{M}_0t^p}{\widetilde{M}_p}\right)=\sup_{p\in\NN}\log\left(\frac{M_0t^p}{M_p}\right)=\omega_{\mathbf{M}}(t).$$ \emph{Thus when involving results from \cite{index} the assumption $M_0=1$ is not restricting the generality since we apply the statements to $\widetilde{\mathbf{M}}$ and by the above equality of associated functions and quotients they transfer to $\mathbf{M}$ immediately.}

\emph{Finally, when $\mathbf{M}$ is log-convex then by applying \cite[Lemmas 2.0.4 \& 2.0.6]{diploma} to $\widetilde{\mathbf{M}}$ we get that $p\mapsto(M_p/M_0)^{1/p}$ is non-decreasing, $M_pM_q\le M_0M_{p+q}$ for all $p,q\in\NN$ and $(M_p/M_0)^{1/p}\le\mu_p$ for all $p\in\NN_{>0}$.}
\end{remark}

\subsection{The conjugate sequence $\mathbf{M}^{\ast}$}\label{conjsequsection}
Let us recall some definitions and properties listed and stated in \cite[Sect. 2.5 \& 2.6]{microclasses}; see also there for more details and explanations. Let $\mathbf{M}\in\RR_{>0}^{\NN}$ be given, then the \emph{conjugate sequence} $\mathbf{M}^{\ast}=(M^{\ast}_p)_{p\in\NN}$ is defined by
\begin{equation}\label{conjsequdef}
M^{\ast}_p:=\frac{p!}{M_p}=\frac{1}{m_p},\;\;\;p\in\NN;
\end{equation}
i.e. $\mathbf{M}^{\ast}=\frac{1}{\mathbf{m}}$. The corresponding sequence of quotients $\mu^{\ast}=(\mu^{\ast}_p)_{p\in\NN}$ is given by $\mu^{\ast}_p=\frac{p}{\mu_p}$, $p\in\NN_{>0}$, and set $\mu^{\ast}_0:=1$. We gather several immediate properties:

\begin{itemize}
\item[$(a)$] Set $\mathbf{M}^{\ast\ast}:=(\mathbf{M}^{\ast})^{\ast}$, then $\mathbf{M}^{\ast\ast}=\mathbf{M}$ and $\mathbf{M}^{\ast}\cdot\mathbf{M}\equiv\mathbf{G}^1$.

\item[$(b)$] The following are equivalent:
\begin{itemize}
\item[$(*)$] $\mathbf{M}^{\ast}$ is log-convex,

\item[$(*)$] $\mathbf{m}$ is \emph{log-concave,} see \cite[$(2.5)$]{microclasses}, i.e.
$$\forall\;p\in\NN_{>0}:\;\;\;m_p^2\ge m_{p+1}m_{p-1}.$$
\item[$(*)$] $p\mapsto\frac{\mu_p}{p}$ is non-increasing.
\end{itemize}
If a sequence $\mathbf{N}$ is log-concave, then following the arguments in \cite[Lemma 2.0.4]{diploma} we get that $p\mapsto(N_p/N_0)^{1/p}$ is non-increasing. Consequently, if $\mathbf{M}^{\ast}$ is log-convex and $M_0^{\ast}=m_0=1$, then $p\mapsto(m_p)^{1/p}$ is non-increasing.

\item[$(c)$] If both $\mathbf{M}$ and $\mathbf{M}^{\ast}$ are log-convex, then $\mathbf{M}^{\ast}\cdot\mathbf{M}\equiv\mathbf{G}^1$ implies $\mathbf{M}_{\iota}\mathbf{M}^{\ast}_{\iota}=+\infty$.

\item[$(d)$] It is immediate that
$$\mathbf{M}\hyperlink{mtriangle}{\vartriangleleft}\mathbf{G}^1\Leftrightarrow\lim_{p\rightarrow+\infty}(m_p)^{1/p}=0\Leftrightarrow \lim_{p\rightarrow+\infty}(M^{\ast}_p)^{1/p}=\mathbf{M}^{\ast}_{\iota}=+\infty,$$
and this is equivalent to the fact that $\omega_{\mathbf{M}^{\ast}}$ is well defined; recall \cite[Lemma 2.2]{regularnew} applied to $\mathbf{M}^{\ast}$. Thus, whenever involving $\omega_{\mathbf{M}^{\ast}}$ it is natural to assume this growth property on $\mathbf{M}$.

If $\mathbf{M}^{\ast}$ is also log-convex, then $\mathbf{M}\hyperlink{mtriangle}{\vartriangleleft}\mathbf{G}^1$ if and only if $\lim_{p\rightarrow+\infty}\mu^{\ast}_p=+\infty$; see \cite[Sect. 2.2]{regularnew}.

\item[$(e)$] $\mathbf{M}\le\mathbf{N}$ if and only if $\mathbf{N}^{\ast}\le\mathbf{M}^{\ast}$ and $\mathbf{M}\hyperlink{preceq}{\preceq}\mathbf{N}$ if and only if $\mathbf{N}^{\ast}\hyperlink{preceq}{\preceq}\mathbf{M}^{\ast}$. Thus $\mathbf{M}$, $\mathbf{N}$ are equivalent if and only if the conjugate sequences $\mathbf{M}^{\ast}$, $\mathbf{N}^{\ast}$ are so.

    Moreover, $\mathbf{M}\hyperlink{mtriangle}{\vartriangleleft}\mathbf{N}$ if and only if $\mathbf{m}\hyperlink{mtriangle}{\vartriangleleft}\mathbf{n}$ if and only if $\mathbf{N}^{\ast}\hyperlink{mtriangle}{\vartriangleleft}\mathbf{M}^{\ast}$.

\item[$(f)$] $\mathbf{M}^{*}$ satisfies \hyperlink{mg}{$(\on{mg})$} if and only if $m_pm_q\le C^{p+q+1}m_{p+q}$ if and only if $\binom{p+q}{p}M_pM_q\le C^{p+q+1}M_{p+q}$ for some $C\ge 1$ and all $p,q\in\NN$. If $\mathbf{M}$ is log-convex, then Remark \ref{indexweightsequrem} gives $\binom{p+q}{p}M_pM_q\le 2^{p+q}M_0M_{p+q}\le\max\{2,M_0\}^{p+q+1}M_{p+q}$ for all $p,q$, so the desired estimate holds with $C:=\max\{2,M_0\}$ and $\mathbf{M}^{*}$ satisfies \hyperlink{mg}{$(\on{mg})$}.
\end{itemize}

Next, we verify the converse of $(f)$ and first, motivated by comment $(b)$ above, we prove the following technical result.

\begin{lemma}\label{rootalmostdecrlemma}
Let $\mathbf{M}\in\RR_{>0}^{\NN}$ be given and recall the notation $m_p:=\frac{M_p}{p!}$, $p\in\NN$. Then the following are equivalent:
\begin{itemize}
\item[$(i)$] $\mathbf{M}$ satisfies
\begin{equation}\label{rootalmostdecrlemmaequ}
\exists\;H\ge 1\;\forall\;p\in\NN_{>0}:\;\;\;(m_{p+1})^{1/(p+1)}\le H^{1/(p+1)}(m_p)^{1/p}.
\end{equation}

\item[$(ii)$] $\mathbf{M}$ satisfies
\begin{equation}\label{rootalmostdecrlemmaequ1}
\exists\;A\ge 1\;\forall\;p\in\NN_{>0}:\;\;\;\mu_{p+1}\le A(M_p)^{1/p}.
\end{equation}
\end{itemize}
The proof shows the following correspondences: If \eqref{rootalmostdecrlemmaequ} holds with some $H$, then in \eqref{rootalmostdecrlemmaequ1} we choose $A:=2eH$ and, conversely, \eqref{rootalmostdecrlemmaequ1} implies \eqref{rootalmostdecrlemmaequ} with $H:=A$.
\end{lemma}

\demo{Proof}
$(i)\Rightarrow(ii)$: Let $p\in\NN_{>0}$, then we have
\begin{align*}
&(m_{p+1})^{1/(p+1)}\le H^{1/(p+1)}(m_p)^{1/p}\Leftrightarrow(m_{p+1})^p\le H^p(m_p)^{p+1}\Leftrightarrow\left(\frac{M_{p+1}}{(p+1)!}\right)^p\le H^p\left(\frac{M_p}{p!}\right)^{p+1}
\\&
\Leftrightarrow\left(\frac{M_{p+1}}{M_p}\right)^p\le H^p\frac{(p+1)!^p}{p!^{p+1}}M_p\Leftrightarrow\mu_{p+1}^p\le H^p(p+1)^p\frac{M_p}{p!}\Leftrightarrow\mu_{p+1}\le H(p+1)(m_p)^{1/p}.
\end{align*}
Moreover, $(p+1)(m_p)^{1/p}\le 2p(m_p)^{1/p}\le 2e(p!)^{1/p}(m_p)^{1/p}=2e(M_p)^{1/p}$ is valid by \eqref{Stirlinglike}. Combining both estimates we have verified \eqref{rootalmostdecrlemmaequ1} with $A:=2eH$.\vspace{6pt}

$(ii)\Rightarrow(i)$: Let $p\in\NN_{>0}$ be arbitrary, and then
\begin{align*}
&\mu_{p+1}\le A(M_p)^{1/p}\Leftrightarrow\frac{(p+1)!m_{p+1}}{p!m_p}\le A(p!m_p)^{1/p}\Leftrightarrow m_{p+1}\le A\frac{(p!)^{1/p}}{p+1}(m_p)^{(p+1)/p}
\\&
\Leftrightarrow(m_{p+1})^{1/(p+1)}\le A^{1/(p+1)}\left(\frac{(p!)^{1/p}}{p+1}\right)^{1/(p+1)}(m_p)^{1/p}.
\end{align*}
Since $\left(\frac{(p!)^{1/p}}{p+1}\right)^{1/(p+1)}\le 1\Leftrightarrow(p!)^{1/p}\le p+1$ is clear for any $p\in\NN_{>0}$ the estimate \eqref{rootalmostdecrlemmaequ} is shown with $H:=A$ as desired.
\qed\enddemo

Concerning the appearing conditions in Lemma \ref{rootalmostdecrlemma} we point out:

\begin{remark}\label{rootalmostdecrlemmarem}
\emph{The choice $H=1$ in \eqref{rootalmostdecrlemmaequ} precisely means that $p\mapsto(m_p)^{1/p}$ is non-increasing, whereas it is known that \eqref{rootalmostdecrlemmaequ1} is equivalent to \hyperlink{mg}{$(\on{mg})$} if $\mathbf{M}\in\RR_{>0}^{\NN}$ is log-convex; we refer to \cite[Lemma 2.2 $(1)\Leftrightarrow(2)$, Rem. 2.3, $(2.10)$]{whitneyextensionweightmatrix}, the proof of \cite[Thm. 1]{matsumoto} and also to the detailed discussion in \cite{modgrowthstrange}. Note that for both implications in this equivalence log-convexity is required in the arguments but, however, $\lim_{p\rightarrow+\infty}(M_p)^{1/p}=+\infty$ is not necessary for this part. On the other side, this characterization breaks down in general in the mixed setting as it is shown in \cite{modgrowthstrange} and \cite{modgrowthstrangeII}.}
\end{remark}

\begin{corollary}\label{rootalmostdecrlemmacor}
If both $\mathbf{M}^{*}$ and $\mathbf{M}$ are log-convex, then $\mathbf{M}$ has to satisfy \hyperlink{mg}{$(\on{mg})$}.
\end{corollary}

\demo{Proof}
Via comment $(b)$ above, if $\mathbf{M}^{\ast}$ is log-convex, then $p\mapsto(m_p/m_0)^{1/p}$ is non-increasing and so $$\forall\;p\in\NN_{>0}:\;\;\;m_0^{-\frac{1}{p(p+1)}}(m_p)^{1/p}=m_0^{\frac{1}{p+1}-\frac{1}{p}}(m_p)^{1/p}\ge(m_{p+1})^{1/(p+1)}.$$
Since $m_0^{-1/p}\le(1-\epsilon)^{-1}\Leftrightarrow(1-\epsilon)^p\le m_0$ holds for any $p\in\NN_{>0}$ when $\epsilon\in(0,1)$ is chosen sufficiently close to $1$, we get that $H^{1/(p+1)}\ge m_0^{-\frac{1}{p(p+1)}}\Leftrightarrow H\ge m_0^{-1/p}$ holds for $H$ chosen sufficiently large and so \eqref{rootalmostdecrlemmaequ} is valid. Then Lemma \ref{rootalmostdecrlemma} and Remark \ref{rootalmostdecrlemmarem} yield the conclusion.
\qed\enddemo

Finally, we briefly recall the technical but crucial and useful comments from \cite[Lemma 2.3 $(a)$, Rem. 2.4, Rem. A.5 $(iv)$]{microclasses}.

\begin{remark}\label{Rem24reloaded}
\emph{Let $\mathbf{M}\in\RR_{>0}^{\NN}$ be given and assume that $\left(\frac{\mu_p}{p}\right)_p$ is \emph{almost decreasing;} i.e.}
\begin{equation}\label{almostdecr}
\exists\;H\ge 1\;\forall\;1\le p\le q:\;\;\;\frac{\mu_q}{q}\le H\frac{\mu_p}{p}.
\end{equation}
\emph{Note that by definition $\alpha(\mu)<1$ implies \eqref{almostdecr} whereas this property implies $\alpha(\mu)\le 1$. Then set}
\begin{equation}\label{almoststuffequ}
\lambda_p:=H^{-1}p\sup_{q\ge p}\frac{\mu_q}{q},\;\;\;p\ge 1,\hspace{20pt}\lambda_0:=1,
\end{equation}
\emph{and introduce the sequence $\mathbf{L}=(L_p)_{p\in\NN}$ via $L_p:=\prod_{i=0}^p\lambda_i$. Then we get (see \cite[Rem. 2.4]{microclasses}):}

\begin{itemize}
\item[$(*)$] \emph{$L_0=1$ holds.}

\item[$(*)$] \emph{$\mathbf{l}=(l_p)_{p\in\NN}$ given by $l_p:=\frac{L_p}{p!}$ is log-concave; i.e. $p\mapsto\frac{\lambda_p}{p}$ is non-increasing.}

\item[$(*)$] \emph{$\mathbf{M}\hyperlink{approx}{\approx}\mathbf{L}$ is valid. If $\lim_{p\rightarrow+\infty}(m_p)^{1/p}=0$ resp. $\lim_{p\rightarrow+\infty}(M_p)^{1/p}=+\infty$ holds, then by this equivalence $\mathbf{l}$ resp. $\mathbf{L}$ shares the particular growth property, too.}

\item[$(*)$] \emph{If $\mathbf{M}$ is log-convex, then $\mathbf{L}$, too.}
\end{itemize}
\emph{Moreover, in view of the general construction in \cite[Rem. A.5 $(iv)$]{microclasses}, if $\mathbf{L}$ is a log-convex weight sequence, then by applying a technical modification for $\mathbf{L}$ for finitely many indices (at the beginning), one can even achieve that there exists $\widetilde{\mathbf{L}}\in\hyperlink{LCset}{\mathcal{LC}}$ such that $\widetilde{L}_0=1=\widetilde{L}_1$, $\widetilde{\mathbf{l}}$ is log-concave, $\mathbf{M}\hyperlink{approx}{\approx}\widetilde{\mathbf{L}}$ and, finally, $\widetilde{\mathbf{L}},\widetilde{\mathbf{L}}^{\ast}\in\hyperlink{LCset}{\mathcal{LC}}$.}

\emph{Indeed, by the above comments this situation applies when $\mathbf{M}$ is log-convex, $\left(\frac{\mu_p}{p}\right)_p$ is almost decreasing, $\lim_{p\rightarrow+\infty}(m_p)^{1/p}=0$ and $\lim_{p\rightarrow+\infty}(M_p)^{1/p}=+\infty$. In this case in order to lighten notation we simply write again $\mathbf{L}$, $\mathbf{L}^{\ast}$ for these sequences.}
\end{remark}

The importance of the previous remark is given by the fact that for given $\mathbf{M}$ satisfying some growth requirements one can find an equivalent and very regular sequence $\mathbf{L}$ such that $\mathbf{L}^{\ast}$ is also very regular (and equivalent to $\mathbf{M}^{\ast}$).

\subsection{Conjugate associated weight functions}\label{conjassofctsection}
We are interested in the situation when $\omega^{\ast}_{\mathbf{M}}$ is well defined and characterize this in terms of $\mathbf{M}$.

\begin{lemma}\label{C2assoweightfctlemma}
Let $\mathbf{M}$ be a log-convex weight sequence. Then the following are equivalent:
\begin{itemize}
\item[$(i)$] $\omega^{\ast}_{\mathbf{M}}$ is well defined.

\item[$(ii)$] $\omega_{\mathbf{M}}$ satisfies $(\mathfrak{C}_2)$; i.e. $t=o(\omega_{\mathbf{M}}(t))$ as $t\rightarrow+\infty$ equivalently $\omega_{\mathbf{M}}\hyperlink{omvartriangle}{\vartriangleleft}\id$.

\item[$(iii)$] $\omega_{\mathbf{G}^1}(t)=o(\omega_{\mathbf{M}}(t))$ as $t\rightarrow+\infty$; i.e. $\omega_{\mathbf{M}}\hyperlink{omvartriangle}{\vartriangleleft}\omega_{\mathbf{G}^1}$.

\item[$(iv)$] $\mathbf{M}\hyperlink{mtriangle}{\vartriangleleft}\mathbf{G}^1$ holds; i.e.
\begin{equation}\label{C2assoweightfctlemmaequ}
\lim_{p\rightarrow+\infty}\left(\frac{M_p}{p!}\right)^{1/p}=\lim_{p\rightarrow+\infty}(m_p)^{1/p}=0.
\end{equation}
\end{itemize}
\end{lemma}

The crucial growth condition \eqref{C2assoweightfctlemmaequ} is non-standard in the theory of ultradifferentiable function classes and related weighted settings.

\demo{Proof}
First, $(i)\Leftrightarrow(ii)$ holds by Lemma \ref{conjugatewelldeflemma} applied to $\omega^{\ast}_{\mathbf{M}}$ and $(ii)\Leftrightarrow(iii)$ is valid because $\id\hyperlink{sim}{\sim}\omega_{\mathbf{G}^1}$. This equivalence is well-known; for an explicit proof we refer e.g. to the case $s=1$ in \cite[Ex. 2.9]{genLegendreconj}.

And this equivalence implies, in particular, that $\omega_{\mathbf{G}^1}$ satisfies both \hyperlink{om1}{$(\omega_1)$} and \hyperlink{om6}{$(\omega_6)$} since these conditions are obviously true for $\id$.\vspace{6pt}

$(iii)\Rightarrow(iv)$: The assumption precisely means that for all $\epsilon>0$ there exists some $C_{\epsilon}\ge 1$ such that $\omega_{\mathbf{G}^1}(t)\le\epsilon\omega_{\mathbf{M}}(t)+C_{\epsilon}$ for all $t\ge 0$. Then let $h>0$ small be given and choose $n\in\NN_{>0}$ such that $h^{-1}\le 2^n$. Since $\omega_{\mathbf{G}^1}$ has \hyperlink{om1}{$(\omega_1)$} by iteration we obtain $\omega_{\mathbf{G}^1}(h^{-1}t)\le\omega_{\mathbf{G}^1}(2^nt)\le L\omega_{\mathbf{G}^1}(t)+L$ for some $L\ge 1$ depending on $h$ and for all $t\ge 0$. We apply the previous relation to $\epsilon:=L^{-1}$ and obtain:
$$\forall\;t\ge 0:\;\;\;\omega_{\mathbf{G}^1}(h^{-1}t)\le L\omega_{\mathbf{G}^1}(t)+L\le\omega_{\mathbf{M}}(t)+LC_{L^{-1}}+L.$$
Set $L_1:=LC_{L^{-1}}+L$, then by using this estimate and involving \eqref{Prop32Komatsu} we obtain for all $p\in\NN$:
\begin{align*}
M_p&=M_0\sup_{t\ge 0}\frac{t^p}{\exp(\omega_{\mathbf{M}}(t))}\le M_0\exp(L_1)\sup_{t\ge 0}\frac{t^p}{\exp(\omega_{\mathbf{G}^1}(h^{-1}t))}
\\&
=M_0\exp(L_1)h^p\sup_{s\ge 0}\frac{s^p}{\exp(\omega_{\mathbf{G}^1}(s))}=M_0\exp(L_1)h^pG^1_p.
\end{align*}
For the last equality note that $G^1_0=0!=1$. Summarizing,
$$\forall\;h>0\;\exists\;L_1\ge 1\;\forall\;p\in\NN:\;\;\;M_p\le M_0\exp(L_1)h^pG^1_p,$$
which gives $\mathbf{M}\hyperlink{mtriangle}{\vartriangleleft}\mathbf{G}^1$.\vspace{6pt}

$(iv)\Rightarrow(iii)$: By assumption
$$\forall\;h>0\;\exists\;C_h\ge 1\;\forall\;p\in\NN\;\forall\;t\ge 0:\;\;\;\frac{M_0t^p}{p!}\le C_h\frac{M_0(ht)^p}{M_p},$$
hence by definition $\omega_{\mathbf{G}^1}(t)\le\omega_{\mathbf{M}}(ht)+\log(C_h/M_0)$ for all $t\ge 0$.

Let now $A\ge 1$ large be given and fixed. Choose $n\in\NN_{>0}$ such that $A\le 2^n$ and by iterating \hyperlink{om6}{$(\omega_6)$} we obtain $A\omega_{\mathbf{G}^1}(t)\le 2^n\omega_{\mathbf{G}^1}(t)\le\omega_{\mathbf{G}^1}(Ht)+H$ for some $H\ge 1$ depending on chosen $A$ and for all $t\ge 0$. When applying the above relation to $h:=H^{-1}$ we arrive at
$$\forall\;t\ge 0:\;\;\;A\omega_{\mathbf{G}^1}(t)\le\omega_{\mathbf{G}^1}(Ht)+H\le\omega_{\mathbf{M}}(t)+\log(C_{H^{-1}}/M_0)+H,$$
hence $\limsup_{t\rightarrow+\infty}\frac{\omega_{\mathbf{G}^1}(t)}{\omega_{\mathbf{M}}(t)}\le\frac{1}{A}$. Since $A\ge 1$ can be chosen arbitrary large this shows $\omega_{\mathbf{G}^1}(t)=o(\omega_{\mathbf{M}}(t))$ and are done.
\qed\enddemo

\emph{Note:} Condition \eqref{C2assoweightfctlemmaequ} admits the possibility to check the desired growth property for given explicit examples directly and remark that log-convexity is only required for $(iii)\Rightarrow(iv)$. By comment $(d)$ in Section \ref{conjsequsection} and in view of Lemmas \ref{C2assoweightfctlemma} and \ref{conjugatewelldeflemma} we obtain the following consequence:

\begin{corollary}\label{conjwelldefcor}
Let $\mathbf{M}$ be a log-convex weight sequence, then the following are equivalent:
\begin{itemize}
\item[$(*)$] $\lim_{p\rightarrow+\infty}(m_p)^{1/p}=0$ is valid.

\item[$(*)$] $\omega^{\ast}_{\mathbf{M}}$ is well defined.

\item[$(*)$] $\omega_{\mathbf{M}^{\ast}}$ is well defined.
\end{itemize}
\end{corollary}
Note that in Corollary \ref{conjwelldefcor} it is not required necessarily to have log-convexity for $\mathbf{M}^{\ast}$, i.e. log-concavity for $\mathbf{m}$.

Next we investigate the growth relations: For given weight sequences $\mathbf{M}$, $\mathbf{N}$ the relation  $\mathbf{M}\hyperlink{preceq}{\preceq}\mathbf{N}$ precisely means
$$\exists\;C,h>0\;\forall\;t\ge 0\;\forall\;p\in\NN:\;\;\;\frac{N_0t^p}{N_p}\le\frac{N_0}{M_0}C\frac{M_0(th)^p}{M_p},$$
which gives by definition
\begin{equation}\label{assofctrel}
\exists\;D>0\;\exists\;h>0\;\forall\;t\ge 0:\;\;\;\omega_{\mathbf{N}}(t)\le\omega_{\mathbf{M}}(ht)+D;
\end{equation}
i.e. $\omega_{\mathbf{M}}\hyperlink{ompreceqc}{\preceq_{\mathfrak{c}}}\omega_{\mathbf{N}}$. Thus $\mathbf{M}\hyperlink{approx}{\approx}\mathbf{N}$ yields $\omega_{\mathbf{M}}\hyperlink{simc}{\sim_{\mathfrak{c}}}\omega_{\mathbf{N}}$. Similarly, $\mathbf{M}\hyperlink{mtriangle}{\vartriangleleft}\mathbf{N}$ implies $\omega_{\mathbf{M}}\hyperlink{omvartrianglec}{\vartriangleleft_{\mathfrak{c}}}\omega_{\mathbf{N}}$. Using this we can show the following:

\begin{lemma}\label{growthrellemma}
Let $\mathbf{M}$, $\mathbf{N}$ be log-convex weight sequences.
\begin{itemize}
\item[$(a)$] Assume that $\mathbf{M}\hyperlink{preceq}{\preceq}\mathbf{N}$ and that $\lim_{p\rightarrow+\infty}(n_p)^{1/p}=0$ holds.  Then $\omega_{\mathbf{N}^{\ast}}$ and $\omega_{\mathbf{M}^{\ast}}$ are well defined and
\begin{equation}\label{assofctrel1}
\exists\;D>0\;\exists\;h>0\;\forall\;t\ge 0:\;\;\;\omega_{\mathbf{M}^{\ast}}(t)\le\omega_{\mathbf{N}^{\ast}}(ht)+D;
\end{equation}
i.e. $\omega_{\mathbf{N}^{\ast}}\hyperlink{ompreceqc}{\preceq_{\mathfrak{c}}}\omega_{\mathbf{M}^{\ast}}$.

Moreover, $\omega^{\ast}_{\mathbf{M}}$ and $\omega^{\ast}_{\mathbf{M}}$ are well defined and
\begin{equation}\label{assofctrel2}
\exists\;D>0\;\exists\;h>0\;\forall\;s\ge 0:\;\;\;\omega^{\ast}_{\mathbf{M}}(s)\le\omega^{\ast}_{\mathbf{N}}(sh)+D;
\end{equation}
i.e. $\omega^{\ast}_{\mathbf{N}}\hyperlink{ompreceqc}{\preceq_{\mathfrak{c}}}\omega^{\ast}_{\mathbf{M}}$.

\item[$(b)$] Assume that $\mathbf{M}\hyperlink{mtriangle}{\vartriangleleft}\mathbf{N}$ and that $\lim_{p\rightarrow+\infty}(n_p)^{1/p}=0$ holds. Then even $\omega_{\mathbf{N}^{\ast}}\hyperlink{omvartriangle}{\vartriangleleft_{\mathfrak{c}}}\omega_{\mathbf{M}^{\ast}}$ and $\omega^{\ast}_{\mathbf{N}}\hyperlink{omvartriangle}{\vartriangleleft_{\mathfrak{c}}}\omega^{\ast}_{\mathbf{M}}$ and all (associate) weight functions under consideration are well defined.
\end{itemize}
\end{lemma}

Therefore, when $\mathbf{M}\hyperlink{approx}{\approx}\mathbf{N}$ and either $\lim_{p\rightarrow+\infty}(m_p)^{1/p}=0$ or $\lim_{p\rightarrow+\infty}(n_p)^{1/p}=0$, then $\omega_{\mathbf{M}}\hyperlink{simc}{\sim_{\mathfrak{c}}}\omega_{\mathbf{N}}$, $\omega_{\mathbf{M}^{\ast}}\hyperlink{simc}{\sim_{\mathfrak{c}}}\omega_{\mathbf{N}^{\ast}}$, $\omega^{\ast}_{\mathbf{N}}\hyperlink{simc}{\sim_{\mathfrak{c}}}\omega^{\ast}_{\mathbf{M}}$ and all (associate) weight functions under consideration are well defined. Concerning this additional information note that equivalence is preserved under taking the conjugate operation; recall comment $(e)$ in Section \ref{conjsequsection}.

\demo{Proof}
$(a)$ First, by comment $(e)$ in Section \ref{conjsequsection} we have $\mathbf{N}^{\ast}\hyperlink{preceq}{\preceq}\mathbf{M}^{\ast}$. Second, by $\lim_{p\rightarrow+\infty}(n_p)^{1/p}=0$ and $\mathbf{M}\hyperlink{preceq}{\preceq}\mathbf{N}$ also $\lim_{p\rightarrow+\infty}(m_p)^{1/p}=0$. Thus both $\omega_{\mathbf{N}^{\ast}}$ and $\omega_{\mathbf{M}^{\ast}}$ are well defined, see Corollary \ref{conjwelldefcor}, and similarly like \eqref{assofctrel} estimate \eqref{assofctrel1} follows.\vspace{6pt}

And Corollary \ref{conjwelldefcor} also implies that both $\omega^{\ast}_{\mathbf{M}}$ and $\omega^{\ast}_{\mathbf{M}}$ are well defined. Since $\omega_{\mathbf{M}}\hyperlink{ompreceqc}{\preceq_{\mathfrak{c}}}\omega_{\mathbf{N}}$ holds we have that $(i)$ in Lemma \ref{equivalencelemma1} implies $\omega^{\ast}_{\mathbf{N}}\hyperlink{ompreceqc}{\preceq_{\mathfrak{c}}}\omega^{\ast}_{\mathbf{M}}$.\vspace{6pt}

$(b)$ The conclusions follow analogously by taking into account $(ii)$ in Lemma \ref{equivalencelemma1}, comment $(e)$ in Section \ref{conjsequsection} and Corollary \ref{conjwelldefcor}.
\qed\enddemo

\subsection{The conjugate sequence and the generalized upper and lower Legendre envelopes}\label{weightsequenceLegendresection}
We study the interaction of the conjugate operation and the Legendre envelopes introduced in Section \ref{generalizedLegendresection} in the weight sequence setting when involving the corresponding conjugate sequence. First we verify an auxiliary result which should be compared with \cite[Lemma 3.4, Thm. 4.14 \& 4.15]{genLegendreconj}.

\begin{lemma}\label{conjsequlemma}
Let $\mathbf{M}$, $\mathbf{P}$, $\mathbf{N}$, $\mathbf{Q}$ be weight sequences.
\begin{itemize}
\item[$(i)$] If $\mathbf{M}\hyperlink{preceq}{\preceq}\mathbf{P}$ and $\mathbf{N}\hyperlink{preceq}{\preceq}\mathbf{Q}$, then
$\omega_{\mathbf{M}}\check{\star}\omega_{\mathbf{N}}\hyperlink{ompreceqc}{\preceq_{\mathfrak{c}}}\omega_{\mathbf{P}}\check{\star}\omega_{\mathbf{Q}}$ holds whereas $\mathbf{M}\hyperlink{mtriangle}{\vartriangleleft}\mathbf{P}$ and $\mathbf{N}\hyperlink{mtriangle}{\vartriangleleft}\mathbf{Q}$ imply
$\omega_{\mathbf{M}}\check{\star}\omega_{\mathbf{N}}\hyperlink{omvartrianglec}{\vartriangleleft_{\mathfrak{c}}}\omega_{\mathbf{P}}\check{\star}\omega_{\mathbf{Q}}$.

\item[$(ii)$] If $\mathbf{M}\hyperlink{preceq}{\preceq}\mathbf{P}$ and $\mathbf{Q}\hyperlink{preceq}{\preceq}\mathbf{N}$, then
$\omega_{\mathbf{M}}\widehat{\star}\omega_{\mathbf{N}}\hyperlink{ompreceqc}{\preceq_{\mathfrak{c}}}\omega_{\mathbf{P}}\widehat{\star}\omega_{\mathbf{Q}}$ holds whereas $\mathbf{M}\hyperlink{mtriangle}{\vartriangleleft}\mathbf{P}$ and $\mathbf{Q}\hyperlink{mtriangle}{\vartriangleleft}\mathbf{N}$ imply
$\omega_{\mathbf{M}}\widehat{\star}\omega_{\mathbf{N}}\hyperlink{omvartrianglec}{\vartriangleleft_{\mathfrak{c}}}\omega_{\mathbf{P}}\widehat{\star}\omega_{\mathbf{Q}}$.
\end{itemize}
\end{lemma}
Concerning $(ii)$ note that in order to give this part a meaning one shall naturally assume that $\omega_{\mathbf{M}}\widehat{\star}\omega_{\mathbf{N}}$ is well defined, see \cite[Prop. 5.2, Def. 5.5]{genLegendreconj}. Recall that by \cite[Prop. 5.2, $(I)(ii)\Rightarrow(i)$]{genLegendreconj} we get that $\mathbf{N}\hyperlink{mtriangle}{\vartriangleleft}\mathbf{M}$ ensures this fact and one obtains a characterization provided $\mathbf{N}$ is log-convex. The relations between the sequences stated in $(ii)$ give that $\mathbf{N}\hyperlink{mtriangle}{\vartriangleleft}\mathbf{M}$ implies $\mathbf{Q}\hyperlink{mtriangle}{\vartriangleleft}\mathbf{P}$ and thus this result is consistent with \cite[Prop. 5.2]{genLegendreconj}.

\demo{Proof}
$(i)$ By assumption we can find $C_1,C_2,h_1,h_2>0$ such that $M_j\le C_1h_1^jP_j$ and $N_j\le C_2h_2^jQ_j$ for all $j\in\NN$. Set $C:=\max\{C_1,C_2\}$ and $h:=\max\{h_1,h_2\}$ and the definition of associated weight functions gives $\omega_{\mathbf{P}}(t)\le\omega_{\mathbf{M}}(ht)+\log(CP_0/M_0)$ and $\omega_{\mathbf{Q}}(t)\le\omega_{\mathbf{N}}(ht)+\log(CQ_0/N_0)$ for all $t\ge 0$. Thus \eqref{wedgeformula} implies for all $t\ge 0$:
\begin{align*}
&\omega_{\mathbf{P}}\check{\star}\omega_{\mathbf{Q}}(t)=\inf_{s>0}\{\omega_{\mathbf{P}}(s)+\omega_{\mathbf{Q}}(t/s)\}\le\inf_{s>0}\{\omega_{\mathbf{M}}(hs)+\omega_{\mathbf{N}}(ht/s)\}+\log(C^2P_0Q_0/(M_0N_0))
\\&
=\inf_{u>0}\{\omega_{\mathbf{M}}(u)+\omega_{\mathbf{N}}(h^2t/u)\}+\log(C^2P_0Q_0/(M_0N_0))=\omega_{\mathbf{M}}\check{\star}\omega_{\mathbf{N}}(h^2t)+\log(C^2P_0Q_0/(M_0N_0)).
\end{align*}
This finishes the first conclusion; the second one dealing with the stronger growth relations is completely analogous.\vspace{6pt}

$(ii)$ Similarly as before, one infers that $\omega_{\mathbf{P}}(t)\le\omega_{\mathbf{M}}(ht)+\log(CP_0/M_0)$ and $\omega_{\mathbf{N}}(t)\le\omega_{\mathbf{Q}}(ht)+\log(CN_0/Q_0)$ for all $t\ge 0$. Then \eqref{widehatformula} yields for any $t>0$:
\begin{align*}
&\omega_{\mathbf{P}}\widehat{\star}\omega_{\mathbf{Q}}(t)=\sup_{s\ge 0}\{\omega_{\mathbf{P}}(s)-\omega_{\mathbf{Q}}(s/t)\}\le\sup_{s\ge 0}\{\omega_{\mathbf{M}}(hs)-\omega_{\mathbf{N}}(s/(ht))\}+\log(C^2P_0N_0/(M_0Q_0))
\\&
=\sup_{u\ge 0}\{\omega_{\mathbf{M}}(u)-\omega_{\mathbf{N}}(u/(h^2t))\}+\log(C^2P_0N_0/(M_0Q_0))=\omega_{\mathbf{M}}\widehat{\star}\omega_{\mathbf{N}}(h^2t)+\log(C^2P_0N_0/(M_0Q_0)).
\end{align*}
Again, the second part is analogous and note that for $t=0$ we even get equality with the value $0$; see Section \ref{generalizedLegendresection}.
\qed\enddemo

We formulate and prove now the main statement in this section and verify the meaning of the conjugate sequence within the weight sequence setting when involving the generalized Legendre envelopes $\check{\star}$ and $\widehat{\star}$.

\begin{theorem}\label{conjsequthm}
Let $\mathbf{M}$ be a weight sequence and $\mathbf{M}^{\ast}$ the corresponding conjugate sequence. Then the following are valid:

\begin{itemize}
\item[$(i)$] Assume that $\lim_{p\rightarrow+\infty}(m_p)^{1/p}=0$ holds and that there exist sequences $\mathbf{P}$, $\mathbf{Q}$ such that $\mathbf{P}\hyperlink{approx}{\approx}\mathbf{M}\hyperlink{approx}{\approx}\mathbf{Q}$, with $\mathbf{P}$ being log-convex whereas $\mathbf{q}$ being log-concave. Then
$$\omega_{\mathbf{M}}\check{\star}\omega_{\mathbf{M}^{\ast}}\hyperlink{sim}{\sim}\id.$$

\item[$(ii)$] Let $\mathbf{M}$ be equivalent to a strong log-convex sequence $\mathbf{S}$ and such that even $\lim_{p\rightarrow+\infty}(m_p)^{1/p}=+\infty$, then $$\omega_{\mathbf{M}}\widehat{\star}\omega_{(\mathbf{M}^{\ast})^{-1}}\hyperlink{sim}{\sim}\id.$$

\item[$(iii)$] Let $\mathbf{M}$ be log-convex, then $$\omega_{\mathbf{M}^2\mathbf{M}^{\ast}}\widehat{\star}\omega_{\mathbf{M}}\hyperlink{sim}{\sim}\id.$$

\item[$(iv)$] Assume that $\lim_{p\rightarrow+\infty}(m_p)^{1/p}=0$ and $\lim_{p\rightarrow+\infty}(M_pm_p)^{1/p}=+\infty$ holds and that there exist sequences $\mathbf{P}$, $\mathbf{Q}$ such that $\mathbf{P}\hyperlink{approx}{\approx}\mathbf{M}\hyperlink{approx}{\approx}\mathbf{Q}$, with $\mathbf{P}$ being log-convex and $\mathbf{q}$ being log-concave. Assume also that $\mathbf{P}\cdot\mathbf{q}$ is log-convex, then
    $$\omega_{\mathbf{M}}\widehat{\star}\omega_{\mathbf{M}^{\ast}}\hyperlink{simc}{\sim_{\mathfrak{c}}}\omega_{\mathbf{M}\cdot\mathbf{m}}=\omega_{\mathbf{G}^1\cdot \mathbf{m}^2}.$$

\item[$(v)$] Let $\mathbf{M}$ be log-convex, $\lim_{p\rightarrow+\infty}(m_p)^{1/p}=0$, $\mathbf{m}$ is log-concave, $\mathbf{G}^1\mathbf{M}^{-2}$ is log-convex and satisfying $\lim_{p\rightarrow+\infty}(p!M_p^{-2})^{1/p}=+\infty$. Then
    $$\omega_{\mathbf{M}^{\star}}\widehat{\star}\omega_{\mathbf{M}}=\omega_{\frac{\mathbf{M}^{\ast}}{\mathbf{M}}}=\omega_{\mathbf{G}^1\mathbf{M}^{-2}}.$$
\end{itemize}
\end{theorem}

\demo{Proof}
$(i)$ First, note that both $\mathbf{P}$ and $\mathbf{Q}^{\ast}$ are log-convex weight sequences by assumption, see $(b)$ and $(d)$ in Section \ref{conjsequsection}, and so \cite[Thm. 5.1]{genLegendreconj} can be applied to $\mathbf{M}\equiv\mathbf{P}$ and $\mathbf{N}\equiv\mathbf{Q}^{\ast}$. Moreover, by relation $\mathbf{M}\cdot\mathbf{M}^{\ast}\equiv\mathbf{G}^1$ we get $\mathbf{P}\cdot\mathbf{Q}^{\ast}=:\widetilde{\mathbf{G}^1}\hyperlink{approx}{\approx}\mathbf{G}^1$ and so $\mathbf{P}_{\iota}\cdot\mathbf{Q}^{\ast}_{\iota}=+\infty$; see $(c)$ in Section \ref{conjsequsection}. Comment $(d)$ there and Section \ref{assofunctionsection} yield that $\omega_{\mathbf{P}}$ and $\omega_{\mathbf{Q}^{\star}}$ are weight functions. Thus, when gathering this information we obtain:
$$\forall\;t\in[0,+\infty):\;\;\;\omega_{\mathbf{P}}\check{\star}\omega_{\mathbf{Q}^{\star}}(t)=\omega_{\mathbf{P}\cdot\mathbf{Q}^{\star}}(t)=\omega_{\widetilde{\mathbf{G}^1}}(t).$$
Next recall that $\omega_{\mathbf{G}^1}\hyperlink{sim}{\sim}\id$ is valid by \cite[Ex. 2.9]{genLegendreconj} (for $s=1$) and that the equivalence between $\widetilde{\mathbf{G}^1}$ and $\mathbf{G}^1$ implies $\omega_{\widetilde{\mathbf{G}^1}}\hyperlink{simc}{\sim_{\mathfrak{c}}}\omega_{\mathbf{G}^1}$; see \eqref{assofctrel}. Moreover, $\omega_{\mathbf{M}}\check{\star}\omega_{\mathbf{M}^{\star}}\hyperlink{simc}{\sim_{\mathfrak{c}}}\omega_{\mathbf{P}}\check{\star}\omega_{\mathbf{Q}^{\star}}$ by taking into account $(e)$ in Section \ref{conjsequsection} and $(i)$ in Lemma \ref{conjsequlemma} and which implies $\omega_{\mathbf{M}}\check{\star}\omega_{\mathbf{M}^{\star}}\hyperlink{simc}{\sim_{\mathfrak{c}}}\omega_{\mathbf{G}^1}$. Finally, since $\omega_{\mathbf{G}^1}$ satisfies \hyperlink{om1}{$(\omega_1)$} by \cite[Prop 2.2 $(ii)$]{ultradifferentiablecomparison} we obtain $\omega_{\mathbf{M}}\check{\star}\omega_{\mathbf{M}^{\star}}\hyperlink{sim}{\sim}\omega_{\mathbf{G}^1}$ and hence the conclusion.\vspace{6pt}

$(ii)$ Apply \cite[Thm. 5.7]{genLegendreconj} to $\mathbf{M}\equiv\mathbf{S}$ and $\mathbf{N}\equiv(\mathbf{S}^{\ast})^{-1}$ and obtain
$$\omega_{\mathbf{S}}\widehat{\star}\omega_{\mathbf{N}}=\omega_{\frac{\mathbf{S}}{\mathbf{N}}}=\omega_{\mathbf{G}^1}\hyperlink{sim}{\sim}\id.$$
Note that in this case the quotient sequence $\frac{\mathbf{S}}{\mathbf{N}}$ is equal to $\mathbf{G}^1$ and thus clearly log-convex. Moreover, in view of \eqref{conjsequdef} we have $N_p=\frac{S_p}{p!}=s_p$ for all $p\in\NN$ and consequently the log-convexity for $\mathbf{N}$ is precisely assumption \emph{strong log-convexity} for $\mathbf{S}$. Therefore, all appearing sequences are log-convex by assumption and definition. $\mathbf{M}_{\iota}=+\infty=\mathbf{N}_{\iota}$ since even $\lim_{p\rightarrow+\infty}(m_p)^{1/p}=+\infty$ holds. Finally, relation $\mathbf{N}\hyperlink{vartriangle}{\vartriangleleft}\mathbf{S}$ amounts to have
$$\forall\;h>0\;\exists\;C_h\ge 1\;\forall\;p\in\NN:\;\;\;\frac{S_p}{p!}\le C_hh^pS_p,$$
which is equivalent to requiring $\lim_{p\rightarrow+\infty}(p!)^{1/p}=+\infty$ and this is clear by \eqref{Stirlinglike}. Thus the result follows by taking into account comment $(e)$ in Section \ref{conjsequsection}, $(ii)$ in Lemma \ref{conjsequlemma} and involving again \cite[Prop 2.2]{ultradifferentiablecomparison}: Since $\omega_{\mathbf{G}^1}$ satisfies both \hyperlink{om1}{$(\omega_1)$} and \hyperlink{om6}{$(\omega_6)$} both equivalence relations \hyperlink{sim}{$\sim$} and \hyperlink{simc}{$\sim_{\mathfrak{c}}$} hold simultaneously.\vspace{6pt}

$(iii)$ Apply \cite[Thm. 5.7]{genLegendreconj} to $\mathbf{M}_1\equiv\mathbf{M}^2\mathbf{M}^{\ast}=\mathbf{G}^1\mathbf{M}$ and $\mathbf{N}_1\equiv\mathbf{M}$ and thus $\frac{\mathbf{M}_1}{\mathbf{N}_1}=\mathbf{G}^1$. All sequences $\mathbf{G}^1\mathbf{M}$, $\mathbf{M}$ and $\mathbf{G}^1$ are clearly log-convex and consequently, analogously as in $(ii)$, we get
$$\omega_{\mathbf{M}^2\mathbf{M}^{\ast}}\widehat{\star}\omega_{\mathbf{M}}=\omega_{\mathbf{G}^1}\hyperlink{sim}{\sim}\id.$$
Note that here $(\mathbf{M}_1)_{\iota}=+\infty=(\mathbf{N}_1)_{\iota}$ holds since $\mathbf{M}$ is a weight sequence and relation $\mathbf{N}_1\hyperlink{vartriangle}{\vartriangleleft}\mathbf{M}_1$ amounts to have
$$\forall\;h>0\;\exists\;C_h\ge 1\;\forall\;p\in\NN:\;\;\;M_p\le C_hh^p\frac{M_pM_pp!}{M_p},$$
and which is again equivalent to $\lim_{p\rightarrow+\infty}(p!)^{1/p}=+\infty$.\vspace{6pt}

$(iv)$ Apply \cite[Thm. 5.7]{genLegendreconj} to $\mathbf{P}$ and $\mathbf{N}\equiv\mathbf{Q}^{\ast}$ and note that, like in $(i)$ before, both are log-convex weight sequences by assumption. Then, by recalling $(e)$ in Section \ref{conjsequsection}, \eqref{conjsequdef}, \eqref{assofctrel} and $(ii)$ in Lemma \ref{conjsequlemma} we obtain $$\omega_{\mathbf{M}}\widehat{\star}\omega_{\mathbf{M}^{\ast}}\hyperlink{simc}{\sim_{\mathfrak{c}}}\omega_{\mathbf{P}}\widehat{\star}\omega_{\mathbf{Q}^{\ast}}=\omega_{\frac{\mathbf{P}}{\mathbf{Q}^{\ast}}}\hyperlink{simc}{\sim_{\mathfrak{c}}}\omega_{\mathbf{M}\cdot\mathbf{m}}=\omega_{\mathbf{G}^1\cdot \mathbf{m}^2}.$$
For this note that $\mathbf{P}\cdot\mathbf{q}\equiv\frac{\mathbf{P}}{\mathbf{Q}^{\ast}}$ is assumed to be log-convex and $\mathbf{Q}^{\ast}\hyperlink{vartriangle}{\vartriangleleft}\mathbf{P}$ if and only if $\mathbf{M}^{\ast}\hyperlink{vartriangle}{\vartriangleleft}\mathbf{M}$ which amounts to have
$$\forall\;h>0\;\exists\;C_h\ge 1\;\forall\;p\in\NN:\;\;\;\frac{p!}{M_p}\le C_hh^pM_p=C_hh^pp!m_p;$$
i.e. requiring $\lim_{p\rightarrow+\infty}(M_pm_p)^{1/p}=+\infty$.\vspace{6pt}

$(v)$ Apply \cite[Thm. 5.7]{genLegendreconj} to $\mathbf{M}\equiv\mathbf{M}^{\ast}$ and $\mathbf{N}\equiv\mathbf{M}$. Again both sequences are log-convex weight sequences by assumption and in view of \eqref{conjsequdef} one has
$$\forall\;t\in[0,+\infty):\;\;\;\omega_{\mathbf{M}^{\ast}}\widehat{\star}\omega_{\mathbf{M}}(t)=\omega_{\frac{\mathbf{M}^{\ast}}{\mathbf{M}}}(t)=\omega_{\mathbf{G}^1\mathbf{M}^{-2}}(t).$$
The log-convexity $\frac{\mathbf{M}^{\ast}}{\mathbf{M}}$ is ensured since $\mathbf{G}^1\mathbf{M}^{-2}$ is log-convex and $\lim_{p\rightarrow+\infty}(p!M^{-2}_p)^{1/p}=+\infty$ precisely gives $\mathbf{M}\hyperlink{vartriangle}{\vartriangleleft}\mathbf{M}^{\ast}$.
\qed\enddemo

Concerning $(iv)$ and $(v)$ in Theorem \ref{conjsequthm} we comment on some explicit example(s); more precisely we consider \emph{small Gevrey-sequences.}

\begin{remark}
\emph{Take $\mathbf{M}\equiv\mathbf{G}^{\alpha}$ with $\alpha\in(0,1)$. Then $\mathbf{m}\equiv\mathbf{G}^{\alpha-1}$, so $\mathbf{M}$ is a log-convex weight sequence and $\mathbf{m}$ is log-concave such that $\lim_{p\rightarrow+\infty}(m_p)^{1/p}=0$. Moreover, $\mathbf{G}^1\cdot\mathbf{m}^2\equiv\mathbf{G}^{2\alpha-1}$ is log-convex and $\lim_{p\rightarrow+\infty}(p!m_p^2)^{1/p}=+\infty$ provided that $\alpha\in(\frac{1}{2},1)$. On the other hand, $\mathbf{G}^1\mathbf{M}^{-2}\equiv\mathbf{G}^{1-2\alpha}$ is log-convex and $\lim_{p\rightarrow+\infty}(p!M_p^{-2})^{1/p}=+\infty$ provided that $\alpha\in(0,\frac{1}{2})$.}

\emph{However, in this setting the index $\alpha=\frac{1}{2}$ has to be excluded: $\mathbf{M}\equiv\mathbf{G}^{1/2}$ yields $\mathbf{M}\equiv\mathbf{M}^{\ast}$ but $\omega_{\mathbf{N}}\widehat{\star}\omega_{\mathbf{N}}$ is not well defined with $\mathbf{N}$ being any log-convex weight sequence (thus $\mathbf{N}_{\iota}=+\infty$); for this recall \cite[Cor. 5.3]{genLegendreconj}.}
\end{remark}

\subsection{Conjugate sequences and conjugate associated weight functions}
Let $\mathbf{M}$ be a log-convex weight sequence satisfying $\lim_{p\rightarrow+\infty}(m_p)^{1/p}=0$. In view of Corollary \ref{conjwelldefcor} we have that both $\omega^{\ast}_{\mathbf{M}}$ and $\omega_{\mathbf{M}^{\ast}}$ are well defined and hence it makes sense to consider two conjugate operations: either the conjugate sequence $\mathbf{M}^{\ast}$ or the conjugate associated weight function $\omega^{\ast}_{\mathbf{M}}$. The aim is to establish a precise connection between both notions, to relate and to compare $\omega^{\ast}_{\mathbf{M}}$ and $\omega_{\mathbf{M}^{\ast}}$ and we prove the following main statement.

\begin{theorem}\label{conjequivthm}
Let $\mathbf{M}$ be a given weight sequence such that $\lim_{p\rightarrow+\infty}(m_p)^{1/p}=0$. Then the following holds:
\begin{itemize}
\item[$(i)$] Assume that there exist sequences $\mathbf{L}$, $\mathbf{N}$ such that $\mathbf{L}\hyperlink{approx}{\approx}\mathbf{M}\hyperlink{approx}{\approx}\mathbf{N}$, $\mathbf{L}$ is log-convex whereas $\mathbf{n}$ is log-concave. Then
$$\exists\;C\ge 1\;\forall\;s\ge 0:\;\;\;\omega^{\ast}_{\mathbf{M}}(s)\le\omega_{\mathbf{M}^{\ast}}(Cs)+1.$$

\item[$(ii)$] Assume that $\mathbf{M}$ is log-convex and that $\left(\frac{\mu_p}{p}\right)_p$ is almost decreasing. Then there exists a log-convex weight sequence $\mathbf{L}$ such that $\mathbf{L}\hyperlink{approx}{\approx}\mathbf{M}$, $\mathbf{l}$ is log-concave and

    $$\forall\;s\in[0,\lambda^{\ast}_1]:\;\;\;\omega_{\mathbf{L}}^{\ast}(s)\ge\omega_{\mathbf{L}^{\ast}}(s)=0,\hspace{15pt}\forall\;q\in\NN_{>0}:\;\;\;\omega_{\mathbf{L}}^{\ast}(\lambda^{\ast}_q)\ge\omega_{\mathbf{L}^{\ast}}(\lambda^{\ast}_q),$$
and
$$\forall\;\lambda^{*}_q<s<\lambda^{*}_{q+1},\;q\in\NN_{>0}:\;\;\;\omega_{\mathbf{L}}^{\ast}(2s)\ge\omega_{\mathbf{L}^{\ast}}(s).$$
\end{itemize}

Summarizing, we infer
\begin{equation}\label{conjequivthmequ0}
\exists\;C\ge 1\;\forall\;s\ge 0:\;\;\;\omega_{\mathbf{L}^{\ast}}(s/2)\le\omega_{\mathbf{L}}^{\ast}(s)\le\omega_{\mathbf{L}^{\ast}}(Cs)+1,
\end{equation}
which gives for $\mathbf{M}$, $\mathbf{M}^{\ast}$ the estimates
\begin{equation}\label{conjequivthmequ}
\exists\;h,D\ge 1\;\forall\;s\ge 0:\;\;\;-D+\omega_{\mathbf{M}^{\ast}}(s/h)\le\omega_{\mathbf{M}}^{\ast}(s)\le\omega_{\mathbf{M}^{\ast}}(hs)+D.
\end{equation}
Therefore, $\omega_{\mathbf{L}^{\ast}}\hyperlink{simc}{\sim_{\mathfrak{c}}}\omega_{\mathbf{L}}^{\ast}$ and $\omega_{\mathbf{M}^{\ast}}\hyperlink{simc}{\sim_{\mathfrak{c}}}\omega_{\mathbf{M}}^{\ast}$ is valid.
\end{theorem}

\demo{Proof}
First, in view of $(d)$ in Section \ref{conjsequsection} the function $\omega_{\mathbf{M}^{\ast}}$ is well defined. Second, for $s=0$ one obtains $\omega^{\ast}_{\mathbf{M}}(0)=\sup_{t\ge 0}\{-\omega_{\mathbf{M}}(t)\}=0$ which immediately gives $\omega^{\ast}_{\mathbf{M}}(0)=0=\omega_{\mathbf{M}^{\ast}}(0)$ and hence equality. So from now on we focus on $s>0$.\vspace{6pt}

$(i)$ By $(i)$ in Theorem \ref{conjsequthm} we have $\omega_{\mathbf{M}}\check{\star}\omega_{\mathbf{M}^{\ast}}\hyperlink{sim}{\sim}\id$ which does imply
$$\exists\;C\ge 1\;\forall\;t\ge 0:\;\;\;-1+\frac{1}{C}t\le\inf_{s>0}\{\omega_{\mathbf{M}}(s)+\omega_{\mathbf{M}^{\ast}}(t/s)\},$$
resp. with $u:=\frac{t}{s}$ that
$$\exists\;C\ge 1\;\forall\;u,s\ge 0:\;\;\;-1+\frac{1}{C}su\le\omega_{\mathbf{M}}(s)+\omega_{\mathbf{M}^{\ast}}(u).$$
Note that this estimate is clear for $s=0$. Now put $v:=\frac{u}{C}$ and so we have shown
$$\exists\;C\ge 1\;\forall\;v,s\ge 0:\;\;\;sv\le\omega_{\mathbf{M}}(s)+\omega_{\mathbf{M}^{\ast}}(Cv)+1,$$
which gives by definition $\omega^{\ast}_{\mathbf{M}}(v)=\sup_{s\ge 0}\{sv-\omega_{\mathbf{M}}(s)\}\le\omega_{\mathbf{M}^{\ast}}(Cv)+1$ for all $v\ge 0$ and therefore the conclusion.\vspace{6pt}

$(ii)$ By Remark \ref{Rem24reloaded} we have that there exists an equivalent and log-convex weight sequence $\mathbf{L}$ such that $\mathbf{l}$ is log-concave and $\lim_{p\rightarrow+\infty}(l_p)^{1/p}=0$ (by equivalence) and so $\mathbf{L}^{\ast}$ is a log-convex weight sequence, too. Moreover, in view of comment $(e)$ in Section \ref{conjsequsection} also $\mathbf{M}^{\ast}$ and $\mathbf{L}^{\ast}$ are equivalent. We start with the proof of the desired estimates:

First, on the one hand $\omega_{\mathbf{L}^{\ast}}(s)=0$ for all $s\in[0,\lambda^{\ast}_1]$ follows by Lemma \ref{lemma1} applied to $\mathbf{L}^{\ast}$ and, on the other hand, since $\omega_{\mathbf{L}}(\lambda_1)=0$ by definition $\omega_{\mathbf{L}}^{\ast}(s)\ge s\lambda_1\ge 0$ for all $s\ge 0$.

Second, let $q\in\NN_{>0}$ and by taking into account definition \eqref{conjsequdef}, Lemma \ref{lemma1} applied to both $\mathbf{L}^{\ast}$ and $\mathbf{L}$, and \eqref{Stirlinglike} (used in the second inequality) we estimate as follows for all $q\in\NN_{>0}$:
\begin{align*}
\omega_{\mathbf{L}}^{\ast}(\lambda^{\ast}_q)&=\omega_{\mathbf{L}}^{\ast}(q/\lambda_q)=\sup_{t\ge 0}\{\frac{q}{\lambda_q}t-\omega_{\mathbf{L}}(t)\}\underbrace{\ge}_{t:=\lambda_q}q-\omega_{\mathbf{L}}(\lambda_q)=q-\log\left(\frac{L_0\lambda_q^q}{L_q}\right)=\log\left(\frac{L_qe^q}{L_0\lambda_q^q}\right)
\\&
\ge\log\left(\frac{L_qq^q}{q!L_0\lambda_q^q}\right)=\log\left(\frac{l_qq^q}{l_0\lambda_q^q}\right)=\log\left(\frac{L_0^{\ast}(\lambda^{\ast}_q)^q}{L^{\ast}_q}\right)=\omega_{\mathbf{L}^{\ast}}(\lambda^{\ast}_q).
\end{align*}

Finally, let $s>0$ be such that $\frac{q}{\lambda_q}<s<\frac{q+1}{\lambda_{q+1}}$ for some $q\in\NN_{>0}$. For each $s>\frac{1}{\lambda_1}$ such $q$ exists since $\lim_{p\rightarrow+\infty}\frac{q}{\lambda_q}=+\infty$ by assumption $\mathbf{L}\hyperlink{mtriangle}{\vartriangleleft}\mathbf{G}^1$ and comment $(d)$ in Section \ref{conjsequsection}. Hence $2s>\frac{2q}{\lambda_q}\ge\frac{q+1}{\lambda_q}\ge\frac{q+1}{\lambda_{q+1}}$ and the last estimate is valid because $\mathbf{L}$ is log-convex. Then, by taking into account the fact that $\omega_{\mathbf{L}^{\ast}}$ and $\omega^{\ast}_{\mathbf{L}}$ are non-decreasing and the estimate shown before we get:
\begin{align*}
\omega_{\mathbf{L}}^{\ast}(2s)&\ge\omega_{\mathbf{L}}^{\ast}\left(\frac{q+1}{\lambda_{q+1}}\right)\ge\omega_{\mathbf{L}^{\ast}}\left(\frac{q+1}{\lambda_{q+1}}\right)\ge\omega_{\mathbf{L}^{\ast}}(s).
\end{align*}
Combining this part with $(i)$ applied to $\mathbf{M}\equiv\mathbf{L}\equiv\mathbf{N}$ we have shown \eqref{conjequivthmequ0} and from this \eqref{conjequivthmequ} follows by taking into account \eqref{assofctrel1} and \eqref{assofctrel2} and recall for this the equivalences between $\mathbf{M}$ and $\mathbf{L}$, $\mathbf{M}^{\ast}$ and $\mathbf{L}^{\ast}$.
\qed\enddemo

We proceed with two immediate consequences of this theorem; concerning the notation in the next corollary we refer to \cite[Sect. 3]{microclasses} and \cite[Sect. 2.5, $(2.14)$]{weightedentireinclusion1}.

\begin{corollary}\label{conjequivthmcor}
Let $\mathbf{M}$ be a log-convex weight sequence such that $\left(\frac{\mu_p}{p}\right)_p$ is almost decreasing and $\lim_{p\rightarrow+\infty}(m_p)^{1/p}=0$. Then in the main result \cite[Thm. 3.4]{microclasses}, instead of the dilatation-type system $v_{\mathbf{M}^{\ast},c}(t):=\exp(-\omega_{\mathbf{M}^{\ast}}(ct))$, $c>0$, one can equivalently involve the (dilatation-type) weight function system
$$v^{\ast}_{\mathbf{M},c}(t):=\exp(-\omega^{\ast}_{\mathbf{M}}(ct)),\;\;\;c>0.$$
\end{corollary}

\demo{Proof}
This follows by \eqref{conjequivthmequ} which transfers to the above defined weight functions in the weighted entire setting as follows:
$$\exists\;h,D\ge 1\;\forall\;c>0\;\forall\;t\ge 0:\;\;\;v_{\mathbf{M}^{\ast},ch}(t)\le e^Dv^{\ast}_{\mathbf{M},c}(t)\le e^{2D}v_{\mathbf{M}^{\ast},c/h}(t).$$
\qed\enddemo

\emph{Note:} In Corollary \ref{conjequivthmcor} the assumptions on $\mathbf{M}$ are stronger compared with \cite[Thm. 3.4 \& Rem. 3.5]{microclasses} since here we require more regularity and growth assumptions on $\mathbf{M}$. Indeed, apart from the inverse normalization condition $1=M_0\ge M_1$ in \cite{microclasses}, which is a minor technical issue, here in addition crucially it is required that $\mathbf{M}$ is a log-convex weight sequence. However, since we are involving $\omega_{\mathbf{M}}$ directly in the definition of the new weight function system, the above Corollary \ref{conjequivthmcor} gives more information and recall: The assumption that $\mathbf{M}$ is a log-convex weight sequence ensures that $\omega_{\mathbf{M}}$ is well defined, see Definition \ref{defweightsequ} and Section \ref{assofunctionsection}, and hence this additional assumption on $\mathbf{M}$ is natural when involving $\omega_{\mathbf{M}}$.

\begin{corollary}\label{conjequivthmcor1}
Let $\mathbf{M}$ be a log-convex weight sequence such that $\lim_{p\rightarrow+\infty}(m_p)^{1/p}=0$ and $\alpha(\mu)<1$. Then:
\begin{itemize}
\item[$(*)$] The assertions listed in Remark \ref{Rem24reloaded} are valid;

\item[$(*)$] Both $\omega_{\mathbf{M}}^{\ast}$ and $\omega_{\mathbf{M}^{\ast}}$ are well defined and satisfy in addition \hyperlink{om1}{$(\omega_1)$};

\item[$(*)$] $\omega_{\mathbf{M}^{\ast}}\hyperlink{sim}{\sim}\omega_{\mathbf{M}}^{\ast}$ is valid;

\item[$(*)$] $\mathbf{M}$ satisfies \hyperlink{mg}{$(\on{mg})$}.
\end{itemize}
\end{corollary}

\demo{Proof}
First, by definition $\alpha(\mu)<1$ obviously implies that $\left(\frac{\mu_p}{p}\right)_p$ is almost decreasing; thus Remark \ref{Rem24reloaded} holds and $(ii)$ in Theorem \ref{conjequivthm} can be used. Second, by applying \cite[Thm. 2.16]{index} to $\omega_{\mathbf{M}}$ and \cite[Cor. 4.6 $(ii)$]{index} to $\mathbf{M}$ we get
\begin{equation}\label{conjequivthmcor1equ1}
\overline{\gamma}(\omega_{\mathbf{M}})=\frac{1}{\beta(\omega_{\mathbf{M}})}=\alpha(\mu)<1,
\end{equation}
where $\beta(\cdot)$ for a given weight function is defined analogously as for a weight sequence; we refer to \cite[Sect. 2.2, Thm. 2.4]{index} and the citations there.

Then $(i)$ in Theorem \ref{indexlemma} applied to $\sigma=\omega_{\mathbf{M}}$ gives $0<1-\overline{\gamma}(\omega_{\mathbf{M}})\le\gamma(\omega_{\mathbf{M}}^{\ast})$ and thus $\omega_{\mathbf{M}}^{\ast}$ has \hyperlink{om1}{$(\omega_1)$}; recall Remark \ref{om1om6indexrem}.

Theorem \ref{conjequivthm} implies $\omega_{\mathbf{M}^{\ast}}\hyperlink{simc}{\sim_{\mathfrak{c}}}\omega_{\mathbf{M}}^{\ast}$ and then $\omega_{\mathbf{M}^{\ast}}\hyperlink{sim}{\sim}\omega_{\mathbf{M}}^{\ast}$ follows by taking into account \cite[Proposition 2.2 $(ii)$]{ultradifferentiablecomparison}. Since \hyperlink{om1}{$(\omega_1)$} is preserved under equivalence also $\omega_{\mathbf{M}^{\ast}}$ satisfies this property. Finally, via \eqref{conjequivthmcor1equ1} we obtain $\overline{\gamma}(\omega_{\mathbf{M}})<1$, which implies by taking into account Remark \ref{om1om6indexrem} that, in particular, $\omega_{\mathbf{M}}$ has to satisfy \hyperlink{om6}{$(\omega_6)$}. And this condition is then equivalent to \hyperlink{mg}{$(\on{mg})$} for $\mathbf{M}$; see \cite[Lemma 2.2]{whitneyextensionweightmatrix}, \cite[Thm. 1]{matsumoto} and \cite{modgrowthstrange}.
\qed\enddemo

\begin{remark}\label{conjequivthmcor1rem}
\emph{Similarly, let $\mathbf{M}\in\RR_{>0}^{\NN}$ be a log-convex weight sequence and assume that $\left(\frac{\mu_p}{p}\right)_p$ is almost decreasing. Then the assertions listed in Remark \ref{Rem24reloaded} are valid, $\alpha(\mu)\le 1$ and via \eqref{conjequivthmcor1equ1} we obtain then $\overline{\gamma}(\omega_{\mathbf{M}})\le 1$. And this implies again that, in particular, $\omega_{\mathbf{M}}$ has to satisfy \hyperlink{om6}{$(\omega_6)$} and \hyperlink{mg}{$(\on{mg})$} holds for $\mathbf{M}$.}
\end{remark}

\section{The Braun-Meise-Taylor weight function setting}\label{BMTsection}
We focus now on a certain and crucial subclass of weight functions, more precisely we consider weight functions in the sense of \emph{Braun-Meise-Taylor (BMT-weights for short);} see \cite{BraunMeiseTaylor90}.

\subsection{Weight functions in the sense of Braun-Meise-Taylor}\label{BMTdefsection}
Let $\omega:[0,+\infty)\rightarrow[0,+\infty)$ be continuous, non-decreasing, $\omega(0)=0$ and $\lim_{t\rightarrow+\infty}\omega(t)=+\infty$. If $\omega$ satisfies in addition $\omega(t)=0$ for all $t\in[0,1]$, then $\omega$ is called \emph{normalized} and this can always be assumed w.l.o.g. For convenience we write that $\omega$ has $\hypertarget{om0}{(\omega_0)}$ if it satisfies all these assumptions; see e.g. \cite[Sect. 2.1]{index} and \cite[Sect. 2.2]{sectorialextensions}.\vspace{6pt}

Apart from \hyperlink{om1}{$(\omega_1)$} and \hyperlink{om6}{$(\omega_6)$} already treated in Section \ref{growthsection} we consider the following conditions; again the abbreviations have already been used in ~\cite{dissertation}.

\begin{itemize}
\item[\hypertarget{om3}{$(\omega_3)$}] $\log(t)=o(\omega(t))$ as $t\rightarrow+\infty$.
	
\item[\hypertarget{om4}{$(\omega_4)$}] $\varphi_{\omega}:t\mapsto\omega(e^t)$ is a convex function on $\RR$.
	
\end{itemize}


In the literature one can find different assumptions for BMT-weights. However, some of the conditions are basic and for convenience we introduce the set
$$\hypertarget{omset0}{\mathcal{W}_0}:=\{\omega:[0,\infty)\rightarrow[0,\infty): \omega\;\text{has}\;\hyperlink{om0}{(\omega_0)},\hyperlink{om3}{(\omega_3)},\hyperlink{om4}{(\omega_4)}\}.$$
In the forthcoming \hyperlink{omset0}{$\mathcal{W}_0$} is understood to be the set of (normalized) weight functions in the sense of Braun-Meise-Taylor. Then, for any $\omega\in\hyperlink{omset0}{\mathcal{W}_0}$ we define the \emph{Legendre-Fenchel-Young-conjugate} of $\varphi_{\omega}$ by
\begin{equation}\label{legendreconjugate}
	\varphi^{*}_{\omega}(x):=\sup\{x y-\varphi_{\omega}(y): y\ge 0\},\;\;\;x\ge 0,
\end{equation}
with the following properties, see e.g. \cite[Rem. 1.3, Lemma 1.5]{BraunMeiseTaylor90} and \cite[Sect. 2.9]{ultradifferentiablecomparison}: It is convex and non-decreasing, $\varphi^{*}_{\omega}(0)=0$, $\varphi^{**}_{\omega}=\varphi_{\omega}$, $\lim_{x\rightarrow+\infty}\frac{x}{\varphi^{*}_{\omega}(x)}=0$ and finally $x\mapsto\frac{\varphi_{\omega}(x)}{x}$ and $x\mapsto\frac{\varphi^{*}_{\omega}(x)}{x}$ are non-decreasing on $[0,+\infty)$. Note that by normalization we can extend the supremum in \eqref{legendreconjugate} from $y\ge 0$ to $y\in\RR$ without changing the value of $\varphi^{*}_{\omega}(x)$ for a given $x\ge 0$.\vspace{6pt}

We mention the following known result, see e.g. \cite[Lemma 2.8]{testfunctioncharacterization} resp. \cite[Lemma 2.4]{sectorialextensions} and the references mentioned in the proofs there.

\begin{lemma}\label{assoweightomega0}
Let $\mathbf{M}\in\hyperlink{LCset}{\mathcal{LC}}$. Then $\omega_{\mathbf{M}}\in\hyperlink{omset0}{\mathcal{W}_0}$. Furthermore, \hyperlink{om6}{$(\omega_6)$} holds for $\omega_{\mathbf{M}}$ if and only if $\mathbf{M}$ satisfies \hyperlink{mg}{$(\on{mg})$}.
\end{lemma}

\subsection{Associated weight matrix}\label{assomatrixsection}
Before recalling some results concerning the associated weight matrix $\mathcal{M}_{\omega}$ we start with the following:
\begin{definition}\label{weightmatrixdef}
A \emph{weight matrix} is the set of sequences
$$\mathcal{M}=\{\mathbf{M}^{(\ell)}: \ell>0\},$$
such that $\mathbf{M}^{(\ell_1)}\le\mathbf{M}^{(\ell_2)}$ for all $0<\ell_1\le\ell_2$. We call $\mathcal{M}$ \emph{standard log-convex} if $\mathbf{M}^{(\ell)}\in\hyperlink{LCset}{\mathcal{LC}}$ for each $\ell>0$ and \emph{constant} if $\mathbf{M}^{(\ell_1)}\hyperlink{approx}{\approx}\mathbf{M}^{(\ell_2)}$ for all $\ell_1,\ell_2>0$.
\end{definition}

Now let us summarize some facts which are shown in \cite[Sect. 5]{compositionpaper} and are needed in this work; all properties listed below are valid for $\omega\in\hyperlink{omset0}{\mathcal{W}_0}$ except \eqref{newexpabsorb} for which \hyperlink{om1}{$(\omega_1)$} is crucial. More basic properties and conditions for abstractly given weight matrices $\mathcal{M}$ can be found e.g. in \cite[Sect. 4]{compositionpaper}.

\begin{itemize}
	\item[$(i)$] The idea was that to each $\omega\in\hyperlink{omset0}{\mathcal{W}_0}$ one can associate a weight matrix $\mathcal{M}_{\omega}:=\{\mathbf{W}^{(\ell)}=(W^{(\ell)}_p)_{p\in\NN}: \ell>0\}$ by\vspace{6pt}
	
	\centerline{$W^{(\ell)}_p:=\exp\left(\frac{1}{\ell}\varphi^{*}_{\omega}(\ell p)\right)$.}\vspace{6pt}
	
	We have that $\mathbf{W}^{(\ell)}\in\hyperlink{LCset}{\mathcal{LC}}$ for each $\ell>0$, i.e. $\mathcal{M}_{\omega}$ is \emph{standard log-convex,} and we even have the stronger order relation
	\begin{equation}\label{quotientorderequ}
		\forall\;\ell_2\ge\ell_1>0:\;\;\;\vartheta^{(\ell_1)}\le\vartheta^{(\ell_2)},
	\end{equation}
	with $\vartheta^{(\ell)}$ denoting the corresponding sequence of quotients; see \cite[Sect. 2.5]{whitneyextensionweightmatrix}. ``Usually'', for abstractly given weight matrices $\mathcal{M}:=\{\mathbf{M}^{(\ell)}: \ell>0\}$ it is sufficient to assume the pointwise order $\mathbf{M}^{(\ell_1)}\le\mathbf{M}^{(\ell_2)}$ for all $0<\ell_1\le\ell_2$; see \cite[Sect. 4]{compositionpaper}.
	
	\item[$(ii)$] $\mathcal{M}_{\omega}$ satisfies
	\begin{equation}\label{newmoderategrowth}
		\forall\;\ell>0\;\forall\;p,q\in\NN:\;\;\;W^{(\ell)}_{p+q}\le W^{(2\ell)}_pW^{(2\ell)}_q.
	\end{equation}
	
	\item[$(iii)$] \hyperlink{om6}{$(\omega_6)$} for $\omega$ holds if and only if some/each $\mathbf{W}^{(\ell)}$ satisfies \hyperlink{mg}{$(\on{mg})$} if and only if $\mathbf{W}^{(\ell_1)}\hyperlink{approx}{\approx}\mathbf{W}^{(\ell_2)}$ for each $\ell_1,\ell_2>0$. Thus \hyperlink{om6}{$(\omega_6)$} is characterizing the case when $\mathcal{M}_{\omega}$ is \emph{constant.}
	
	\item[$(iv)$] In case $\omega$ has in addition \hyperlink{om1}{$(\omega_1)$}, then $\mathcal{M}_{\omega}$ also satisfies
	\begin{equation}\label{newexpabsorb}
		\forall\;h\ge 1\;\exists\;d\ge 1\;\forall\;\ell>0\;\exists\;D\ge 1\;\forall\;p\in\NN:\;\;\;h^pW^{(\ell)}_p\le D W^{(d\ell)}_p.
	\end{equation}
 This estimate is crucial for verifying the equality of ultradifferentiable classes $\mathcal{E}_{[\mathcal{M}_{\omega}]}=\mathcal{E}_{[\omega]}$ (as locally convex vector spaces). Concerning the definition of these spaces we refer to \cite{BraunMeiseTaylor90} and \cite{compositionpaper} and the notation $[\cdot]$ denotes the convention meaning either the \emph{Roumieu case} $\{\cdot\}$ or the \emph{Beurling case} $(\cdot)$.
	
	\item[$(v)$] We have $\omega\hyperlink{sim}{\sim}\omega_{\mathbf{W}^{(\ell)}}$ for each $\ell>0$, more precisely
	\begin{equation}\label{goodequivalenceclassic}
		\forall\;\ell>0\,\,\exists\,D_{\ell}>0\;\forall\;t\ge 0:\;\;\;\ell\omega_{\mathbf{W}^{(\ell)}}(t)\le\omega(t)\le 2\ell\omega_{\mathbf{W}^{(\ell)}}(t)+D_{\ell};
	\end{equation}
	for a proof see \cite[Theorem 4.0.3, Lemma 5.1.3]{dissertation}, \cite[Lemma 5.7]{compositionpaper} and also \cite[Lemma 2.5]{sectorialextensions}. Note that, on the one hand, for proving \eqref{goodequivalenceclassic} the convexity condition \hyperlink{om4}{$(\omega_4)$} is indispensable but, on the other hand, \hyperlink{om4}{$(\omega_4)$} is only required for the second estimate in \eqref{goodequivalenceclassic}.
\end{itemize}

We finish this section by recalling some comments summarized in detail in \cite[Sect. 2.9]{ultradifferentiablecomparison}.

\begin{remark}\label{matrixadmissiblerem}
\emph{By inspecting the proofs in \cite{dissertation} and \cite{compositionpaper} it turns out that for any non-decreasing $\omega$ satisfying \hyperlink{om3}{$(\omega_3)$} it makes sense to associate the well-defined matrix $\mathcal{M}_{\omega}$; we refer to \cite[Sect. 2.9]{ultradifferentiablecomparison} and also to the comments in \cite[Sect. 3.1]{dissertation}, \cite[Sect. 2.7, $(a)-(c)$]{genLegendreconjBMT} and \cite[Sect. 6 $(III)$]{modgrowthstrangeII}. And via \hyperlink{om1}{$(\omega_1)$} and \eqref{newexpabsorb} also the equality between the corresponding weighted spaces holds when the classes are introduced directly (by definition) via weighting the derivatives in terms of $\varphi^{*}_{\omega}$. Indeed, note that the proof of \eqref{newexpabsorb} does not require necessarily \hyperlink{om4}{$(\omega_4)$} but this condition is crucial for showing \eqref{goodequivalenceclassic} (for obtaining the second estimate). Moreover, \hyperlink{om4}{$(\omega_4)$} is needed for the equivalences listed in $(iii)$. Inspired by these comments any non-decreasing function $\omega$ satisfying \hyperlink{om3}{$(\omega_3)$} is called \emph{matrix admissible} and the equivalence between matrix admissible functions transfers in a precise way to the corresponding associated weight matrices (yielding equivalent weight matrices) and which ensures the equality of the corresponding (matrix) weighted spaces; see the relations introduced in \cite[Sect. 4.2, p. 111]{compositionpaper} and the estimate \cite[Sect. 6 $(III)$, $(6.6)$]{modgrowthstrangeII}.}
\end{remark}


\subsection{Conjugate associated weight matrices}\label{conjassomatrixsection}
We investigate the effects of taking the conjugate operations from the previous sections in the BMT-weight function setting: Let $\omega$ be a given BMT-weight function, then consider $\omega^{\ast}$, $\omega^{\ast}_{\mathbf{W}^{(\ell)}}$ according to \eqref{conjugate} and the conjugate weight sequences $(\mathbf{W}^{(\ell)})^{\ast}$ w.r.t. \eqref{conjsequdef}. Our aim now is to study the correspondences of these notions in detail. First we start with some new definitions and observations which are based on these comments.

In this section let us consider $\omega\in\hyperlink{omset0}{\mathcal{W}_0}$ with associated weight matrix $\mathcal{M}_{\omega}=\{\mathbf{W}^{(\ell)}: \ell>0\}$. Since any such weight satisfies by definition $\omega(0)=0$ we always have $(\mathfrak{C}_1)$.

In order to deal with $\omega^{\ast}$ we require that this function is well defined and hence we verify now when $(\mathfrak{C}_2)$ is valid.

\begin{lemma}\label{firstBMTlemma}
Let $\omega\in\hyperlink{omset0}{\mathcal{W}_0}$ be given with associated weight matrix $\mathcal{M}_{\omega}=\{\mathbf{W}^{(\ell)}: \ell>0\}$. Then the following are equivalent:
\begin{itemize}
\item[$(i)$] $\omega$ satisfies $(\mathfrak{C}_2)$; i.e. $\omega^{\ast}$ is a well defined weight function.

\item[$(ii)$] Some/each $\omega^{\ast}_{\mathbf{W}^{(\ell)}}$ is well defined.

\item[$(iii)$] Some/each $\omega_{\mathbf{W}^{(\ell)}}$ satisfies $(\mathfrak{C}_2)$.

\item[$(iv)$] Some/each $\omega_{(\mathbf{W}^{(\ell)})^{\ast}}$ is well defined.

\item[$(v)$] For some/each $\ell>0$ one has the relation $\mathbf{W}^{(\ell)}\hyperlink{mtriangle}{\vartriangleleft}\mathbf{G}^1$; i.e.
\begin{equation}\label{firstBMTlemmaequ}
    \exists\;\ell>0/\forall\;\ell>0:\;\;\;\lim_{p\rightarrow+\infty}\left(\frac{W^{(\ell)}_p}{p!}\right)^{1/p}=\lim_{p\rightarrow+\infty}(w^{(\ell)}_p)^{1/p}=0.
\end{equation}
\end{itemize}
\end{lemma}

\demo{Proof}
$(i)\Leftrightarrow(iii)$ holds by \eqref{goodequivalenceclassic} since $(\mathfrak{C}_2)$ is obviously preserved under equivalence of weight functions. And $(ii)\Leftrightarrow(iii)\Leftrightarrow(iv)\Leftrightarrow(v)$ follows from Lemma \ref{C2assoweightfctlemma} and Corollary \ref{conjwelldefcor} applied to some/any sequence $\mathbf{W}^{(\ell)}$.
\qed\enddemo

Let $\omega\in\hyperlink{omset0}{\mathcal{W}_0}$ be given with associated weight matrix $\mathcal{M}_{\omega}=\{\mathbf{W}^{(\ell)}: \ell>0\}$, then in view of Lemmas \ref{indexlemma0}, \ref{C2assoweightfctlemma} and \ref{firstBMTlemma}, and Corollary \ref{conjwelldefcor} we summarize:

\begin{itemize}
\item[$(a)$] Assume that $t=o(\omega(t))$ as $t\rightarrow+\infty$ holds.

Then $\omega^{\ast}$, $\omega_{\mathbf{W}^{(\ell)}}^{\ast}$ and also $\omega_{(\mathbf{W}^{(\ell)})^{\ast}}$ are well defined for any $\ell>0$ and thus are at least weight functions in the sense of Definition \ref{weightfctdef}. However, in this situation $\omega^{\ast}$ satisfies \hyperlink{om6}{$(\omega_6)$}.

\item[$(b)$] Let us consider the \emph{conjugate associated weight matrix}
\begin{equation}\label{conjumatrix}
\mathcal{M}^{\ast}_{\omega}:=\{(\mathbf{W}^{(1/\ell)})^{\ast}: \ell>0\},
\end{equation}
and note: Comment $(e)$ in Section \ref{conjsequsection} implies that $(\mathbf{W}^{(\ell_2)})^{\ast}\le(\mathbf{W}^{(\ell_1)})^{\ast}$ for all $0<\ell_1\le\ell_2$ and hence to make this fit with Definition \ref{weightmatrixdef} we have to invert the matrix parameter for ensuring the monotonic pointwise order $(\mathbf{W}^{(1/\ell_1)})^{\ast}\le(\mathbf{W}^{(1/\ell_2)})^{\ast}$ for all $0<\ell_1\le\ell_2$.

\item[$(c)$] Comment $(e)$ in Section \ref{conjsequsection} also gives that $\mathcal{M}_{\omega}$ is constant if and only if $\mathcal{M}^{\ast}_{\omega}$ is so. Since each $\mathbf{W}^{(\ell)}\in\hyperlink{LCset}{\mathcal{LC}}$ one has that each $(\mathbf{W}^{(\ell)})^{\ast}$ satisfies \hyperlink{mg}{$(\on{mg})$}; see comment $(f)$ in Section \ref{conjsequsection}.
\end{itemize}

We are interested in the question whether $\mathcal{M}^{\ast}_{\omega}$ is associated with an matrix admissible (BMT) weight function since this means that taking the conjugate operation (on the matrix level) ``preserves the weight function setting'' and one might expect the identity $\mathcal{M}^{\ast}_{\omega}=\mathcal{M}_{\omega^{\ast}}$.

\begin{proposition}\label{matrixadmissibleprop}
We get the following:
\begin{itemize}
\item[$(i)$] Let $\omega$ be a weight function (in the sense of Definition \ref{weightfctdef}) and assume that $(\mathfrak{C}_1)$ and $(\mathfrak{C}_2)$ holds. Then $\omega^{\ast}$ is matrix admissible and hence $\mathcal{M}_{\omega^{\ast}}$ is well defined.

\item[$(ii)$] Let $\omega\in\hyperlink{omset0}{\mathcal{W}_0}$ with associated weight matrix $\mathcal{M}_{\omega}$ and assume that $(\mathfrak{C}_2)$ holds. Moreover, let $\mathcal{M}^{\ast}_{\omega}$ be the corresponding conjugate matrix from \eqref{conjumatrix}.
\begin{itemize}
\item[$(a)$] Assume that there exist $\sigma\in\hyperlink{omset0}{\mathcal{W}_0}$ with associated weight matrix $\mathcal{M}_{\sigma}=\{\mathbf{S}^{(\ell)}: \ell>0\}$ and $\ell_1,\ell_2>0$ such that $(\mathbf{W}^{(1/\ell_1)})^{\ast}\hyperlink{approx}{\approx}\mathbf{S}^{(\ell_2)}$. Then $\sigma$ has to satisfy \hyperlink{om6}{$(\omega_6)$} and $\mathcal{M}_{\sigma}$ is constant.

\item[$(b)$] Assume that there exists $\sigma\in\hyperlink{omset0}{\mathcal{W}_0}$ such that $\mathcal{M}^{\ast}_{\omega}$ is associated with $\sigma$. Then (also) $\omega$ has to satisfy \hyperlink{om6}{$(\omega_6)$} and $\mathcal{M}_{\omega}$ is constant.
\end{itemize}
\end{itemize}
\end{proposition}

\demo{Proof}
$(i)$ By assumption $\omega^{\ast}$ is a well defined weight function and \hyperlink{om3}{$(\omega_3)$} follows by taking into account \eqref{conjugatefastgrowthequ}.

$(ii)(a)$ Each $(\mathbf{W}^{(\ell)})^{\ast}$ satisfies \hyperlink{mg}{$(\on{mg})$} as mentioned in $(c)$ before and hence $\mathbf{S}^{(\ell_2)}$ has \hyperlink{mg}{$(\on{mg})$}, too. Then the assertions follow by $(iii)$ in Section \ref{assomatrixsection}.

$(ii)(b)$ By $(ii)(a)$ the matrix $\mathcal{S}\equiv\mathcal{M}^{\ast}_{\omega}$ has to be constant, so $\mathcal{M}_{\omega}$ is constant too and this fact implies \hyperlink{om6}{$(\omega_6)$} for $\omega$; again see $(iii)$ in Section \ref{assomatrixsection}.
\qed\enddemo

However, note that \hyperlink{om6}{$(\omega_6)$} is a (too) restrictive assumption for general BMT-weight functions in the sense that then $\mathcal{M}_{\omega}$ has to be constant and the corresponding weighted spaces can already be described by a single weight sequence (take any arbitrary $\mathbf{W}^{(\ell)}$); for this recall the comparison results from \cite{BonetMeiseMelikhov07} and the statements from \cite[Sect. 5]{compositionpaper}.

\hyperlink{om6}{$(\omega_6)$} is equivalent to having $\overline{\gamma}(\omega)<+\infty$; if we assume even a stronger resp. more specific bound on this growth index then we can show more crucial properties and so we formulate and prove the first main result in this section.

\begin{theorem}\label{mainBMTthm}
Let $\omega\in\hyperlink{omset0}{\mathcal{W}_0}$ be given such that $t=o(\omega(t))$ as $t\rightarrow+\infty$ (i.e. $(\mathfrak{C}_2)$). Let $\mathcal{M}_{\omega}$ be the associated weight matrix and $\mathcal{M}^{\ast}_{\omega}$ the corresponding conjugate matrix.
\begin{itemize}
\item[$(i)$] One has
\begin{equation}\label{mainBMTthmequ}
\forall\;\ell>0\;\exists\;D_{\ell}\ge 1\;\forall\;s\ge 0:\;\;\;\frac{1}{\ell}\omega^{\ast}(s\ell)\le\omega_{\mathbf{W}^{(\ell)}}^{\ast}(s)\le\frac{1}{2\ell}\omega^{\ast}(2\ell s)+\frac{D_{\ell}}{2\ell},
\end{equation}
with $D_{\ell}$ being the constants from \eqref{goodequivalenceclassic}.

\item[$(ii)$] Assume in addition that $\overline{\gamma}(\omega)<1$, then:
\begin{itemize}
\item[$(a)$] $\omega^{\ast}$ satisfies \hyperlink{om1}{$(\omega_1)$}.

\item[$(b)$] $\omega^{\ast}\hyperlink{sim}{\sim}\omega_{\mathbf{W}^{(\ell)}}^{\ast}$ for any $\ell>0$ and each $\omega_{\mathbf{W}^{(\ell)}}^{\ast}$ is well defined and satisfies \hyperlink{om1}{$(\omega_1)$}.

\item[$(c)$] Each $\omega_{(\mathbf{W}^{(\ell)})^{\ast}}$ is well defined, satisfies \hyperlink{om1}{$(\omega_1)$} and the following equivalences hold:
\begin{equation}\label{mainBMTthmequ1}
    \forall\;\ell>0:\;\;\;\omega^{\ast}\hyperlink{sim}{\sim}\omega_{\mathbf{W}^{(\ell)}}^{\ast}\hyperlink{sim}{\sim}\omega_{(\mathbf{W}^{(\ell)})^{\ast}}.
\end{equation}
\item[$(d)$] The matrices $\mathcal{M}_{\omega}$ and $\mathcal{M}^{\ast}_{\omega}$ are constant and all sequences $\mathbf{W}^{(\ell)}$ and $(\mathbf{W}^{(\ell)})^{\ast}$ have \hyperlink{mg}{$(\on{mg})$}.

\item[$(e)$] For each $\ell>0$ there exists $\mathbf{V}^{(\ell)}$ such that $\mathbf{V}^{(\ell)}\hyperlink{approx}{\approx}\mathbf{W}^{(\ell)}$, $\mathbf{V}^{(\ell)}$ is log-convex and $\mathbf{v}^{(\ell)}$ is log-concave. Indeed, one can even achieve that $\mathbf{V}^{(\ell)},(\mathbf{V}^{(\ell)})^{\ast}\in\hyperlink{LCset}{\mathcal{LC}}$ for each $\ell>0$.

\item[$(f)$] If $\omega$ satisfies in addition \hyperlink{om1}{$(\omega_1)$} then $\overline{\gamma}(\omega^{\ast})<1$ and hence, in particular, $\omega^{\ast}$ has \hyperlink{om6}{$(\omega_6)$}.
\end{itemize}

\item[$(iii)$] Both matrices $\mathcal{M}_{\omega}$ and $\mathcal{M}^{\ast}_{\omega}$ are constant (and thus, in particular, $\omega$ satisfies \hyperlink{om6}{$(\omega_6)$}) if one assumes in addition at least one of the following assertions:
\begin{itemize}
    \item[$(a)$] There exists some $\ell_0>0$ such that $\left(\frac{\vartheta^{(\ell_0)}_p}{p}\right)_p$ is almost decreasing.

    \item[$(b)$] $\omega$ satisfies \hyperlink{om1}{$(\omega_1)$} and there exists some $\ell_0>0$ and a log-convex weight sequence $\mathbf{N}$ such that $\mathbf{W}^{(\ell_0)}\hyperlink{approx}{\approx}\mathbf{N}$ and such that $\left(\frac{\nu_p}{p}\right)_p$ is almost decreasing.

    \item[$(c)$] There exists some $\ell_0>0$ and a log-convex sequence $\mathbf{N}$ such that $\mathbf{N}\hyperlink{approx}{\approx}\mathbf{W}^{(\ell_0)}$ and $p\mapsto(n_p)^{1/p}$ is non-increasing.
\end{itemize}
All assertions $(a)-(c)$ imply $\overline{\gamma}(\omega)\le 1$ and, moreover, $(c)$ yields $(a)$.
\end{itemize}
\end{theorem}

When we put $\mathcal{V}:=\{\mathbf{V}^{(\ell)}: \ell>0\}$, then $(ii)(e)$ shows that in general the pointwise growth order $\mathbf{V}^{(\ell_1)}\le\mathbf{V}^{(\ell_2)}$ might fail in general but the corresponding weighted spaces given by $\mathcal{V}$ and $\mathcal{M}_{\omega}$ coincide and apart from this technical detail $\mathcal{V}$ is standard log-convex.

\demo{Proof}
$(i)$ We use \eqref{goodequivalenceclassic} and follow the proof of $(i)$ in Lemma \ref{equivalencelemma}: This gives $\omega^{\ast}(s)\ge 2\ell\omega_{\mathbf{W}^{(\ell)}}^{\ast}(\frac{s}{2\ell})-D_{\ell}$ and $\omega_{\mathbf{W}^{(\ell)}}^{\ast}(s)\ge\frac{1}{\ell}\omega^{\ast}(s\ell)$ for all $s\ge 0$ and $\ell>0$.\vspace{6pt}

$(ii)(a)$ By $(i)$ in Theorem \ref{indexlemma} one has $0<1-\overline{\gamma}(\omega)\le\gamma(\omega^{\ast})$ and so $\omega^{\ast}$ satisfies \hyperlink{om1}{$(\omega_1)$}; recall Remark \ref{om1om6indexrem}.\vspace{6pt}

$(ii)(b)$ Let $\ell>0$ be arbitrary but fixed. Then the first estimate in \eqref{mainBMTthmequ} gives $\ell\omega_{\mathbf{W}^{(\ell)}}^{\ast}(s/\ell)\ge\omega^{\ast}(s)$ for all $s\ge 0$. Thus, for any $\ell\ge 1$ the relation $\omega^{\ast}(s)=O(\omega_{\mathbf{W}^{(\ell)}}^{\ast}(s))$ as $s\rightarrow+\infty$ is clear. If $0<\ell<1$, then choose $n\in\NN_{>0}$ such that $\frac{1}{\ell}\le 2^n$ holds and iterating \hyperlink{om1}{$(\omega_1)$} for $\omega^{\ast}$ gives
$$\omega^{\ast}(s/\ell)\le\omega^{\ast}(2^ns)\le L\omega^{\ast}(s)+L\le L\ell\omega_{\mathbf{W}^{(\ell)}}^{\ast}(s/\ell)+L,$$
for some $L\ge 1$ (depending on $\ell$ via $n$) and all $s\ge 0$ and so $\omega^{\ast}(s)=O(\omega_{\mathbf{W}^{(\ell)}}^{\ast}(s))$, too.

On the other hand, the second estimate in \eqref{mainBMTthmequ} implies $\omega_{\mathbf{W}^{(\ell)}}^{\ast}(s)=O(\omega^{\ast}(s))$ as $s\rightarrow+\infty$ for all $0<\ell\le\frac{1}{2}$. If $\ell>\frac{1}{2}$, then we choose $n\in\NN_{>0}$ to ensure $2\ell\le 2^n\Leftrightarrow\ell\le 2^{n-1}$, iterate \hyperlink{om1}{$(\omega_1)$} and estimate as follows for some $L\ge 1$ and all $s\ge 0$:
$$\omega_{\mathbf{W}^{(\ell)}}^{\ast}(s)\le\frac{1}{2\ell}\omega^{\ast}(2^ns)+\frac{D_{\ell}}{2\ell}\le\frac{L}{2\ell}\omega^{\ast}(s)+\frac{L}{2\ell}+\frac{D_{\ell}}{2\ell};$$
thus also in this case $\omega_{\mathbf{W}^{(\ell)}}^{\ast}(s)=O(\omega^{\ast}(s))$.

Summarizing, the equivalence between $\omega^{\ast}$ and each $\omega_{\mathbf{W}^{(\ell)}}^{\ast}$ is verified and this implies the rest since $(\mathfrak{C}_2)$ and the growth index $\gamma(\cdot)$ resp. \hyperlink{om1}{$(\omega_1)$} are preserved under equivalence.\vspace{6pt}

$(ii)(c)$ It remains to prove the second equivalence and we apply Corollary \ref{conjequivthmcor1} to $\mathbf{W}^{(\ell)}$, $\ell>0$ arbitrary. First, by assumption $t=o(\omega(t))$ and Lemma \ref{firstBMTlemma} one has that each $\mathbf{W}^{(\ell)}$ is a log-convex weight sequence such that $\lim_{p\rightarrow+\infty}(w^{(\ell)}_p)^{1/p}=0$. Second, by \eqref{goodequivalenceclassic} and the fact that the index $\overline{\gamma}(\cdot)$ is preserved under equivalence we have
\begin{equation}\label{mainBMTthmproofequ}
\forall\;\ell>0:\;\;\;\alpha(\vartheta^{(\ell)})=\overline{\gamma}(\omega_{\mathbf{W}^{(\ell)}})=\overline{\gamma}(\omega)<1;
\end{equation}
for the first equality recall \eqref{conjequivthmcor1equ1}. Then Corollary \ref{conjequivthmcor1} applied to $\mathbf{W}^{(\ell)}$ yields the conclusion.\vspace{6pt}

$(ii)(d)$ By Remark \ref{om1om6indexrem} assumption $\overline{\gamma}(\omega)<1$ implies, in particular, that $\omega$ satisfies \hyperlink{om6}{$(\omega_6)$}. Thus $\mathcal{M}_{\omega}$ is constant, recall $(iii)$ in Section \ref{assomatrixsection}, and each $\mathbf{W}^{(\ell)}$ has \hyperlink{mg}{$(\on{mg})$} (recall also Corollary \ref{conjequivthmcor1}). Comment $(e)$ in Section \ref{conjsequsection} yields that $\mathcal{M}^{\ast}_{\omega}$ is constant too and each $(\mathbf{W}^{(\ell)})^{\ast}$ has \hyperlink{mg}{$(\on{mg})$} by $(f)$ in Section \ref{conjsequsection}.\vspace{6pt}

$(ii)(e)$ \eqref{mainBMTthmproofequ} implies that each $\left(\frac{\vartheta^{(\ell)}_p}{p}\right)_p$ is almost decreasing and then apply Remark \ref{Rem24reloaded} to $\mathbf{W}^{(\ell)}$.\vspace{6pt}

$(ii)(f)$ Remark \ref{om1om6indexrem} and the assumptions give $0<\gamma(\omega)\le\overline{\gamma}(\omega)<1$. Then $(ii)$ in Theorem \ref{indexlemma} implies $\overline{\gamma}(\omega^{\ast})\le 1-\gamma(\omega)<1$, hence \hyperlink{om6}{$(\omega_6)$} for $\omega^{\ast}$ is valid again via Remark \ref{om1om6indexrem}.\vspace{6pt}

$(iii)(a)$ We apply Remark \ref{conjequivthmcor1rem} to $\mathbf{W}^{(\ell_0)}$ and get $\overline{\gamma}(\omega)=\overline{\gamma}(\omega_{\mathbf{W^{(\ell_0)}}})=\alpha(\vartheta^{(\ell_0)})\le 1$. Then \hyperlink{om6}{$(\omega_6)$} follows for $\omega_{\mathbf{W^{(\ell_0)}}}$ and by \eqref{goodequivalenceclassic} also $\omega$ (and each $\omega_{\mathbf{W^{(\ell)}}}$) satisfies \hyperlink{om6}{$(\omega_6)$}. This implies that $\mathcal{M}_{\omega}$ is constant, see $(iii)$ in Section \ref{assomatrixsection} and consequently $\mathcal{M}^{\ast}_{\omega}$ is constant as well.\vspace{6pt}

$(iii)(b)$ By equivalence we get $\omega_{\mathbf{W^{(\ell_0)}}}\hyperlink{simc}{\sim_{\mathfrak{c}}}\omega_{\mathbf{N}}$, see \eqref{assofctrel}, and via \hyperlink{om1}{$(\omega_1)$} for $\omega$ and \eqref{goodequivalenceclassic} we get \hyperlink{om1}{$(\omega_1)$} for $\omega_{\mathbf{W^{(\ell_0)}}}$ too. Then $\omega_{\mathbf{W^{(\ell_0)}}}\hyperlink{sim}{\sim}\omega_{\mathbf{N}}$ follows, see again \cite[Prop. 2.2 $(ii)$]{ultradifferentiablecomparison}. So $(iii)(a)$ applied to $\mathbf{N}$ gives $\overline{\gamma}(\omega)=\overline{\gamma}(\omega_{\mathbf{W^{(\ell_0)}}})=\overline{\gamma}(\omega_{\mathbf{N}})\le 1$ and the rest follows as before.\vspace{6pt}

$(iii)(c)$ Since $\mathbf{N}$ is log-convex, by Lemma \ref{rootalmostdecrlemma} and Remark \ref{rootalmostdecrlemmarem} it follows that $\mathbf{N}$ satisfies \hyperlink{mg}{$(\on{mg})$} and since this condition is preserved under equivalence the sequence $\mathbf{W}^{(\ell_0)}$ has \hyperlink{mg}{$(\on{mg})$}, too. This implies via $(iii)$ in Section \ref{assomatrixsection} that the matrix $\mathcal{M}_{\omega}$ is constant and hence $\mathcal{M}^{\ast}_{\omega}$ is constant, too. Finally, concerning the growth index note: $\mathbf{N}\hyperlink{approx}{\approx}\mathbf{W}^{(\ell_0)}$ holds if and only if $\mathbf{n}\hyperlink{approx}{\approx}\mathbf{w}^{(\ell_0)}$ and so there exists $A\ge 1$ such that for any $1\le p\le q$ one has $(w^{(\ell_0)}_q)^{1/q}\le A(n_q)^{1/q}\le A(n_p)^{1/p}\le A^2(w^{(\ell_0)}_p)^{1/p}$. Since $\mathbf{W}^{(\ell_0)}$ satisfies \hyperlink{mg}{$(\on{mg})$} the estimate $\vartheta^{(\ell_0)}_q\le B(W^{(\ell_0)}_q)^{1/q}=B(q!w^{(\ell_0)}_q)^{1/q}\le Bq(w^{(\ell_0)}_q)^{1/q}$ is valid for some $B\ge 1$ and all $q\in\NN_{>0}$ (recall again Remark \ref{rootalmostdecrlemmarem}). Moreover, $\frac{p}{e}(w^{(\ell_0)}_p)^{1/p}\le(p!w^{(\ell_0)}_p)^{1/p}=(W^{(\ell_0)}_p)^{1/p}\le\vartheta^{(\ell_0)}_p$ holds for all $p\in\NN_{>0}$ by \eqref{Stirlinglike}, by log-convexity and since $W^{(\ell_0)}_0=1$. Summarizing, $\frac{\vartheta^{(\ell_0)}_q}{q}\le A^2Be\frac{\vartheta^{(\ell_0)}_p}{p}$ for all $1\le p\le q$ is shown and since both constants $A$ and $B$ are not depending on $p,q$ the sequence $\left(\frac{\vartheta^{(\ell_0)}_p}{p}\right)_p$ is almost decreasing and the rest follows as in $(iii)(a)$.
\qed\enddemo

For the sake of completeness we remark: Let $\omega$ be given as in Theorem \ref{mainBMTthm} and assume that in addition $0<\gamma(\omega)<1$. Then by $(ii)$ in Theorem \ref{indexlemma} we get $\overline{\gamma}(\omega^{\ast})\le 1-\gamma(\omega)<1$ and hence, in particular, that $\omega^{\ast}$ satisfies \hyperlink{om6}{$(\omega_6)$}. By \eqref{mainBMTthmequ} and \cite[Prop. 2.2 $(i)$]{ultradifferentiablecomparison} it follows that $\omega^{\ast}\hyperlink{simc}{\sim_{\mathfrak{c}}}\omega_{\mathbf{W}^{(\ell)}}^{\ast}$ for all $\ell>0$.

We are now in position to prove the second main result.

\begin{theorem}\label{mainBMTthm1}
Let $\omega\in\hyperlink{omset0}{\mathcal{W}_0}$ be given such that $t=o(\omega(t))$ as $t\rightarrow+\infty$ (i.e. $(\mathfrak{C}_2)$) and that $\overline{\gamma}(\omega)<1$ holds. Let $\mathcal{M}_{\omega}$ be the associated weight matrix and $\mathcal{M}^{\ast}_{\omega}$ the corresponding conjugate matrix.

Then as l.c.v.s.
\begin{equation}\label{mainBMTthm1equ}
\forall\;\ell>0:\;\;\;\mathcal{E}_{[\omega^{\ast}]}=\mathcal{E}_{[\omega_{\mathbf{W}^{(\ell)}}^{\ast}]}=\mathcal{E}_{[\omega_{(\mathbf{W}^{(\ell)})^{\ast}}]}=\mathcal{E}_{[\omega_{(\mathbf{V}^{(\ell)})^{\ast}}]}=\mathcal{E}_{[\mathcal{V}^{\ast}]}=\mathcal{E}_{[\mathcal{M}^{\ast}_{\omega}]}=\mathcal{E}_{[(\mathbf{W}^{(\ell)})^{\ast}]},
\end{equation}
with $\mathcal{V}^{\ast}:=\{(\mathbf{V}^{(1/\ell)})^{\ast}: \ell>0\}$ and $\mathbf{V}^{(\ell)}$ denoting the sequence from $(ii)(e)$ in Theorem \ref{mainBMTthm}.
\end{theorem}

\demo{Proof}
The first two equalities in \eqref{mainBMTthm1equ} follow from $(ii)(c)$ in Theorem \ref{mainBMTthm} (via \eqref{mainBMTthmequ1}), Remark \ref{matrixadmissiblerem} and $(i)$ in Proposition \ref{matrixadmissibleprop}: Note that $\omega^{\ast}$ and each $\omega_{\mathbf{W}^{(\ell)}}^{\ast}$ are automatically matrix admissible and the same holds for each $\omega_{(\mathbf{W}^{(\ell)})^{\ast}}$ since associated weight functions have \hyperlink{om3}{$(\omega_3)$}, recall Lemma \ref{assoweightomega0}. (Alternatively, by \eqref{mainBMTthmequ1} condition \hyperlink{om3}{$(\omega_3)$} can be transferred between all appearing weight functions.)\vspace{6pt}

Concerning the third equality, first note that by $(ii)(e)$ in Theorem \ref{mainBMTthm} we get $\mathbf{V}^{(\ell)}\hyperlink{approx}{\approx}\mathbf{W}^{(\ell)}$, equivalently $(\mathbf{V}^{(\ell)})^{\ast}\hyperlink{approx}{\approx}(\mathbf{W}^{(\ell)})^{\ast}$, for all $\ell>0$. Then \eqref{assofctrel} yields $\omega_{(\mathbf{V}^{(\ell)})^{\ast}}\hyperlink{simc}{\sim_{\mathfrak{c}}}\omega_{(\mathbf{W}^{(\ell)})^{\ast}}$ for any $\ell>0$. Finally, since (each) $\omega_{(\mathbf{W}^{(\ell)})^{\ast}}$ satisfies \hyperlink{om1}{$(\omega_1)$}, see $(ii)(c)$ in Theorem \ref{mainBMTthm}, by using \cite[Prop. 2.2 $(ii)$]{ultradifferentiablecomparison} this relation implies $\omega_{(\mathbf{V}^{(\ell)})^{\ast}}\hyperlink{sim}{\sim}\omega_{(\mathbf{W}^{(\ell)})^{\ast}}$ for any $\ell>0$.\vspace{6pt}

For the fourth equality we apply \cite[Thm. 5.4]{equalitymixedOregular} to the matrix $\mathcal{V}^{\ast}$ and note: By $(ii)(c)$ in Theorem \ref{mainBMTthm} and the previously verified equivalence each $\omega_{(\mathbf{V}^{(\ell)})^{\ast}}$ satisfies \hyperlink{om1}{$(\omega_1)$} too and $\omega_{(\mathbf{V}^{(\ell_1)})^{\ast}}\hyperlink{sim}{\sim}\omega_{(\mathbf{V}^{(\ell_2)})^{\ast}}$ for any $\ell_1,\ell_2>0$. Moreover, recall that each $(\mathbf{V}^{(\ell)})^{\ast}$ satisfies \hyperlink{mg}{$(\on{mg})$} because $\mathbf{V}^{(\ell)}$ is log-convex; see comment $(f)$ in Section \ref{conjsequsection}. Finally note that $\mathcal{V}^{\ast}$ is apart from the pointwise order of the sequences formally standard log-convex as required in \cite[Thm. 5.4]{equalitymixedOregular}. For this note that $(\mathbf{V}^{(\ell)})^{\ast}\in\hyperlink{LCset}{\mathcal{LC}}$ for each $\ell>0$; recall again Remark \ref{Rem24reloaded}. But the failure of the pointwise order is only a minor technical issue and can be achieved when switching to a standard log-convex matrix $\mathcal{U}^{\ast}:=\{(\mathbf{U}^{(\ell)})^{\ast}: \ell>0\}$ such that $(\mathbf{V}^{(\ell)})^{\ast}\hyperlink{approx}{\approx}(\mathbf{U}^{(\ell)})^{\ast}$ for each $\ell>0$.

The fifth equality is clear by $(ii)(e)$ in Theorem \ref{mainBMTthm} whereas the last equality holds by the fact that $\mathcal{M}_{\omega}$ is constant; see $(ii)(d)$ in this result.
\qed\enddemo

\begin{remark}\label{rigidrem}
\emph{The equalities (as l.c.v.s.) stated in \eqref{mainBMTthm1equ} are also valid when replacing the symbol/functor $\mathcal{E}$ by some other weighted and analogously defined spaces (e.g. weighted sequences of complex numbers). Similarly, the function classes can also be assumed to be $H$-valued; see the definitions in \cite[Sect. 2.3 \& 2.4]{microclasses} for the weight sequence case. This is due to the fact that in the proof above we exclusively deal with weights and their growth and regularity properties, while no information concerning the underlying weighted category is required.}

\emph{In view of Corollary \ref{conjequivthmcor1} and the identities \eqref{conjequivthmcor1equ1} and \eqref{mainBMTthmproofequ} in order to deal with the conjugate matrix $\mathcal{M}^{\ast}_{\omega}$ and the corresponding (associated) weight functions $\omega^{\ast}$, $\omega_{\mathbf{W}^{(\ell)}}^{\ast}$, $\omega_{(\mathbf{W}^{(\ell)})^{\ast}}$, it seems to be useful and natural to assume $\overline{\gamma}(\omega)<1$. However, in this case both matrices $\mathcal{M}_{\omega}$ and $\mathcal{M}^{\ast}_{\omega}$ are constant and hence this assumption is a crucial restriction within the weight function setting. And $(iii)$ in Theorem \ref{mainBMTthm} provides more conditions under which both matrices are constant.}
\end{remark}

\subsection{Ultradifferentiable Braun-Meise-Taylor classes as weighted entire spaces}\label{Thm34trivialremsection}
We are interested in establishing an isomorphism between $\mathcal{E}_{[\omega]}$, defined via fast growing (non-standard) weights $\omega$, and weighted spaces of entire functions expressed in terms of $\omega^{\ast}$ or related conjugate (associated) weight functions. Thus the aim is to transfer the main statement \cite[Thm. 3.4]{microclasses} from the weight sequence setting to the BMT-weight function setting. As in \cite[Thm. 3.4]{microclasses} this isomorphism is provided via the extension operator $E: f\mapsto\sum_{k=0}^{+\infty}\frac{f^{(k)}(x_0)}{k!}(z-x_0)^k$, $x_0\in I$ fixed, $z\in\CC$. Here $I\subseteq\RR$ denotes an interval and $E$ is defined on $\mathcal{E}_{[\omega]}(I,H)$; i.e. the functions $f$ are considered to be $H$-valued with $H$ a Hilbert space. These spaces are defined analogously to the usual scalar-valued Braun-Meise-Taylor classes from \cite{BraunMeiseTaylor90} and for precise definitions concerning the weight sequence spaces $\mathcal{E}_{[\mathbf{M}]}(I,H)$ see \cite[Sect. 2.3 \& 2.4]{microclasses}. For the forthcoming comments and notions concerning weighted entire spaces, as in Corollary \ref{conjequivthmcor}, we refer to \cite[Sect. 3]{microclasses} and \cite{weightedentireinclusion1}.\vspace{6pt}

When inspecting the proof of \cite[Thm. 3.4]{microclasses} it turns out that when choosing $h:=1$ there the first part can be transferred for any $\omega\in\hyperlink{omset0}{\mathcal{W}_0}$ satisfying $t=o(\omega(t))$ as $t\rightarrow+\infty$: Indeed, by Lemma \ref{firstBMTlemma} each $\omega_{(\mathbf{W}^{(\ell)})^{\ast}}$ is well defined and then, in the \emph{Roumieu case,} for some $A,\ell>0$ and all $z\in\CC$ we have the estimate $\|E(f)\|\le 2A\exp(\omega_{(\mathbf{W}^{(\ell)})^{\ast}}(2|z|))$. Similarly, in the \emph{Beurling case,} for all $\ell>0$ there exists $A>0$ such that this estimate is valid for any $z\in\CC$. Thus the operator $E$ is continuous having a similar target space compared to the one from \cite[Thm. 3.4]{microclasses}: More precisely, the parameter is subject to the weight sequences in the matrix $\mathcal{M}_{\omega}$.

Note that for $\omega$ one might expect the corresponding \emph{``exponential-type weight function system''} as target space and which seems to be more natural: In the Roumieu case that there exist $C,\ell>0$ such that $C\exp(\ell\omega^{\ast}(|z|))$ for all $z\in\CC$ whereas in Beurling case for all $\ell>0$ there exists $C>0$ such that the previous estimate holds for all $z\in\CC$. Again the crucial parameter $\ell>0$ is subject to the associated weight sequences.

But in view of Theorem \ref{conjequivthm} resp. Corollary \ref{conjequivthmcor1}, $(ii)$ in Theorem \ref{mainBMTthm} (especially \eqref{mainBMTthmequ1}) a connection between both weighted entire notions without any additional assumption on $\omega$ is unclear.

However, in any case, by inspecting the second part of the proof where the continuity of the inverse (restriction) mapping is studied, it turns out that the log-convexity of the conjugate sequence $\mathbf{M}^{\ast}$ is indispensable, see \eqref{Prop32Komatsu} and Lemma \ref{lemma1}, and by $(b)$ in Section \ref{conjsequsection} in the weight function case one requires then the log-concavity for (each) $\mathbf{w}^{(\ell)}$. But $(iii)(a)$ in Theorem \ref{mainBMTthm} yields the fact that under this assumption both $\mathcal{M}_{\omega}$ and $\mathcal{M}^{\ast}_{\omega}$ are constant and $\omega$ satisfies \hyperlink{om6}{$(\omega_6)$} and this is too restrictive in the general case. Indeed, for applying $(iii)(a)$ in Theorem \ref{mainBMTthm} even the more general assumption that \emph{some $\left(\frac{\vartheta^{(\ell)}_p}{p}\right)_p$ is almost decreasing} is sufficient to conclude.

\section{On regularity results by M. Markin for BMT-weight functions}\label{Markinsection}
The aim is now to combine all information and to see how the results from \cite[Sect. 4]{microclasses} transfer to the BMT-weight function setting. For more explanations and the background concerning the small Gevrey setting we refer to \cite[Sect. 4]{microclasses}, the original work \cite{Markin01} and also to \cite{Gorbachuk}. So let $\omega\in\hyperlink{omset0}{\mathcal{W}_0}$ be given with associated weight matrix $\mathcal{M}_{\omega}:=\{\mathbf{W}^{(\ell)}: \ell>0\}$. (Indeed, for the next basic definitions it suffices to assume that $\omega$ is matrix admissible.)

First, for a \emph{densely defined (and closed) operator $A$ on a Hilbert space $H$} we set
$$C^{\infty}(A):=\bigcap_{n \in \NN} D(A^n),$$
where $D(A^n)$ is the domain of definition of $A^n$, the $n$-fold iteration of $A$; we refer to \cite[Sect. 2]{Markin01} and \cite[Sect. 1.1]{Gorbachuk} for this definition. Then consider the \emph{Roumieu-type class}
$$\mathcal{E}_{\{\mathcal{M}_{\omega}\}}(A):= \{f \in C^{\infty}(A): \exists\;C,h,\ell >0\;\forall\; p \in \NN:\;\;\; \|A^pf\|\le C h^p W^{(\ell)}_p\},$$
and the corresponding \emph{Beurling-type class} is defined by
$$\mathcal{E}_{(\mathcal{M}_{\omega})}(A):= \{f \in C^{\infty}(A): \forall\;h,\ell >0\;\exists\;C>0\;\forall\; p \in \NN:\;\;\; \|A^pf\|\le C h^p W^{(\ell)}_p\}.$$
Thus
$$\mathcal{E}_{\{\mathcal{M}_{\omega}\}}(A)=\bigcup_{\ell>0}\mathcal{E}_{\{\mathbf{W}^{(\ell)}\}}(A),\hspace{15pt}\mathcal{E}_{(\mathcal{M}_{\omega})}(A)=\bigcap_{\ell>0}\mathcal{E}_{(\mathbf{W}^{(\ell)})}(A),$$
where $\mathcal{E}_{\{\mathbf{W}^{(\ell)}\}}(A)$ and $\mathcal{E}_{(\mathbf{W}^{(\ell)})}(A)$ are defined accordingly; see \cite[Sect. 4.3]{microclasses}. These definitions are inspired by the Gevrey case treated in \cite[Sect. 2]{Markin01}, see also \cite[Sect. 1.2]{Gorbachuk} for more general sequences, and see the settings in \cite[Sect. 4.3]{microclasses} and in \cite[Sect. 4 \& 5]{compositionpaper}. Indeed, note that the analogous spaces can be introduced for an abstractly given weight matrix $\mathcal{M}=\{\mathbf{M}^{(\ell)}: \ell>0\}$. Moreover, inspired by the ultradifferentiable weight function setting in the sense of Braun, Meise, and Taylor from \cite{BraunMeiseTaylor90}, we set
$$\mathcal{E}_{\{\omega\}}(A):= \{f \in C^{\infty}(A): \exists\;C,\ell >0\;\forall\; p \in \NN:\;\;\; \|A^pf\|\le CW^{(\ell)}_p\},$$
and
$$\mathcal{E}_{(\omega)}(A):= \{f \in C^{\infty}(A): \forall\;\ell >0\;\exists\;C>0\;\forall\; p \in \NN:\;\;\; \|A^pf\|\le CW^{(\ell)}_p\}.$$
Next recall that via \cite[Sect. 1.3, p. 74 \& 75]{Gorbachuk} a different description of $\mathcal{E}_{\{\mathbf{W}^{(\ell)}\}}(A)$ resp. of $\mathcal{E}_{(\mathbf{W}^{(\ell)})}(A)$ (for each $\ell>0$) in terms of $E_A$, the spectral measure associated to $A$, can be deduced; see the first part of \cite[Thm. 1.3]{Gorbachuk} and its proof on \cite[p. 75]{Gorbachuk}:
\begin{equation}\label{RoumGorbachuck}
\mathcal{E}_{\{\mathbf{W}^{(\ell)}\}}(A)=\{ f \in H: \exists\;t>0~ \int_\CC e^{2 \omega_{\mathbf{W}^{(\ell)}}(t|\lambda|)} \langle dE_A(\lambda) f,f\rangle <+\infty\},
\end{equation}
and
\begin{equation}\label{BeurGorbachuck}
\mathcal{E}_{(\mathbf{W}^{(\ell)})}(A)=\{ f \in H: \forall\;t>0~ \int_\CC e^{2 \omega_{\mathbf{W}^{(\ell)}}(t|\lambda|)} \langle dE_A(\lambda) f,f\rangle <+\infty\}.
\end{equation}
Indeed, these equivalent representations are valid as l.c.v.s. when being equipped with their canonical topologies and this statement holds for any weight sequence $\mathbf{M}$ in the sense of Definition \ref{defweightsequ}: For this note that $\lim_{p\rightarrow+\infty}(M_p)^{1/p}=+\infty$ is precisely condition \cite[$(1.4)$]{Gorbachuk} and the function $\rho$ in \cite[$(1.5)$]{Gorbachuk} corresponds to $\exp\circ\omega_{\mathbf{M}}$. We also remark that the assumption for $\mathbf{M}$ to be non-decreasing in \cite[Thm. 1.3]{Gorbachuk} is not required in the proof and also log-convexity is not necessary for this statement. (However, each $\mathbf{W}^{(\ell)}$ satisfies these additional properties; see $(i)$ in Section \ref{assomatrixsection}.) Finally, recall that $\lim_{p\rightarrow+\infty}(M_p)^{1/p}=+\infty$ is natural when working with $\omega_{\mathbf{M}}$ in order to ensure that this function is well defined; see \cite[Lemma 2.2]{regularnew}.\vspace{6pt}

Write again $[\cdot]$ as a common notation for $\{\cdot\}$ and $(\cdot)$ and then, by involving the analogous growth controls in the definitions, \eqref{newexpabsorb} and $(iii)$ in Section \ref{assomatrixsection} we infer:

\begin{theorem}\label{Gorbachuckweightmatrixthm}
Let $\omega\in\hyperlink{omset0}{\mathcal{W}_0}$ be given and let $\mathcal{M}_{\omega}$ be the corresponding associated weight matrix. If $\omega$ satisfies in addition \hyperlink{om1}{$(\omega_1)$}, then as l.c.v.s.
$$\mathcal{E}_{[\mathcal{M}_{\omega}]}(A)=\mathcal{E}_{[\omega]}(A),$$
and \hyperlink{om6}{$(\omega_6)$} implies as l.c.v.s.
$$\forall\;\ell>0:\;\;\;\mathcal{E}_{[\mathcal{M}_{\omega}]}(A)=\mathcal{E}_{[\mathbf{W}^{(\ell)}]}(A).$$
Finally, the identities \eqref{RoumGorbachuck} and \eqref{BeurGorbachuck} hold for any $\ell>0$.
\end{theorem}

By taking into account the idea that in the weight function setting canonically one deals with a certain one-parameter family of weight sequences belonging to the set \hyperlink{LCset}{$\mathcal{LC}$}, namely the associated weight matrix $\mathcal{M}_{\omega}$, it is interesting to study the case $\mathfrak{F}=\mathcal{M}_{\omega}$ with $\mathfrak{F}$ denoting the abstractly given family of weight sequences treated in \cite[Sect. 4]{microclasses}. This is due to the fact that in this situation a precise controlled loss of regularity is possible in the sense that one describes all appearing sequences uniformly with a given single weight function $\omega\in\hyperlink{omset0}{\mathcal{W}_0}$. We summarize:

\begin{itemize}
\item[$(a)$] \cite[Thm. 4.1]{microclasses} transfers immediately to each sequence $\mathbf{W}^{(\ell)}$ since this result holds for arbitrary $\mathbf{M}\in\RR_{>0}^{\NN}$.

\item[$(b)$] Having $\mathfrak{F}=\mathcal{M}_{\omega}$ in the crucial result \cite[Lemma 4.3]{microclasses}, then by \eqref{goodequivalenceclassic} it follows that \cite[$(4.2)$]{microclasses}, i.e.
    \begin{equation}\label{42}
    \forall\;\mathbf{N}\in\mathfrak{F}\;\exists\;\mathbf{M}\in\mathfrak{F}:\;\;\;\omega_{\mathbf{M}}(2t)=O(\omega_{\mathbf{N}}(t)),\;\;\;t\rightarrow+\infty,
    \end{equation}
    implies \hyperlink{om1}{$(\omega_1)$} for $\omega$.

    Moreover, via the properties of the technical but crucial \emph{uniform bound sequence} $\mathbf{a}=(a_p)_{p\in\NN}$ stated in \cite[Lemma 4.3 (i) \& (ii)]{microclasses} if such a sequence exists for $\mathcal{M}_{\omega}$ then one has, in particular, that condition \eqref{firstBMTlemmaequ} has to be satisfied and so Lemma \ref{firstBMTlemma} implies the fact that all conjugates $\omega^{\ast}$, $\omega_{\mathbf{W}^{(\ell)}}^{\ast}$, $\omega_{(\mathbf{W}^{(\ell)})^{\ast}}$ are well defined.

  Conversely, \hyperlink{om1}{$(\omega_1)$} for $\omega$ implies via \eqref{goodequivalenceclassic} condition \hyperlink{om1}{$(\omega_1)$} for each $\omega_{\mathbf{W}^{(\ell)}}$ and so \eqref{42} (resp. \cite[$(4.2)$]{microclasses}) with $\mathbf{M}=\mathbf{W}^{(\ell)}=\mathbf{N}$.

\item[$(c)$] However, (also) for $\mathfrak{F}=\mathcal{M}_{\omega}$ it is crucial to see whether such a uniform bound sequence $\mathbf{a}$ exists and one is interested in constructing this sequence needed in \cite[Lemma 4.3]{microclasses} via the technical result \cite[Prop. 4.5]{microclasses}. Note that $\mathbf{a}$ should satisfy
    \begin{itemize}
    \item[$(*)$] $\lim_{p\rightarrow+\infty}(a_p)^{1/p}=0$,

    \item[$(*)$] $\mathbf{a}$ is a \emph{uniform bound for $\mathcal{M}_{\omega}$;} i.e.
    $$\forall\;\ell>0\;\exists\;C\ge 1\;\forall\;p\in\NN:\;\;\;\frac{W^{(\ell)}_p}{p!}=w^{(\ell)}_p\le Ca_p.$$
    \end{itemize}

\item[$(d)$] Indeed, normalization of each sequence $\mathbf{W}^{(\ell)}$ and pointwise order, i.e. \cite[Prop. 4.5 (i) \& (ii)]{microclasses}, are clear and \cite[Prop. 4.5 (iii)]{microclasses} amounts to have \eqref{firstBMTlemmaequ} which holds if $t=o(\omega(t))$ is valid; see again Lemma \ref{firstBMTlemma}.

But assertion \cite[Prop. 4.5 (iv)]{microclasses}, i.e. $p\mapsto(w^{(\ell)}_p)^{1/p}$ is non-increasing for each $\ell>0$, is not clear in general and in fact it is too strong resp. too restrictive in the general BMT-weight function setting: By
$(iii)(c)$ in Theorem \ref{mainBMTthm} we obtain that both matrices $\mathcal{M}_{\omega}$ and $\mathcal{M}^{\ast}_{\omega}$ are constant and hence in this case technically one deals with the single weight sequence situation (recall that even $\overline{\gamma}(\omega)\le 1$ holds). More precisely, one can choose any arbitrary $\mathbf{W}^{(\ell)}$ resp. $(\mathbf{W}^{(\ell)})^{\ast}$. Additionally, the fact that $\mathcal{M}_{\omega}$ is constant clearly violates assumption \cite[Prop. 4.5 (v)]{microclasses}.

(This problem is also suggested by $(ii)(e)$ in Theorem \ref{mainBMTthm}: When involving $\mathcal{V}$ we have that $p\mapsto(v^{(\ell)}_p)^{1/p}$ is non-increasing for any $\ell>0$ provided that $\overline{\gamma}(\omega)<1$. However, by taking into account $(ii)(d)$ in Theorem \ref{mainBMTthm}, see also Remark \ref{rigidrem}, under this assumption on the growth index again both matrices $\mathcal{M}_{\omega}$ and $\mathcal{M}^{\ast}_{\omega}$ are constant.)\vspace{6pt}

Summarizing, the assumptions on $\mathfrak{F}=\mathcal{M}_{\omega}$ in \cite[Prop. 4.5]{microclasses} are too restrictive in the general BMT-weight function setting and, indeed, they have been inspired by the family of all \emph{small Gevrey sequences;} i.e. $\mathbf{G}^{s}$ with $0<s<1$.
\end{itemize}

In order to overcome the problem described in the last comment we have to generalize \cite[Prop. 4.5]{microclasses} and, of course, this new version can also be applied to the weight sequence setting from \cite[Sect. 4]{microclasses}, i.e. to abstractly given families $\mathfrak{F}$ with no underlying weight function $\omega$.

\begin{proposition}\label{Prop45new}
Let $\mathfrak{F}:=\{\mathbf{N}^{(\beta)}\in\RR_{>0}^{\NN}: \beta>0\}$ be a one-parameter family of
sequences $\mathbf{N}^{(\beta)}$ which satisfies the following properties:
\begin{itemize}
\item[$(i)$] $\mathbf{N}^{(\beta_1)}\le\mathbf{N}^{(\beta_2)}\Leftrightarrow\mathbf{n}^{(\beta_1)}\le
    \mathbf{n}^{(\beta_2)}$ for all $0<\beta_1\le\beta_2$ (pointwise order),

\item[$(ii)$] $\lim_{p\rightarrow+\infty}(n^{(\beta)}_p)^{1/p}=0$ for each
    $\beta>0$.
\end{itemize}
Then there exists a sequence $\mathbf{a}=(a_p)_{p\in\NN}$ such that $p\mapsto(a_p)^{1/p}$ is non-increasing, $\lim_{p\rightarrow+\infty}(a_p)^{1/p}=0$, and that
\begin{equation}\label{afastvanishing}
\forall\;\beta>0:\;\;\;\lim_{p\rightarrow+\infty}\left(\frac{a_p}{n_p^{(\beta)}}\right)^{1/p}=+\infty,
\end{equation}
which implies, in particular, that $\mathbf{a}$ is a \emph{uniform bound for $\mathfrak{F}$.}

If in addition $\mathfrak{F}$ satisfies
\begin{itemize}
\item[$(iii)$] $\lim_{j\rightarrow+\infty}\left(\frac{N^{(\beta_2)}_j}{N^{(\beta_1)}_j}\right)^{1/j}=\lim_{j\rightarrow+\infty}\left(\frac{n^{(\beta_2)}_j}{n^{(\beta_1)}_j}\right)^{1/j}=+\infty$
for all $0<\beta_1<\beta_2$ (large growth difference between the sequences), then $\mathfrak{F}$ satisfies \eqref{42}.
\end{itemize}
\end{proposition}

\demo{Proof}
Fix two strictly increasing sequences of positive real numbers $(c_k)_{k\in\NN_{>0}}$ and $(d_k)_{k\in\NN_{>0}}$ such that $c_1=1=d_1$ and $\lim_{k\rightarrow+\infty}c_k=+\infty=\lim_{k\rightarrow+\infty}d_k$. Put $j_1:=1$ and then choose iteratively values $j_{k+1}\in\NN_{>0}$, with $j_{k+1}>j_k$, such that
\begin{equation}\label{uniformboundequ1}
\forall\;k\in\NN_{>0}:\;\;\;(n_{j_k}^{(k)})^{1/j_k}>c_k(n_{j_{k+1}}^{(k+1)})^{1/j_{k+1}},
\end{equation}
and
\begin{equation}\label{uniformboundequ2}
\forall\;k\in\NN_{>0}\;\forall\;j\ge j_{k+2}:\;\;\,\frac{(n^{(k+1)}_{j_{k+1}})^{1/j_{k+1}}}{(n^{(k)}_j)^{1/j}}\ge d_k.
\end{equation}
For \eqref{uniformboundequ1} and \eqref{uniformboundequ2} we have only used properties $(i)$ and $(ii)$.

Now introduce $\mathbf{a}:=(a_p)_{p\in\NN}$ as follows: We set $a_0:=1$ and
\begin{equation}\label{uniformboundequ3}
(a_j)^{1/j}:=(n_{j_{1}}^{(1)})^{1/j_{1}},\;\;\;1=j_1\le j<j_2,\hspace{15pt}(a_j)^{1/j}:=(n_{j_{k}}^{(k)})^{1/j_{k}},\;\;\;j_{k+1}\le j<j_{k+2},\;k\in\NN_{>0}.
\end{equation}
Thus we have by definition and \eqref{uniformboundequ1} that $j\mapsto(a_j)^{1/j}$ is non-increasing and tending to $0$. Indeed, the rate of approaching $0$ is expressed in terms of the auxiliary sequence $(c_k)_{k\in\NN_{>0}}$ and can be chosen as fast as desired.

Let $k_0\in\NN_{>0}$ be given and from now on fixed. For $j\ge
j_{k_0+2}$ we can find $k\ge k_0$ such that $j_{k+2}\le
j<j_{k+3}$. Thus, by \eqref{uniformboundequ3} for all such $j$ we can estimate as follows:
\begin{align*}
\frac{a_j^{1/j}}{(n^{(k_0)}_j)^{1/j}}=\frac{(n^{(k+1)}_{j_{k+1}})^{1/j_{k+1}}}{(n^{(k_0)}_j)^{1/j}}\ge\frac{(n^{(k+1)}_{j_{k+1}})^{1/j_{k+1}}}{(n^{(k)}_j)^{1/j}}\ge d_k,
\end{align*}
hence $\lim_{j\rightarrow+\infty}\frac{a_j^{1/j}}{(n^{(k_0)}_j)^{1/j}}=+\infty$. The first inequality follows from the pointwise order and the second one by \eqref{uniformboundequ2}. The rate of tending to infinity is expressed in terms of $(d_k)_{k\in\NN_{>0}}$ and can be chosen as fast as desired; concerning the arbitrariness of $(c_k)_{k\in\NN_{>0}}$ and $(d_k)_{k\in\NN_{>0}}$ we also refer to Remark \ref{counterexrem} below.

Finally, for each given $\beta>0$ by the pointwise order we can find some $k_0\in\NN_{>0}$ such that $\frac{a_j^{1/j}}{(n^{(\beta)}_j)^{1/j}}\ge\frac{a_j^{1/j}}{(n^{(k_0)}_j)^{1/j}}$ for all $j\ge 1$ and hence $\mathbf{a}$ satisfies \eqref{afastvanishing}, too. Note that the definition of $\mathbf{a}$ differs from the one in \cite[Prop. 4.5]{microclasses} by an index shift (in \eqref{uniformboundequ2}). This is because Proposition \ref{Prop45new} is more general and does not use \cite[Prop. 4.5 (iv) \& (v)]{microclasses}.\vspace{6pt}

The last part, i.e. showing $(iii)$, follows as in \cite[Prop. 4.5]{microclasses}.
\qed\enddemo

\begin{remark}\label{Prop45newrem}
\emph{We recognize that $(i)$ and $(iii)$ in Proposition \ref{Prop45new} are preserved when multiplying each $\mathbf{N}^{(\beta)}$ pointwise with an arbitrary given (and fixed) $\mathbf{M}\in\RR_{>0}^{\NN}$; i.e. consider $\mathfrak{F}\cdot\mathbf{M}:=\{\mathbf{N}^{(\beta)}\cdot\mathbf{M}: \beta>0\}$. This technical fact has the following meaning: Assume that $\mathfrak{F}:=\{\mathbf{N}^{(\beta)}: \beta>0\}$ and $\mathbf{M}$ are given such that $(i)$ and $(iii)$ in Proposition \ref{Prop45new} are valid for $\mathfrak{F}$ and such that $\mathbf{N}^{(\beta)}\hyperlink{mtriangle}{\vartriangleleft}\mathbf{M}$ for all $\beta>0$ which should be compared with \eqref{firstBMTlemmaequ} and note that $(ii)$ is precisely this relation with $\mathbf{M}=\mathbf{G}^1$. Then $\mathfrak{F}\cdot\mathbf{m}^{-1}$ satisfies all assumptions $(i)-(iii)$ and Proposition \ref{Prop45new} applied to $\mathfrak{F}\cdot\mathbf{m}^{-1}$ gives a uniform bound $\mathbf{a}$ for this matrix. For $\mathfrak{F}$ the sequence $\mathbf{a}\cdot\mathbf{m}$ is then relevant and hence in \eqref{afastvanishing} one should write $a_pm_p$ in the numerator.}
\end{remark}

By taking into account Proposition \ref{Prop45new}, Lemma \ref{firstBMTlemma} and comment $(b)$ in this section we infer:

\begin{theorem}\label{uniformboundthm}
Let $\omega\in\hyperlink{omset0}{\mathcal{W}_0}$ be given with associated weight matrix $\mathcal{M}_{\omega}=\{\mathbf{W}^{(\ell)}: \ell>0\}$ and assume that $\omega$ satisfies $t=o(\omega(t))$ as $\rightarrow+\infty$ (i.e. $(\mathfrak{C}_2)$). Then there exists $\mathbf{a}=(a_p)_{p\in\NN}$ such that $p\mapsto(a_p)^{1/p}$ is non-increasing, $\lim_{p\rightarrow+\infty}(a_p)^{1/p}=0$, and $\mathbf{a}$ is a uniform bound for $\mathcal{M}_{\omega}$, more precisely even
$$\forall\;\ell>0:\;\;\;\lim_{p\rightarrow+\infty}\left(\frac{a_p}{w_p^{(\ell)}}\right)^{1/p}=+\infty.$$
If in addition \hyperlink{om1}{$(\omega_1)$} for $\omega$ holds, then for each $\omega_{\mathbf{W}^{(\ell)}}$ this condition is valid too and so \eqref{42} (i.e. \cite[$(4.2)$]{microclasses}) with $\mathbf{M}=\mathbf{W}^{(\ell)}=\mathbf{N}$, $\ell>0$ arbitrary.
\end{theorem}

Using this result we immediately get the following variant of the crucial \cite[Lemma 4.3]{microclasses}:

\begin{corollary}\label{uniformboundcor}
Let $\omega\in\hyperlink{omset0}{\mathcal{W}_0}$ be given with associated weight matrix $\mathcal{M}_{\omega}=\{\mathbf{W}^{(\ell)}: \ell>0\}$ and assume that $\omega$ satisfies $t=o(\omega(t))$ as $\rightarrow+\infty$ (i.e. $(\mathfrak{C}_2)$) and \hyperlink{om1}{$(\omega_1)$}.

If one has as sets the equality
\begin{equation}\label{uniformboundcorequ}
(\mathcal{E}_{\{\omega\}}(A)=)\bigcup_{\ell>0}\mathcal{E}_{\{\mathbf{W}^{(\ell)}\}}(A)=\mathcal{E}_{(\mathbf{G}^1)}(A),
\end{equation}
then $A$ is bounded.
\end{corollary}

\demo{Proof}
By Theorem \ref{uniformboundthm} we can apply the proof of \cite[Lemma 4.3]{microclasses} to $\mathfrak{F}=\mathcal{M}_{\omega}$. Recall that the first equality in \eqref{uniformboundcorequ} holds by Theorem \ref{Gorbachuckweightmatrixthm}.
\qed\enddemo

Thus, the crucial result \cite[Thm. 4.8]{microclasses} turns into the following statement:

\begin{theorem}\label{Thm48new}
Let $\omega\in\hyperlink{omset0}{\mathcal{W}_0}$ be given with associated weight matrix $\mathcal{M}_{\omega}=\{\mathbf{W}^{(\ell)}: \ell>0\}$ and assume that $\omega$ satisfies $t=o(\omega(t))$ as $\rightarrow+\infty$ (i.e. $(\mathfrak{C}_2)$) and \hyperlink{om1}{$(\omega_1)$}.

If for any weak solution $y$ of \eqref{evolequ} on $[0,+\infty)$ there exists $\ell>0$ such that $y\in\mathcal{E}_{\{\mathbf{W}^{(\ell)}\}}([0,+\infty),H)$, then the operator $A$ is bounded.
\end{theorem}

\demo{Proof}
We follow the proof of \cite[Thm. 4.8]{microclasses} and apply it to the family $\mathfrak{F}=\mathcal{M}_{\omega}$. For this we involve \cite[Thm. 4.1]{microclasses} applied to a/the sequence $\mathbf{W}^{(\ell)}$ and Corollary \ref{uniformboundcor}. Concerning the definition of the $H$-valued (local) weighted function class $\mathcal{E}_{\{\mathbf{W}^{(\ell)}\}}([0,+\infty),H)$ we refer again to \cite[Sect. 2.3 \& 2.4]{microclasses}.
\qed\enddemo

\begin{remark}\label{Thm48newrem}
\emph{One shall note that the assumption in Theorem \ref{Thm48new} does not mean that each weak solution $y$ belongs to $\mathcal{E}_{\{\omega\}}([0,+\infty),H)=\mathcal{E}_{\{\mathcal{M}_{\omega}\}}([0,+\infty),H)$, where this equality is valid as l.c.v.s. by \eqref{newexpabsorb} via \hyperlink{om1}{$(\omega_1)$} and the arguments from \cite[Thm. 5.14 (2)]{compositionpaper}. Recall that for $\mathcal{E}_{\{\mathbf{W}^{(\ell)}\}}([0,+\infty),H)$ it is required to have the defining estimates for each compact set (interval) $K\subseteq[0,+\infty)$ with the \emph{same weight sequence} $\mathbf{W}^{(\ell)}$ and hence uniformity w.r.t. the matrix parameter $\ell>0$ on whole $[0,+\infty)$, whereas for the local classes $\mathcal{E}_{\{\omega\}}([0,+\infty),H)$ naturally this parameter is allowed to depend on chosen $K$.}
\end{remark}

Finally, we transfer \cite[Sect. 4.4, Thm. 4.9]{microclasses} where a connection to the weighted entire setting is established: More precisely, that each weak solution $y$ can be extended to an entire function satisfying a precise prescribed growth. For this again the previously described technical issue for $\mathfrak{F}=\mathcal{M}_{\omega}$ appears: Assertions \cite[Sect. 4.4 (i)-(iii)]{microclasses} yield no problems, concerning (iii) recall Proposition \ref{Prop45new}. But \cite[Sect. 4.4 (iv)]{microclasses} means \emph{log-concavity for each $\mathbf{w}^{(\ell)}$} and which is required in order to involve \cite[Thm. 3.4]{microclasses}; see also Corollary \ref{conjequivthmcor} and the explanations in Section \ref{Thm34trivialremsection}. As shown before this requirement is too restrictive in the sense that all sequences are equivalent and hence one loses the flexibility of the matrix parameter. However, for the sake of completeness when combining Theorem \ref{Thm48new} with \cite[Thm. 3.4]{microclasses} (and Corollary \ref{conjequivthmcor}) let us formulate the final statement:

\begin{theorem}\label{Thm49new}
Let $\omega\in\hyperlink{omset0}{\mathcal{W}_0}$ be given with associated weight matrix $\mathcal{M}_{\omega}=\{\mathbf{W}^{(\ell)}: \ell>0\}$. Assume that $\omega$ satisfies
\begin{itemize}
\item[$(i)$] $t=o(\omega(t))$ as $\rightarrow+\infty$ (i.e. $(\mathfrak{C}_2)$),

\item[$(ii)$] \hyperlink{om1}{$(\omega_1)$},

\item[$(iii)$] each $\left(\frac{\vartheta^{(\ell)}_p}{p}\right)_p$ is almost decreasing.
\end{itemize}

If for any weak solution $y$ of \eqref{evolequ} there is $\ell>0$, $C>0$ and $d>0$ such that $y$ can be extended to an entire function satisfying
\begin{equation}\label{Thm49newequ}
\|y(z)\|\le C\exp(\omega_{(\mathbf{W}^{(\ell)})^{\ast}}(d|z|)),
\end{equation}
then the operator $A$ is bounded. Alternatively, in \eqref{Thm49newequ} one can involve the bound $C\exp(\omega^{\ast}_{\mathbf{W}^{(\ell)}}(d|z|))$ and recall that $(iii)$ implies $\overline{\gamma}(\omega)\le 1$ and hence the restrictive condition \hyperlink{om6}{$(\omega_6)$}.

If assertion $(iii)$ is replaced by the stronger condition $\overline{\gamma}(\omega)<1$, then in \eqref{Thm49newequ} alternatively one can use the bound $C\exp(\omega^{\ast}(d|z|))$ and also the \emph{``exponential-type weight function system bounds''} $C\exp(d\omega_{(\mathbf{W}^{(\ell)})^{\ast}}(|z|))$, $C\exp(d\omega^{\ast}_{\mathbf{W}^{(\ell)}}(|z|))$, $C\exp(d\omega^{\ast}(|z|))$.
\end{theorem}

\demo{Proof}
Since $\mathbf{W}^{(\ell)}\in\hyperlink{LCset}{\mathcal{LC}}$ for each $\ell>0$, in view of assumption $(iii)$ and Remark \ref{Rem24reloaded} we have that for each $\ell>0$ there exists $\mathbf{V}^{(\ell)}\in\hyperlink{LCset}{\mathcal{LC}}$ such that $\mathbf{W}^{(\ell)}\hyperlink{approx}{\approx}\mathbf{V}^{(\ell)}$, $(\mathbf{W}^{(\ell)})^{\ast}\hyperlink{approx}{\approx}(\mathbf{V}^{(\ell)})^{\ast}$ and $(\mathbf{V}^{(\ell)})^{\ast}\in\hyperlink{LCset}{\mathcal{LC}}$. Equivalence $\mathbf{W}^{(\ell)}\hyperlink{approx}{\approx}\mathbf{V}^{(\ell)}$ preserves the (local) ultradifferentiable classes, i.e. $\mathcal{E}_{[\mathbf{W}^{(\ell)}]}(I,H)=\mathcal{E}_{[\mathbf{V}^{(\ell)}]}(I,H)$ for each $\ell>0$, and then apply \cite[Thm. 3.4]{microclasses} (resp. Corollary \ref{conjequivthmcor}) to (each) $\mathbf{V}^{(\ell)}$ and recall \eqref{assofctrel}: $(\mathbf{W}^{(\ell)})^{\ast}\hyperlink{approx}{\approx}(\mathbf{V}^{(\ell)})^{\ast}$ implies $\omega_{(\mathbf{W}^{(\ell)})^{\ast}}\hyperlink{simc}{\sim_{\mathfrak{c}}}\omega_{(\mathbf{V}^{(\ell)})^{\ast}}$ and this relation preserves the bound in \eqref{Thm49newequ}. Thus Theorem \ref{Thm48new} gives the conclusion. The supplement concerning the alternative bound follows by Corollary \ref{conjequivthmcor} applied to $\mathbf{W}^{(\ell)}$.

If even $\overline{\gamma}(\omega)<1$ holds, then $(ii)$ in Theorem \ref{mainBMTthm} implies the additional information on the bounds: In this case all conjugate functions are equivalent (see \eqref{mainBMTthmequ1}) and satisfy both \hyperlink{om1}{$(\omega_1)$} and \hyperlink{om6}{$(\omega_6)$}.
\qed\enddemo

We finish by constructing explicitly a weight function $\omega$ satisfying all requirements in Theorem \ref{Thm48new} and such that the corresponding associated weight matrix is non-constant which means that under the assumptions on $\omega$ in this result the flexibility w.r.t. the matrix parameter can become crucial. This also illustrates then the difference between Theorems \ref{Thm48new} and \ref{Thm49new}. In order to proceed we focus on $\omega=\omega_{\mathbf{M}}$.

\begin{proposition}\label{Thm48newexample}
There exists $\omega\in\hyperlink{omset0}{\mathcal{W}_0}$ satisfying $(\mathfrak{C}_2)$, \hyperlink{om1}{$(\omega_1)$} and such that the associated weight matrix $\mathcal{M}_{\omega}$ is non-constant (equivalently \hyperlink{om6}{$(\omega_6)$} for $\omega$ fails).
\end{proposition}

\emph{Note:} In view of Remark \ref{om1om6indexrem} and Lemma \ref{indexlemma1} the weight $\omega$ constructed in this proposition satisfies
$$0<\gamma(\omega)\le 1<\overline{\gamma}(\omega)=+\infty.$$
This also illustrates that in general in $(i)$ in Theorem \ref{indexlemma} one has $1-\overline{\gamma}(\sigma)<\gamma(\sigma^{\ast})$.

\demo{Proof}
We construct a sequence $\mathbf{M}$ such that $\omega_{\mathbf{M}}$ satisfies all required properties. For this consider the corresponding sequence of quotients $\mu=(\mu_p)_{p\in\NN}$ with $\mu_0:=1$ and set $M_p:=\prod_{i=0}^p\mu_i$. Let us fix $\epsilon\in(0,1)$ and an arbitrary non-decreasing sequence of positive real numbers $(d_j)_{j\in\NN_{>0}}$ satisfying $d_j\ge 1$ for all $j$ and $\lim_{j\rightarrow+\infty}d_j=+\infty$. Next introduce several auxiliary sequences: First, consider $(a_p)_{p\in\NN}$ by
\begin{equation}\label{exampleequ1}
a_p:=2^p,\;\;\;p\in\NN.
\end{equation}
Second, let $(c_j)_{j\in\NN_{>0}}$ be a strictly increasing sequence in $\NN_{\ge 2}$ such that
\begin{equation}\label{exampleequ5}
\forall\;j\in\NN_{>0}:\;\;\;\frac{2}{1+\epsilon}>(jd_j)^{1/c_j}
\end{equation}
and define the corresponding sequence $(p_j)_{j\in\NN}$ by
\begin{equation}\label{exampleequ2}
p_0:=0,\hspace{15pt}p_{j+1}:=c_{j+1}+p_j,\;\;\;j\in\NN.
\end{equation}
Note that the choice in \eqref{exampleequ5} is possible since $1+\epsilon<2$. Moreover, put
\begin{equation}\label{exampleequ3}
b_0:=1,\hspace{15pt}b_{p+1}:=(1+\epsilon)b_{p},\;\;\;p_j\le p\le p_{j+1}-2,\;j\in\NN,
\end{equation}
\begin{equation}\label{exampleequ3var}
b_{p+1}:=d_jb_{p},\;\;\;p=p_j-1,\;j\in\NN_{>0}.
\end{equation}
Finally, the sequence $\mu$ of quotients is given by
\begin{equation}\label{exampleequ4}
\mu_0:=1,\hspace{15pt}\mu_q:=b_p,\;\;\;a_p\le q<a_{p+1},\;p\in\NN.
\end{equation}
\emph{Claim I} $\mathbf{M}\in\hyperlink{LCset}{\mathcal{LC}}$ holds and hence $\omega_{\mathbf{M}}\in\hyperlink{omset0}{\mathcal{W}_0}$ (see Lemma \ref{assoweightomega0}).

By \eqref{exampleequ3}, \eqref{exampleequ3var} and the choice of $(d_j)_j$ it follows that $p\mapsto b_p$ is non-decreasing and $\lim_{p\rightarrow+\infty}b_p=+\infty$. Normalization $1=M_0\le M_1$ amounts to have $\mu_1\ge 1$ and which holds since $\mu_1=\mu_{a_0}=b_0=1$.\vspace{6pt}

\emph{Claim II} $\lim_{p\rightarrow+\infty}\frac{\mu_p}{p}=0$ holds.

By \eqref{exampleequ4} we see that $\frac{\mu_{q+1}}{q+1}<\frac{\mu_q}{q}$ for all $a_p\le q<a_{p+1}-1$, $p\in\NN$. When $q=a_{p+1}-1$, then in view of \eqref{exampleequ3} and \eqref{exampleequ3var} we distinguish: If $p+1\neq p_{j+1}$, $j\in\NN$, then $$\frac{\mu_{q+1}}{q+1}=\frac{\mu_{a_{p+1}}}{a_{p+1}}=\frac{b_{p+1}}{a_{p+1}}=\frac{(1+\epsilon)b_p}{2a_p}<\frac{b_p}{a_p}=\frac{\mu_{a_p}}{a_p}.$$ Finally, if $p+1=p_{j+1}$ for some $j\in\NN$, then
\begin{align*}
\frac{\mu_{q+1}}{q+1}&=\frac{\mu_{a_{p_{j+1}}}}{a_{p_{j+1}}}=\frac{b_{p_{j+1}}}{2^{c_{j+1}}a_{p_j}}=\frac{d_{j+1}b_{p_{j+1}-1}}{2^{c_{j+1}}a_{p_j}}=\frac{d_{j+1}(1+\epsilon)^{c_{j+1}-1}b_{p_j}}{2^{c_{j+1}}a_{p_j}}\le\frac{d_{j+1}(1+\epsilon)^{c_{j+1}}b_{p_j}}{2^{c_{j+1}}a_{p_j}}
\\&
\le\frac{1}{j+1}\frac{b_{p_j}}{a_{p_j}}=\frac{1}{j+1}\frac{\mu_{a_{p_{j}}}}{a_{p_{j}}},
\end{align*}
by the choice of the sequence $(c_j)_{j\in\NN_{>0}}$ in \eqref{exampleequ5}. Summarizing, the claim is shown.\vspace{6pt}

\emph{Claim III} $\lim_{p\rightarrow+\infty}(m_p)^{1/p}=0$ holds and hence $\omega_{\mathbf{M}}$ satisfies $(\mathfrak{C}_2)$ (see Lemma \ref{C2assoweightfctlemma} and Corollary \ref{conjwelldefcor}).

Since $\mathbf{M}\in\hyperlink{LCset}{\mathcal{LC}}$ we get $(M_p)^{1/p}\le\mu_p\Leftrightarrow\mu_1\cdots\mu_p\le\mu_p^p$ for all $p\in\NN_{>0}$ (and the second estimate is valid for $p=0$, too). Now, by \eqref{Stirlinglike} we obtain
$$\frac{p}{e}(m_p)^{1/p}\le(p!m_p)^{1/p}=(M_p)^{1/p}\le\mu_p=\frac{\mu_p}{p}p,$$
and so $(m_p)^{1/p}\le e\frac{\mu_p}{p}$ for all $p\in\NN_{>0}$ which yields the conclusion by using \emph{Claim II.}\vspace{6pt}

\emph{Claim IV} $\omega_{\mathbf{M}}$ satisfies \hyperlink{om1}{$(\omega_1)$}. In order to conclude it suffices to check that $\liminf_{p\rightarrow+\infty}\frac{\mu_{Qp}}{\mu_p}>1$ for some $Q\in\NN_{\ge 2}$ (see \eqref{beta3}); for a direct proof of this implication we refer to \cite[Lemma 12, $(2)\Rightarrow(4)$]{BonetMeiseMelikhov07} (cf. also \cite[Prop. 3.4]{subaddlike}) and see \cite[Thm. 2.11, Cor. 2.14, Thm. 3.11, Cor. 4.6 $(i)$]{index} involving growth index arguments.

By definition in \eqref{exampleequ4} one has $\frac{\mu_{2p}}{\mu_p}\ge 1+\epsilon$ for all $p\in\NN_{>0}$ sufficiently large, since $d_j\ge 1+\epsilon$ for all $j\in\NN_{>0}$ large enough and so the desired property holds with $Q:=2$.\vspace{6pt}

\emph{Claim V} $\omega_{\mathbf{M}}$ does not satisfy \hyperlink{om6}{$(\omega_6)$} which is equivalent to the fact that $\mathbf{M}$ violates \hyperlink{mg}{$(\on{mg})$} (see Lemma \ref{assoweightomega0}) and equivalently $\mathcal{M}_{\omega_{\mathbf{M}}}$ is non-constant (recall $(iii)$ in Section \ref{assomatrixsection}).

We use the characterization shown in \cite[Lemma 2.2]{whitneyextensionweightmatrix}, recall again the discussion in \cite{modgrowthstrange}, and verify $\sup_{p\in\NN}\frac{\mu_{2p}}{\mu_p}=+\infty$: For this choose $p=a_{p_j-1}$, $j\in\NN_{>0}$, then by \eqref{exampleequ3var} and \eqref{exampleequ4} we get $\frac{\mu_{2p}}{\mu_p}=\frac{\mu_{a_{p_j}}}{\mu_{a_{p_j-1}}}=\frac{d_jb_{p_j-1}}{b_{p_j-1}}=d_j$ and recall that $d_j\rightarrow+\infty$.
\qed\enddemo

\begin{remark}\label{counterexrem}
\emph{The proof of Proposition \ref{Thm48newexample} is quite technical but it can directly be generalized when replacing $j$ in \eqref{exampleequ5} by $f_j$ such that $(f_j)_{j\in\NN_{>0}}$ is a further given strictly (fast) increasing sequence. Of course, $(f_j)_{j\in\NN_{>0}}$ via the choice of $(c_j)_{j\in\NN_{>0}}$ will effect then the definition of $\mathbf{M}$ but we get that the conclusions in the above result even hold when $p\mapsto\frac{\mu_p}{p}$ resp. $p\mapsto(m_p)^{1/p}$ is required to tend faster to $0$ and hence the rate of approaching $0$ is not effecting the validity of this (counter-)example.}

\emph{On the other hand, the sequence $(d_j)_{j\in\NN_{>0}}$ is used to achieve that $\mathbf{M}$ violates \hyperlink{mg}{$(\on{mg})$} and hence $\mathcal{M}_{\omega_{\mathbf M}}$ fails $(iii)$ in Theorem \ref{Thm49new}.}
\end{remark}

\appendix
\section{On slowly varying associated weight functions}\label{slowlyvaryingappendix}
Motivated by the comment in Section \ref{growthsection} that for a given weight function $\sigma$ the reflexivity $\sigma\hyperlink{omvartrianglec}{\vartriangleleft_{\mathfrak{c}}}\sigma$ implies the fact that $\sigma$ has to be \emph{slowly varying,} the goal of this appendix is to investigate this condition for $\omega_{\mathbf{M}}$ in terms of the defining sequence $\mathbf{M}$.

\begin{proposition}\label{lowvarassofctcharact}
Let $\mathbf{M}$ be a log-convex weight sequence and recall the notation $\mu_p:=\frac{M_p}{M_{p-1}}$, $p\in\NN_{>0}$.
\begin{itemize}
\item[$(i)$] When $\mathbf{M}$ satisfies
\begin{equation}\label{beta3}
\exists\;Q\in\NN_{\ge 2}:\;\;\;\liminf_{p\rightarrow+\infty}\frac{\mu_{Qp}}{\mu_p}>1
\end{equation}
and
\begin{equation}\label{lowvarassofctcharactequ1}
\lim_{p\rightarrow+\infty}\frac{\mu_p}{(M_{p-1}/M_0)^{1/(p-1)}}=\lim_{p\rightarrow+\infty}\frac{\mu_p}{(M_{p-1})^{1/(p-1)}}=+\infty,
\end{equation}
then $\omega_{\mathbf{M}}$ is \emph{slowly varying;} i.e. \eqref{slowlyvarying}:
$$\forall\;u>0:\;\;\;\lim_{t\rightarrow+\infty}\frac{\omega_{\mathbf{M}}(ut)}{\omega_{\mathbf{M}}(t)}=1.$$

\item[$(ii)$] When $\omega_{\mathbf{M}}$ is slowly varying, then \eqref{lowvarassofctcharactequ1} has to be satisfied.
\end{itemize}
\end{proposition}

\demo{Proof}
First recall that for any weight function $\omega$ it suffices to check \eqref{slowlyvarying} for all $u>1$: When $0<u<1$, then set $s:=ut$ and get $\lim_{t\rightarrow+\infty}\frac{\omega(ut)}{\omega(t)}=\lim_{s\rightarrow+\infty}\frac{\omega(s)}{\omega(s/u)}=1$.\vspace{6pt}

$(i)$ Let $u>1$ large be given and fixed and let $t\ge\mu_1$. Then there exist $p,q\in\NN_{>0}$, $q\ge p$, such that $\mu_p\le t<\mu_{p+1}$ and $\mu_q\le ut<\mu_{q+1}$. For this note that $p\mapsto\mu_p$ is non-decreasing and $\lim_{p\rightarrow+\infty}\mu_p=+\infty$. In view of Lemma \ref{lemma1} we have
\begin{align*}
\frac{\omega_{\mathbf{M}}(ut)}{\omega_{\mathbf{M}}(t)}&=\frac{\log(M_0)+q\log(ut)-\log(M_{q})}{\log(M_0)+p\log(t)-\log(M_p)}
\\&
=\frac{\log(M_0)+p\log(t)+(q-p)\log(t)+q\log(u)-\log(M_p)+\log(M_p)-\log(M_q)}{\log(M_0)+p\log(t)-\log(M_p)}
\\&
=1+\frac{(q-p)\log(t)+q\log(u)-\log(M_q/M_p)}{\log(M_0)+p\log(t)-\log(M_p)}\le 1+\frac{q\log(u)}{\log(M_0)+p\log(t)-\log(M_p)}.
\end{align*}
The estimate is valid since $\mu_p\le t<\mu_{p+1}$ and so by log-convexity
$$(q-p)\log(t)-\log(M_q/M_p)=(q-p)\log(t)-\sum_{i=p+1}^q\log(\mu_i)\le(q-p)\log(\mu_{p+1})-(q-p)\log(\mu_{p+1})=0.$$
Let $\epsilon>0$ be given and then, in order to ensure that $\omega_{\mathbf{M}}$ is slowly varying, we put $t=\mu_p$ and so it suffices to require for all sufficiently large $p\ge 2$ that
\begin{align*}
&\frac{q\log(u)}{\log(M_0)+p\log(\mu_p)-\log(M_p)}\le\epsilon\Leftrightarrow\log(u^q)\le\log((M_0\mu_p^p/M_p)^{\epsilon})
\\&
\Leftrightarrow u^{q/\epsilon}M_p/M_0=u^{q/\epsilon}\mu_1\cdots\mu_p\le\mu_p^p\Leftrightarrow u^{q/\epsilon}\mu_1\cdots\mu_{p-1}\le\mu_p^{p-1}\Leftrightarrow u^{\frac{q}{(p-1)\epsilon}}(M_{p-1}/M_0)^{1/(p-1)}\le\mu_p.
\end{align*}

Now, since $\mathbf{M}$ satisfies \eqref{beta3} we have the following:
$$\exists\;Q\in\NN_{\ge 2}\;\exists\;\eta>0\;\exists\;p_{\eta}\in\NN_{>0}\;\forall\;p\ge p_{\eta}:\;\;\;\frac{\mu_{Qp}}{\mu_p}\ge 1+\eta.$$
Iteration of this estimate yields $\frac{\mu_{Q^np}}{\mu_p}\ge(1+\eta)^n$ for all $p\ge p_{\eta}$ and all $n\in\NN_{>0}$. Hence, when $u>1$ is given, then for $n$ chosen sufficiently large to ensure $u<(1+\eta)^n$ we obtain: For any $t\ge\mu_{p_{\eta}}$ one can find $p\ge p_{\eta}$ such that $\mu_p\le t<\mu_{p+1}$ and then $ut<u\mu_{p+1}<(1+\eta)^n\mu_{p+1}\le\mu_{Q^n(p+1)}$. Thus in the above computation, since $\mu$ is non-decreasing, for the index $q$ we have the upper bound $Q^n(p+1)-1$ and therefore $u^{\frac{q}{(p-1)\epsilon}}\le u^{\frac{Q^n(p+1)}{(p-1)\epsilon}}\le u^{\frac{2Q^n}{\epsilon}}$ for all $p\ge\max\{3,p_{\eta}\}$. Note that both $Q$ and $\eta$ are not depending on given $u$ and $\epsilon$, but $n$ is depending on $u$.

Finally, by \eqref{lowvarassofctcharactequ1} we have that $u^{\frac{2Q^n}{\epsilon}}(M_{p-1}/M_0)^{1/(p-1)}\le\mu_p$ is valid for all $p\ge p_{u,\epsilon}(\ge\max\{3,p_{\eta}\})$ sufficiently large enough.\vspace{6pt}

$(ii)$ By involving \eqref{assointrepr} for any $u>1$ and $t>\mu_1$ we rewrite and estimate the crucial quotient under consideration as follows:
$$\frac{\omega_{\mathbf{M}}(ut)}{\omega_{\mathbf{M}}(t)}=\frac{\int_{\mu_1}^{tu}\frac{\Sigma_{\mathbf{M}}(\lambda)}{\lambda}d\lambda}{\int_{\mu_1}^{t}\frac{\Sigma_{\mathbf{M}}(\lambda)}{\lambda}d\lambda}=1+\frac{\int_{t}^{tu}\frac{\Sigma_{\mathbf{M}}(\lambda)}{\lambda}d\lambda}{\int_{\mu_1}^{t}\frac{\Sigma_{\mathbf{M}}(\lambda)}{\lambda}d\lambda}\ge 1+\frac{\Sigma_{\mathbf{M}}(t)\log(u)}{\omega_{\mathbf{M}}(t)}.$$

Since $\omega_{\mathbf{M}}$ is slowly varying, we get
$$\forall\;u>1\;\forall\;\epsilon\in(0,1)\;\exists\;t_{u,\epsilon}>0\;\forall\;t\ge t_{u,\epsilon}:\;\;\;\frac{\omega_{\mathbf{M}}(ut)}{\omega_{\mathbf{M}}(t)}\le 1+\epsilon.$$
By combining both estimates, taking into account (again) Lemma \ref{lemma1} and finally the facts that $p\mapsto\mu_p$ is non-decreasing and $\lim_{p\rightarrow+\infty}\mu_p=+\infty$, one infers
$$\forall\;u>1\;\forall\;\epsilon\in(0,1)\;\exists\;p_{u,\epsilon}\in\NN_{\ge 2}\;\forall\;p\ge p_{u,\epsilon}:\;\;\;
\frac{p\log(u)}{\omega_{\mathbf{M}}(\mu_p)}\le\frac{\Sigma_{\mathbf{M}}(\mu_p)\log(u)}{\omega_{\mathbf{M}}(\mu_p)}=\frac{\Sigma_{\mathbf{M}}(\mu_p)\log(u)}{\log(M_0\mu_p^p/M_p)}\le\epsilon,$$
thus $\log(u^p)\le\epsilon\log(M_0\mu_p^p/M_p)$. Similarly, as in part $(i)$ this estimate is equivalent to having $u^{p/\epsilon}\frac{M_p}{M_0}=u^{p/\epsilon}\mu_1\cdots\mu_p\le\mu_p^p$ and so, finally, we arrive at
$$\forall\;u>1\;\forall\;\epsilon\in(0,1)\;\exists\;p_{u,\epsilon}\in\NN_{\ge 2}\;\forall\;p\ge p_{u,\epsilon}:\;\;\;u^{\frac{p}{(p-1)\epsilon}}(M_{p-1}/M_0)^{1/(p-1)}\le\mu_p.$$
As $u>1$ one has $u^{\frac{p}{(p-1)\epsilon}}\ge u^{1/\epsilon}$ for all $p\ge 2$ and so \eqref{lowvarassofctcharactequ1} is shown.
\qed\enddemo

We finish with the following observations ($\mathbf{M}$ being a log-convex weight sequence):

\begin{itemize}
\item[$(i)$] \eqref{lowvarassofctcharactequ1} obviously violates \eqref{rootalmostdecrlemmaequ1} and hence \hyperlink{mg}{$(\on{mg})$} fails, see the discussion in Remark \ref{rootalmostdecrlemmarem}. Indeed, the formal negation of \eqref{rootalmostdecrlemmaequ1} is $\sup_{p\in\NN_{>0}}\frac{\mu_{p+1}}{(M_p)^{1/p}}=+\infty$ and \eqref{lowvarassofctcharactequ1} is much stronger.

\item[$(ii)$] Note that any slowly varying weight function $\omega$ cannot have \hyperlink{om6}{$(\omega_6)$} since this condition implies
    $$\exists\;H_0\ge 1\;\forall\;H\ge H_0:\;\;\;\liminf_{t\rightarrow+\infty}\frac{\omega(Ht)}{\omega(t)}\ge 2.$$
    Moreover, recall that condition \hyperlink{om6}{$(\omega_6)$} for $\omega_{\mathbf{M}}$ is equivalent to \hyperlink{mg}{$(\on{mg})$}; again we refer to \cite[Lemma 2.2]{whitneyextensionweightmatrix}, \cite[Thm. 1]{matsumoto} and \cite{modgrowthstrange}. Therefore, $(ii)$ in Proposition \ref{lowvarassofctcharact} is consistent with the known information concerning \hyperlink{mg}{$(\on{mg})$}.

\item[$(iii)$] Let $\omega$ be a slowly varying weight function, then it is straightforward to see $\gamma(\omega)=+\infty$. Thus $\overline{\gamma}(\omega)=+\infty$ too and again \hyperlink{om6}{$(\omega_6)$} fails by taking into account Remark \ref{om1om6indexrem}.

\item[$(iv)$] Finally, concerning $(i)$ in Proposition \ref{lowvarassofctcharact} recall that \eqref{beta3} precisely corresponds to the condition appearing in \cite[Thm. 3.11 $(v)$]{index} with $\beta=0$ (apart from the non-effecting facts that $M_0=1$ might fail, see Remark \ref{indexweightsequrem}, and the appearing index shift when considering the sequences of quotients as described at the end of Section \ref{sequdefsection}). Thus via \cite[Thm. 3.11]{index} assumption \eqref{beta3} precisely means for the growth index from \cite{Thilliezdivision} that $\gamma(\mathbf{M})>0$, and under the additional assumption \eqref{lowvarassofctcharactequ1} even $\gamma(\omega_{\mathbf{M}})=+\infty$ follows as commented in $(iii)$. Since $\gamma(\mathbf{M})=+\infty$ is not clear, this part and the above comments $(i)$, $(ii)$ and $(iii)$ are consistent with \cite[Cor. 4.6 $(i)$ \& $(iii)$]{index}.
\end{itemize}

\bibliographystyle{plain}
\bibliography{Bibliography}

\end{document}